\documentclass[pdflatex,sn-mathphys-num]{sn-jnl}

\usepackage[normalem]{ulem}
\usepackage[pdf]{pstricks}
\usepackage[utf8]{inputenc}
\RequirePackage{ifthen}   
\RequirePackage{pgfplots}
\RequirePackage{caption}
\RequirePackage{subcaption}
\RequirePackage{graphicx}
\RequirePackage{amssymb}
\RequirePackage{amsthm}
\usetikzlibrary{matrix}
\usepackage{relsize}
\usepgfplotslibrary{groupplots}
\pgfplotsset{compat=newest}
\usepackage{bigints}
\pgfplotsset{
    every axis plot/.append style = {font = \scriptsize}
  }

\newtheorem{theorem}{Theorem}
\newtheorem{remark}{Remark}%
\newtheorem{definition}{Definition}
\numberwithin{equation}{section}
\newtheorem{lemma}[theorem]{Lemma}
\newtheorem{assumption}{Assumption}

\newcommand{\dofs}{\text{DoF}}

\newcommand{\PiO}{\Pi^{0}}

\newcommand{\RefCoordinate}{\boldsymbol{\xi}}
\newcommand{\PhysCoordinate}{\mathbf{x}}
\newcommand{\mesh}{\mathcal{T}_h}
\newcommand{\rr}{\mathbb{R}}
\newcommand{\nn}{\mathbb{N}}
\newcommand{\RefDomain}{\hat{\Omega}}
\newcommand{\Jacobian}{\mathbf{J}}
\renewcommand{\det}[1]{\text{det} \left( #1 \right)}
\newcommand{\RefGrad}{\hat{\nabla}}
\newcommand{\pp}{\mathbb{P}}

\newcommand{\norm}[1]{\left\| #1 \right\|}
\newcommand{\abs}[1]{\left| #1 \right|}
\newcommand{\F} {\mathbf{F}}

\newcommand{\JacobianDiff}{\mathcal{E}_{\Jacobian}}
\newcommand{\DetDiff}{\mathcal{E}_{\text{det}}}
\newcommand{\InvDiff}{\mathcal{E}_{\text{inv}}}
\newcommand{\CofDiff}{\mathcal{E}_{\text{cof}}}

\newcommand{\RefMesh}{\hat{\mathcal{T}}_h}
\newcommand{\VirtualMesh}{\mathcal{T}_h}
\newcommand{\ExtensionOp}{\mathfrak{E}}
\newcommand{\Exten}[1]{\tilde{#1}}
\newcommand{\pertB}{\tilde{\mathcal{A}}}
\newcommand{\Ah}{\mathcal{A}_h}
\newcommand{\projE}{E_h^\pi}
\newcommand{\holder}{H\"{o}lder}

\graphicspath{{Fig/}}

\begin{document}


\title[Isoparametric Virtual Element Methods]{Isoparametric Virtual Element Methods}

\author[1]{\fnm{Andrea} \sur{Cangiani}}

\author[2]{\fnm{Andreas} \sur{Dedner}}

\author[3]{\fnm{Matthew} \sur{Hubbard}}

\author*[3]{\fnm{Harry} \sur{Wells}}\email{h.wells.research@gmail.com}

\affil[1]{\orgdiv{Mathematics Area}, \orgname{International School for Advanced Studies (SISSA)}, \orgaddress{\street{Via Bonomea 265}, \postcode{34136}, \state{Trieste}, \country{Italy}}}

\affil[2]{\orgdiv{Department of Mathematics}, \orgname{University of Warwick}, \orgaddress{ \postcode{CV4 7AL}, \state{Coventry}, \country{United Kingdom}}}

\affil[3]{\orgdiv{School of Mathematical Sciences}, \orgname{University of Nottingham}, \orgaddress{\street{University Park}, \postcode{NG7 2RD}, \state{Nottingham}, \country{United Kingdom}}}





\abstract{We present two approaches to constructing isoparametric Virtual Element Methods of arbitrary order for linear elliptic partial differential equations on general two-dimensional domains. The first method approximates the variational problem transformed onto a computational reference domain. The second method computes a \emph{virtual} domain and uses bespoke polynomial approximation operators to construct a computable method. Both methods are shown to converge optimally, a behaviour confirmed in practice for the solution of problems posed on curved domains.
}

\maketitle

\section{Introduction}

In recent years there has been an increasing interest in exploiting the mesh flexibility allowed by the Virtual Element Method (VEM)~\cite{basicprinciples} for the efficient discretisation of problems posed on general transformed and time-dependent domains \cite{Lipnikov2019AMeshes,Lipnikov2020ConservativeMeshes,Mazzia2020VirtualEnvironments,Wells2024AMethod}, including a broader interest on combining the VEM with other numerical methods
(e.g. discontinuous Galerkin) \cite{Antonietti2019NumericalGrids,Gaburro2020HighChanges,Gaburro2021AChange,Huang2020AnisotropicMeshes}.
While these studies have already demonstrated the effective utilisation of mesh flexibility permitted by the VEM through compelling numerical experiments, a rigorous theoretical foundation to support the VEM's applicability in this context is still lacking.

Motivated by the need to fill such a theoretical gap, we analyse the effect of domain discretisation errors  in the solution of {\em steady} problems posed on general domains, as a starting point to develop the required theory to analyse the aforementioned approaches on time-dependent domains. In particular, we propose and analyse \emph{Isoparametric Virtual Element Methods} (IsoVEMs) for the solution of linear elliptic problems posed on general curved domains in two dimensions given as the image of a smooth transformation of a \emph{reference polygonal domain}.

Under the classical definition \cite{CiarletElliptic}, isoparametric FEMs are characterised by the use of the same finite element space for the solution as that used to approximate the domain geometry. 
Typically, the geometric elements are constructed from a single reference element via finite element transformations of a given degree. Then, in order to obtain balanced domain and solution approximations, the problem is discretised using a finite element space of the same degree. 
In the case of the VEM, we have no obvious concept of a reference element as virtual element discretisations are constructed directly on the \emph{physical} domain to allow for heterogeneous elements. Thus, we consider, instead, the concept of a \emph{reference polygonal domain} from which a computational domain and mesh approximating the physical domain can be constructed. For instance, the reference polygonal domain may simply be obtained by a fixed linear interpolation of the physical domain's boundary. We note that this idea naturally extends to time-dependent problems whereby the reference polygonal domain may be constructed on the reference (initial) configuration.
In principle, once the reference polygonal domain is given, the isoparametric FEM's idea extends to the VEM as we could use a standard virtual element space to define a continuous transformation of the reference domain on which a same order VEM may be used to solve the problem given on the transformed domain. However, the virtual nature of this transformation makes this naive approach impractical in general, as  the evaluation of virtual element functions on the resulting isoparametric elements may not be available.
Instead, we introduce two approaches to construct a computable isoparametric VEM, herein referred to as \emph{reference} and \emph{physical} IsoVEM. These methods are designed in such a way that they emulate the classical isoparametric finite element method \cite{Ciarlet1972InterpolationMethods,Lenoir1986OptimalBoundaries} extended to polygonal meshes. 

To motivate our consideration of two isoparametric methods, we emphasize the role of PDE data in the context of a given problem. Often, the PDE data may be constrained to a specific computational reference frame, as exemplified in Arbitrary Lagrangian-Eulerian schemes for computational fluid dynamics (CFD) applications \cite{Donea2004,Richter2017}. In other scenarios, the PDE data might need to be computed or approximated within the physical reference frame and then extended to the PDE domain. These varying requirements of PDE data handling motivates the consideration of separate IsoVEMs, accommodating the framework of different computational and physical contexts.

In the reference IsoVEM, we consider the problem mapped back to the reference polygonal domain by a discrete VEM mapping. A VEM is then used to approximate the terms resulting from this transformation such as the Jacobian operator, its determinant, and inverse Jacobian operator. Using these approximations results in a computable VEM on the reference polygonal domain. This framework closely aligns with the approach taken in \cite{Lipnikov2019AMeshes,Lipnikov2020ConservativeMeshes} in which Arbitrary Lagrangian-Eulerian schemes are constructed on a computational reference domain using approximations of the transformed weak formulations. We note, however, that in~\cite{Lipnikov2019AMeshes,Lipnikov2020ConservativeMeshes} the solution itself is computed by a discontinuous Galerkin method and no analysis is presented to support the proposed method.

For the physical IsoVEM, we consider the problem on the domain generated by the VEM approximation of the true reference-to-physical transformation. As such, the computational mesh is made of curved polygons with polynomial edges. Here we take advantage of the curved virtual element space presented in \cite{daVeigaCurvedVEM} for the solution of the constant coefficients Poisson problem and generalised here to a broader class of linear elliptic problems. To ensure the method is computatable, we introduce bespoke polynomial approximation operators such that the consistency and stability terms can be easily computed to a sufficient degree of accuracy. 
Note that we could have instead considered directly the curved VEM of \cite{daVeigaCurvedVEM} defined on the true physical domain, hence avoiding the introduction of a geometrical error.
Our, alternative, approach is  motivated by computability considerations. 
While in curved VEM the mapping used to generate the curved elements is generally unknown and integration over curved elements must be approximated, the polynomial approximation employed in our isoparametric approach is available and uses the minimum amount of quadrature point evaluations whilst still maintaining the same degree of accuracy. Moreover, our approach becomes the only one available for three-dimensional problems, as in this case mesh faces are themselves virtual.
Our approach is also relevant to time-dependent problems, when only the discrete domain is available as it is itself obtained as the numerical solution of a problem determining the evolution of the domain such as in fluid-structure interaction.
On the other hand, our approach requires a more involved analysis, accounting for the domain approximation error.
 
To support the proposed methods, we present $H^1$-norm \textit{a priori} error estimates which show that the same degree of accuracy can be achieved as the classical isoparametric finite element method.
For brevity, we limit the scope of the analysis to two-dimensional problems but emphasise that the analysis can be extended to three-dimensional problems with additional work. 

In addition, we implement both methods using the DUNE-VEM module \cite{Dedner2022ASpaces} which is part of the Distributed and Unified Numerics Environment \cite{Bastian2021TheDevelopments,Dedner2012Dune-Fem:Computing}. Numerical experiments within DUNE validate the $H^1$-norm estimates and obtain optimal rates of convergence in the $H^1$-norm and $L^2$-norm error.

Although the motivation of this paper is the transformation of two-dimensional domains, the framework relates to the field of surface PDEs. Indeed there has been recent publications in developing VEMs for surface problems such as the Laplace-Beltrami equation and bulk-surface PDEs \cite{Frittelli2021Bulk-surfaceDimensions,Frittelli2018VirtualSurfaces}, elliptic \cite{Bachini2021Arbitrary-orderSurfaces} and diffusion-reaction \cite{Li2022LocalSurface} equations on surfaces. In particular, the methods presented in \cite{Bachini2021Arbitrary-orderSurfaces} and \cite{Li2022LocalSurface} closely relate to the proposed methods in this paper. The former relating to our first proposed discretization and the latter being equivalent to a two-dimensional, linear case of our second method. This paper can be viewed as an extension on these ideas, providing a foundation to develop higher-order methods and derive error estimates.

The layout of this paper is as follows. In Section \ref{sec::EllipticProblem} we define the PDE and domain transformation considered in the analysis and present the weak formulation of this problem. In Section \ref{sec::VEM} we recall the fundamentals of the VEM and present the curved VEM from \cite{daVeigaCurvedVEM} for general elliptic problems. We then formulate  the reference and physical IsoVEMs in Section~\ref{sec::MethodI} and detail their a priori analysis in Section \ref{sec::MethodII}. Numerical results are presented in Section \ref{sec::NumericalTests} for a set of transformations and PDEs. Concluding remarks are made in Section \ref{sec::Conclusion}.

\section{Preliminaries}\label{sec::EllipticProblem}
Let $\omega \subset \rr^n$, $n =1,2$, denote an open bounded subset of  $\rr^n$.  For $m$ a non-negative integer and $1\le p \le \infty$, we denote by $W_m^p(\omega)$ the respective Sobolev space   endowed with norm $\norm{\cdot}_{m,p,\omega}$ and semi-norm $\abs{\cdot}_{m,p,\omega}$. In the case of $p=2$ we denote the corresponding Hilbert space by $H^{m}(\omega)$ with norm and semi-norm denoted by $\norm{\cdot}_{m,\omega}$ and $\abs{\cdot}_{m,\omega}$, respectively. Finally, with $\mathbb{P}_k(\omega)$ we denote the space of polynomials of degree $k\in\nn$ over $\omega$ and with $\mathbb{P}_k$ the polynomial space over $\rr^2$.

We use the relation $x \lesssim y$ to indicate the existence of an hidden constant $C>0$ such that $x \leq Cy$. Unless stated explicitly otherwise, it is assumed that this hidden constant is independent of the solution and the discretisation parameter $h$.

\subsection{Model problem}

On $\Omega \subset \rr^2$ a bounded Lipschitz domain, we consider the boundary value problem
\begin{align}
    -\nabla \cdot (a \nabla u) + \mathbf{b} \cdot \nabla u + c u &= f \ \ \  \text{in } \Omega,\label{eq::prob1} 
\\
    u &= 0 \ \  \text{on } \partial \Omega,\label{eq::prob2} 
\end{align}
with given coefficients  $a \in L^{\infty}(\Omega)$, $\mathbf{b} \in \left[ L^{\infty}(\Omega) \right]^2$, $c \in L^\infty(\Omega)$, and forcing function $f\in L^2(\Omega)$. We assume the existence of two constants $a_0>0$ and $\mu_0 \geq 0$ and $\mu:=c - \frac{1}{2} \nabla \cdot \mathbf{b} \in L^{\infty}(\Omega)$ such that
\[
    a_0\le a \quad \text{and}\quad
 \mu_0\le     \mu \quad \text{a.e. in } \Omega. 
\]
From this, we consider the weak formulation with convective term in skew-symmetric form: find $u \in H^1_0(\Omega)$ such that
\begin{equation}\label{eq::EllipticModelPDE}
    \mathcal{A}(u,v) := A(u,v) + B(u,v) + C(u,v) = l(v) \ \ \ \forall v \in H^1_0(\Omega),
\end{equation}
where, for all $u,v \in H^1_0(\Omega)$,
\begin{align}\label{eq::modelAB}
     A(u,v) &:= \int_{\Omega} a \nabla u \cdot \nabla v\ d\PhysCoordinate, \quad
     B(u,v) := \frac{1}{2} \int_{\Omega} \mathbf{b} \cdot (v \nabla u - u \nabla v)\ d\PhysCoordinate,\\
     \label{eq::modelCl}
     C(u,v) &:= \int_{\Omega} \mu\  u v\ d\PhysCoordinate, \qquad\,\,\,\,
     l(v) := \int_{\Omega}   fv\ d\PhysCoordinate.
\end{align}
It is well known that this weak formulation has a unique solution in $H^1_0(\Omega)$ by the Lax-Milgram Theorem.

\begin{remark}
It is by all means possible to slightly relax the above assumptions on the coefficients and to consider the weak form obtained without writing the  convective term in skew-symmetric form by following e.g.~\cite{ellipticVEM}. We refrain from doing so here to avoid inessential technicalities.
\end{remark}

\subsection{Reference domain}

We assume we are given a polygonal reference domain $\RefDomain \subset \rr^2$ such that there exists a bijective  reference-to-physical transformation $\F : \RefDomain \rightarrow \Omega$.
We further suppose that $\F$ is bi-Lipschitz, in line with assumptions made within ALE analysis \cite{Bonito2013Time-DiscreteStability,Bonito2013Time-discreteAnalysis,Formaggia2004StabilityALEFEM,Gastaldi2001AElements} and that $ \F \in \left[W^{m+1}_\infty(\RefDomain) \right]^2$ for some integer $m \in \mathbb{N}$ to be specified later. The bi-Lipschitz condition allows for the analysis to be conducted independently of regularity conditions on $\F^{-1}$ and its corresponding VEM approximation. This contrasts with the classical isoparametric FEM work in which an inverse mapping between isoparametric finite elements is easily computable \cite{Ciarlet1972InterpolationMethods,Lenoir1986OptimalBoundaries}.

As a consequence, the Jacobian matrix
\[
    \Jacobian_\F := \nabla \F,
\]
 is invertible~\cite{Bonito2013Time-DiscreteStability,Bonito2013Time-discreteAnalysis}. Moreover,  the determinant of $\Jacobian_\F$  is uniformly bounded and strictly positive; namely,  there exists a $j_0 \in (0,1)$ such that
\[
 j_0 \le j := \det{\Jacobian_\F}  \in L^{\infty}(\RefDomain).
\]

The domain $\RefDomain$  is the basis of the IsoVEM approach. The polygonal meshes defined over $\RefDomain$ will serve as collections of polygonal reference elements. Corresponding virtual maps constructed over the reference elements will provide the computational meshes approximating  $\Omega$.
\begin{remark} 
For instance, $\RefDomain$ could just be defined as a polygon whose boundary interpolates $\partial \Omega$. In this case, the IsoVEM introduced below provides a direct generalisation of the classical isoparametric finite elements on general meshes. 
\end{remark}

To account for the discrepancy between the problem domain and the computational domain, we require the availability of extension operators, which is guaranteed by the well-known Stein's Extension Theorem \cite{Stein1971SingularFunctions,Adams2003}.
\begin{theorem}[The Stein Extension Theorem]\label{theorem::SteinExtension}
Let $\omega \subset \rr^2$ be bounded and Lipschitz. For all $p \in [1,\infty]$ and integers $s \geq 0$ there exists an extension operator $\ExtensionOp : W_p^s(\omega) \rightarrow W_p^s(\rr^2)$ such that, for all $v \in W^s_p(\omega)$,
$\ExtensionOp v = v$ a.e. in  $\omega$ and $\norm{\ExtensionOp v}_{s,p,\rr^2}\lesssim \norm{v}_{s,p,\omega}$,
where the hidden constant depends only on $s$ and $p$.
\end{theorem}
In this paper, extensions will be required for the PDE data and exact solutions; for short we shall denote these by $\Exten{v}:= \ExtensionOp v$ and so on.

\section{$C^0$-conforming VEM}\label{sec::VEM}
In this section we present an enhanced version of the $C^0$-conforming VEM on two-dimensional polygonal meshes with curved edges from~\cite{daVeigaCurvedVEM}. We present these spaces in a general form and independent of the IsoVEM frames of reference for ease of reading, noting that the standard VEM spaces on polygonal meshes can be considered as a special case of this presentation \cite{basicprinciples, Ahmad2013EquivalentMethods, DaVeiga2014}. The implementation details for VEMs are well documented in the literature~\cite{DaVeiga2014,Sutton2017TheMATLAB,Dedner2022ASpaces} and thus are omitted. 

A mesh $\mesh$ is a partitioning of a general domain $\Omega$ into simple, non-overlapping polygons 
 $E \in \mesh$ with curved edges of class $C^{m+1}$. The mesh $\mesh$ is characterised by the global mesh size $ h := \max_{E \in \mesh} h_{E}$ with $h_{E}$ denoting the diameter of $E$.  We further denote by $e$ a generic (possibly curved)  edge belonging to $\partial E$.

For $k\ge 0$ integer we denote by $\pp_k(E)$, $E\in\mesh$, the space of polynomials of degree up to $k$ and by $\pp_k(\mesh)$ the corresponding space of piecewise polynomials with respect to $\mesh$.

The mesh is required to satisfy the following shape regularity assumptions which are standard in the VEM literature~\cite{basicprinciples,DaVeiga2014,Ahmad2013EquivalentMethods}. 

\begin{assumption}[Shape regularity]\label{assumption::ShapeRegular} 
Every $E \in \mesh$ is a star-shaped domain or a finite union of star shaped domains with respect to a ball of radius greater than $\alpha h_{E}$ for some uniform $\alpha > 0$. Additionally,   the length of any edge $e \in \partial E$ is greater than $\alpha h_{E}$.
\end{assumption}

The above assumption readily ensures the availability of polynomial best approximation results, see e.g.~Theorem 1.45 in \cite{HHOBook}.
\begin{theorem}[$L^2$-projection accuracy]\label{theorem::Pi0Accuracy}
Let a polynomial degree $k \geq 0$, an integer $s \in \{ 0,...,k+1\}$, and a real number $p \in [1,\infty]$ be given. Then, the piecewise \emph{$L^2$-projection operator } $ \PiO_k:L^1(\Omega)\rightarrow \pp_k(\mesh)$
satisfies, on each element $E \in \mesh$,
\begin{equation*}
    \abs{v - \PiO_k v}_{l,p,E} \lesssim h_{E}^{s-l} \abs{v}_{s,p,E},
\end{equation*}
for all $v\in W_p^s(E)$ and $l \in \{0,...,s\}$.
\end{theorem}

\subsection{Enhanced curved VEM spaces}
Let $k\ge 1$ be fixed. We introduce the enhanced version of the curved-edged polygonal mesh VEM of order $k$ from~\cite{daVeigaCurvedVEM}.

Given $E \in\mesh$, we define for each edge $e \subset \partial E$ its corresponding reference interval $I_e = [0,l_e]$ where $l_e = \abs{e}$. We denote the (invertible) arc length parameterisation of the edge $e$ by $\gamma_{e}:I_e \rightarrow e$ and introduce over $e$ the mapped polynomial space of degree $k\ge 1$ as
\begin{equation*}
    \pp^{\gamma}_k(e) = \{ q \in L^2(e)\ :\ q = \hat{q} \circ \gamma_{e}^{-1},\ \hat{q} \in \pp_k(I_e) \}.
\end{equation*}
Note that if $e$ is straight then $\gamma_{e}$ is affine and it follows that $\pp^{\gamma}_k(e)\equiv \pp_{k}(e)$. Otherwise, $\pp^{\gamma}_k(e)$ is not made of polynomials in general.
In line with the construction of polygonal enhanced virtual element spaces~\cite{Ahmad2013EquivalentMethods}, we first introduce the enlarged local space 
\begin{equation*}
    W_k(E) := \{ v_h \in H^1(E) \cap C^0(\overline{E})\ :\ \Delta v_h \in \pp_{k}(E),\ v_h|_{e} \in \pp^\gamma_k(e)\ \forall e \subset \partial E \}.
\end{equation*}
\begin{remark}
    We note that the parametrisation of each curved edge by $I=[0,l_e]$ is a choice we make here to simplify the presentation. Alternatives include parametrising on the unit interval $[0,1]$ or a reference edge $\hat{e}$ from a computational reference polygon $\hat{E}$ as in the original presentation~\cite{daVeigaCurvedVEM}.
\end{remark}
We define a set of linear forms which will be used as degrees of freedom for the curved VEM space defined below in~\eqref{eq::CurvedVEMSpace}.
\begin{definition}[Virtual element $\dofs$]
\label{definition::dofs}
For a regular enough function $v$ of $E$, we introduce the following degrees of freedom:
\begin{itemize}
    \item the point value of $v$ at all vertices of $E$;
    \item for all $e\subset \partial E$, the edge moments 
    $\displaystyle\frac{1}{l_e}\int_{I_e} (v_h|_{e} \circ \gamma_{e}) q\ dx$, for all $q\in \mathcal{M}_{k-2}(I_e)$;
    \item the internal moments $\displaystyle\frac{1}{\abs{E}}\int_{E} v q\ dx$ for all $q\in \mathcal{M}_{k-2}(E)$.
\end{itemize}
 Here, $\mathcal{M}_m(\omega)$, for $\omega=I_e$ or $E$ and $m$ an integer, represents a set of \emph{scaled} monomials which constitute a basis for  $\mathbb{P}_m(\omega)$ with the convention that $\mathcal{M}_m(\omega)=\emptyset$ if $m<0$. The scaled monomials have the important property that $\norm{q}_{0,\infty,E} \sim 1$ for all $q \in \mathcal{M}_m(\omega)$ and will be frequently used in our analysis; see e.g. \cite{basicprinciples,Ahmad2013EquivalentMethods,DaVeiga2014} for more details.
\end{definition}

The forms of Definition~\ref{definition::dofs}  constitute also an (incomplete) set of $\dofs$ for the space $W_k(E)$.
We let  $\Pi^*_k: W_k(E)\rightarrow \mathbb{P}_l(E)$ denote any $L^2$-stable \emph{auxiliary} projection operator that is computable by only accessing such $\dofs$.
The existence of auxiliary projection operators is a consequence of the fact that the above set of $\dofs$ is sufficient to determine the elements of $\mathbb{P}_k(E)$. 

\begin{remark}\label{remark:proj}
 In~\cite{Ahmad2013EquivalentMethods} it is shown that the $H^1$-seminorm  projection known in the VEM literature as  $\Pi^\nabla$ can serve as an auxiliary projection operator. This approach is extended in~\cite{daVeigaCurvedVEM} to the curved polygonal elements  discussed below, where $\Pi^*$ is chosen as an approximation of $\Pi^\nabla$ by using a sufficiently accurate quadrature scheme. However, this choice is not essential. For instance, a complete family of auxiliary projection operators is showcased in~\cite{ConNonConVEM} for the straight VEM case. As mentioned in the main text,  in our implementation we make use of the least square fitting technique from~\cite{Dedner2022RobustCoefficients}. 
\end{remark}

We are now ready to introduce the enhanced local virtual element space as
\begin{align}
    V_k(E) := \Big\{ v_h \in W_k(E) \ :\ &\int_{E} (\Pi^*_k v_h -v_h)q\ d\PhysCoordinate = 0\ \forall q \in \pp_{k}(E)\backslash\pp_{k-2}(E) \Big\}.
    \label{eq::CurvedVEMSpace}
\end{align}

Following \cite{daVeigaCurvedVEM,Ahmad2013EquivalentMethods} it can be shown that the $\dofs$ of Definition~\ref{definition::dofs}  represent a set of unisolvent degrees of freedom for $V_k(E)$. 
 Moreover, such $\dofs$ allow for the evaluation of 
$
\Pi^0_{k}  v_h$ and $ \Pi^1_{k-1}  v_h := \Pi^0_{k-1} \nabla  v_h$ 
 as well as of $\Pi^*_k v_h$, which is computable by construction for all   $v_h\in V_k(E)$.

A global $C^0$-conforming virtual element space is then defined by gluing together the local spaces in~\eqref{eq::CurvedVEMSpace}; thus, we introduce the global (curved) VEM space
\begin{equation}\label{eq::Method2GlobalVEMSpace}
    V_h^k := \{ v_h \in H^1(\Omega)\ :\ v_h|_{E} \in V_{k}(E)\ \ \ \forall E \in \VirtualMesh \},
\end{equation}
together with the space $V_{h,0}^k:= V_h^k \cap H^1_0(\Omega)$ of virtual functions with zero trace on $\partial \Omega$.

It is easy to show that if $E \in \VirtualMesh$ is straight-edged, the  polynomial space  $\pp_{k}(E)$ is a subspace of $V_k(E)$. However, this is not the case if curved edges are present. Still, optimal approximation results can be shown in all cases. The following best approximation result 
can be easily derived by combining 
Theorem 3.7 from \cite{daVeigaCurvedVEM}, dealing with the original curved VEM space, and Theorem 11 in~\cite{ConNonConVEM}, which proves similar results for the enhanced VEM space with straight edges, 
cf. also~\cite{BeiraodaVeiga2020PolynomialEdges}.
\begin{theorem}[VEM interpolation error estimate]\label{theorem::GlobalCurvedVEMInterp}
Let $\mesh$ satisfying Assumption~\ref{assumption::ShapeRegular} and $k \geq 1$  be given. For all $v \in H^s(\Omega)\cap H^1_0(\Omega)$, $s \in \{2,...,k+1\}$,
there exists a $v_I \in V_{h,0}^k$ such that,
for any element $E \in \mesh$, it holds that
\begin{equation*}
    \norm{v-v_I}_{0,E} + h_{E} \abs{v - v_I}_{1,E} \lesssim h_{E}^{s} \abs{v}_{s,E},
\end{equation*}
where the hidden constant depends on the shape-regularity of $\mesh$ and on $k$.
\end{theorem}

\begin{remark}\label{remark:curved_vs_straight}
We stress once again that, whenever the mesh $\mesh$ is made of polygonal elements, the virtual element space~\eqref{eq::Method2GlobalVEMSpace} reduces to the standard $C^0$-conforming space of \cite{Ahmad2013EquivalentMethods}. Further, in the definition of the IsoVEM below, physical elements will be made at most of edges given as polynomial curves, thus simplifying implementation considerably, as we shall see. Nevertheless, here we have introduced the overly general curved VEM to exploit  relevant estimates which are readily available from~\cite{daVeigaCurvedVEM}.
\end{remark}

\subsection{Curved VEM}
The original presentation of the curved VEM in~\cite{daVeigaCurvedVEM} was limited to the Poisson problem with constant coefficients. For completeness, we report here its generalisation to the  elliptic problem~\eqref{eq::EllipticModelPDE} which appears to be new. The results shown are easily obtained by combining~\cite{daVeigaCurvedVEM} with the framework in~\cite{ConNonConVEM}.

 On each ${E}\in\mesh$, we define the local virtual element forms by
\begin{align*}
    A_h^{{E}}({u},{v}) &:= \int_{{E}} {a}\  \Pi^1_{k-1} {u}_h \cdot \Pi^1_{k-1} {v}_h\  d\PhysCoordinate,
    \quad\qquad\qquad
    C_h^{{E}}({u},{v}) := \int_{{E}} {\mu}\  \PiO_{k} {u}_h\ \PiO_{k} {v}_h\  d\PhysCoordinate,
    \\
    B_h^{{E}}({u},{v}) &:= \frac{1}{2} \int_{{E}} {\mathbf{b}} \cdot (\PiO_{k} {v}_h\  \Pi^1_{k-1} {u}_h - \PiO_{k} {u}_h\  \Pi^1_{k-1} {v}_h)\  d\PhysCoordinate,
    \quad
    l_h^{{E}}({v}) := \int_{{E}} {f}\ \PiO_{k}\ {v}_h\  d\PhysCoordinate,
\end{align*}
for all ${u}, {v}\in V_k({E})$. 
Further, as customary in the definition of VEMs, we also require a VEM stabilisation  form. That is, any bilinear form   $S^{{E}}(\cdot,\cdot)$ satisfying the following.
\begin{assumption}\label{assumption::Stab::Stab}
For each ${E}\in\mesh$, a VEM stabilisation form $S^{{E}}:V_k({E})\times V_k({E})\rightarrow \rr$ is available which is computable from the degrees of freedom of Definition~\ref{definition::dofs} and is such that
    \begin{align*}
        \mathcal{A}^{{E}}({v},{v}) \lesssim S^{{E}}({v} , {v})  \lesssim \mathcal{A}^{{E}}({v},{v}) \qquad \forall  {v} \in V_k({E})\, \text{ such that } \, \Pi^0_{k}  {v}=0.
    \end{align*}
Here, $\mathcal{A}^{{E}}$ denotes the restriction to ${E}$ of the bilinear form $\mathcal{A}$ from~\eqref{eq::EllipticModelPDE}.
\end{assumption}
We then define the discrete analogue of the global  forms in ~\eqref{eq::EllipticModelPDE} as follows.  For ${u}, {v}\in  {V}^k_{h}$, we let
\begin{align}
    \mathcal{A}_h({u},{v}) &:= \sum_{{E} \in \mesh}\Big[ A_h^{{E}}({u},{v}) + B_h^{{E}}({u},{v}) + C_h^{{E}}({u},{v}) \nonumber \\
    &\hspace{2cm}+ S^{{E}}({u}_h - \Pi^0_k {u}_h, {v}_h - \Pi^0_k {v}_h)\Big], 
    \label{eq::bilinear_curved}\\
    l_h({v}) &:= \sum_{{E} \in \mesh} l_h^{{E}}({v}).
    \label{eq::linear_curved}
\end{align}
The curved VEM now reads: find ${u}_h \in {V}^k_{h,0}$ such that
\begin{align}\label{eq:curvedVEM}
    \mathcal{A}_h({u}_h,{v}_h) =l_h({v}_h) \qquad \forall {v}_h \in {V}^k_{h,0}.
\end{align}

The implementation of the above curved VEM requires the ability to compute integrals of polynomials over curved domains. We refer to~\cite{daVeigaCurvedVEM} for more details.

\begin{theorem}[Strang-type bound]\label{theorem::StrangBound1}
    Under Assumption~\ref{assumption::Stab::Stab}  the curved VEM~\eqref{eq:curvedVEM} is well-posed. Moreover,  letting ${u}\in H^1_0(\Omega)$ be the solution of~\eqref{eq::EllipticModelPDE} and ${u}_h \in {V}_{h,0}^k$ the curved VEM solution  of~\eqref{eq:curvedVEM}, we have
    \begin{align}
        \norm{{u} - {u}_h}_{1,\Omega}
    &\lesssim  \inf_{{v}_h \in {V}_{h,0}^k} \norm{{u}-{v}_h}_{1,\Omega} + \inf_{{p} \in \mathbb{P}_k(\mesh)} \norm{{u}-{p}}_{1,\Omega} \nonumber \\
    &\hspace{-.4cm}+ \sup_{{w}_h \in {V}_{h,0}^k \backslash \{ 0 \}} \frac{\abs{ l_h({w}_h) - l({w}_h) }}{\norm{{w}_h}_{1,\Omega}}\nonumber\\
    &\hspace{-.4cm}+ \inf_{{p} \in \mathbb{P}_k(\mesh)}\sup_{{w}_h \in {V}_{h,0}^k \backslash \{ 0 \}}\displaystyle \frac{\sum_{{E} \in \mesh} \abs{\mathcal{A}^{{E}}_h({p},{w}_h) - \mathcal{A}^{{E}}({p},{w}_h)}}{\norm{{w}_h}_{1,{\Omega}}}.
    \label{eq::StrangBound}
    \end{align}
    The hidden constant only depends on the continuity and coercivity of the bilinear form~\eqref{eq::bilinear_curved}.
\end{theorem}

The following error bound for the curved VEM is easily proven  following the steps detailed in~\cite{ConNonConVEM} to bound each term in the Strang-type bound~\eqref{eq::StrangBound} exploiting also 
 Theorem~\ref{theorem::Pi0Accuracy} and Theorem~\ref{theorem::GlobalCurvedVEMInterp}.

\begin{theorem}[Curved VEM error estimate]\label{theorem::CurvedH1Error} 
    Let  Assumptions~\ref{assumption::ShapeRegular},~\ref{assumption::Stab::Stab} hold true. Suppose also that, for some $s\in\{1,\dots, k\}$, the solution ${u}$ of~\eqref{eq::EllipticModelPDE} satisfies ${u} \in H^{s+1}(\Omega)$ with data  ${a},{\mathbf{b}}, {\mu} \in W^{s+1}_{\infty}(\RefDomain)$ and $f\in H^{s-1}(\Omega)$.   Then, for $k\le m$, the curved VEM solution  of~\eqref{eq:curvedVEM} satisfies 
    \begin{equation}
        \norm{{u} - {u}_h}_{1,\Omega}+ \norm{\nabla {u} - \Pi^1_{k-1} {u}_h}_{0,\Omega} \lesssim h^s \left( \norm{{u}}_{s+1,\Omega} + \norm{{f}}_{s-1,\Omega} \right).
    \end{equation}
    The hidden constant depends on the shape regularity of $\mesh$, on $k$, on $a_0$ and $\mu_0$, and the $W^{s+1}_\infty (\Omega)$-norms of ${a}$, ${\mathbf{b}}$, and ${\mu}$.
    \end{theorem}

\section{Virtual element mapping}\label{sec:mapping}

The construction of the IsoVEMs depends on a (equal order) virtual element approximation of the reference-to-physical transformation $\F$ rather than on  the transformation itself, the knowledge of which will not be required in practice.
The idea of using the VEM to construct conforming approximations of the computational domain was already used in \cite{Lipnikov2019AMeshes,Lipnikov2020ConservativeMeshes}.

Firstly, we require a virtual element space on the polygonal reference domain $\hat{\Omega}$. Let $\RefMesh$ be a polygonal mesh of $\RefDomain$, containing polygonal elements $\hat{E} \in \RefMesh$ satisfying Assumption \ref{assumption::ShapeRegular}, and $\hat{V}_{h}^l$ the corresponding enhanced virtual element space of degree $l\ge 1$. By design, $\hat{V}_h^l$ is the standard \emph{enhanced} virtual element space on straight-edged meshes \cite{basicprinciples,Ahmad2013EquivalentMethods,DaVeiga2014}, cf. Remark~\ref{remark:curved_vs_straight} above.

A \emph{virtual computational domain} approximating the physical domain $\Omega$ may be defined as the image of the reference domain $\RefDomain$ through a vectorial virtual element approximation of $\F$ with components in the virtual element space  $ \hat{V}_{h}^l$. That is, we let  $ \F_h\in  [\hat{V}_{h}^l]^2$ and  for each $\hat{E}\in\RefMesh$ we define the corresponding physical element $E_h = \F_h(\hat{E})$. We denote the virtual domain and mesh so obtained by
\begin{equation}\label{eq:omega_h}
\Omega_h:=\F_h(\RefDomain),\quad \mesh:=\{E_h:=\F_h(\hat{E})\}.
\end{equation}
The analysis below requires that $\Omega_h$ is a sufficiently accurate approximation of the physical domain $\Omega$.

\begin{figure}[htbp]
\centering
\begin{minipage}[t]{.2\textwidth}
  \centering
  \includegraphics[width=.9\linewidth]{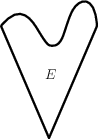}
\end{minipage}
~
\begin{minipage}[t]{.2\textwidth}
  \centering
  \includegraphics[width=.9\linewidth]{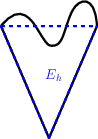}
\end{minipage}
~
\begin{minipage}[t]{.2\textwidth}
  \centering
  \includegraphics[width=0.9\linewidth]{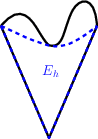}
\end{minipage}
~
\begin{minipage}[t]{.2\textwidth}
  \centering
  \includegraphics[width=0.9\linewidth]{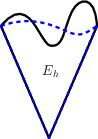}
\end{minipage}
\captionof{figure}{From left to right, a physical element $E$ and its virtual interpolating polygons $E_h$ obtained with virtual map of order $l=1$, $l=2$, and $l=3$.}\label{fig::mappedEDemo}
\end{figure}

\begin{assumption}[Virtual element mapping accuracy]\label{assumption::F_hAccuracy}
Let $\RefMesh$ be a polygonal mesh of $\RefDomain$ and let $m,l \in \mathbb{N}$ with $l\le m$. The reference-to-physical transformation $\F \in \left[ W^{\infty}_{m+1}(\RefDomain) \right]^2$ has a VEM approximation $\F_h \in \left[ \hat{V}_h^l \right]^2$ of degree $l\ge 1$ such that for all $s \in \{0,1,...,l \}$
    \begin{align*}
        \norm{\F - \F_h}_{0,\RefDomain} + h \abs{\F - \F_h}_{1,\RefDomain} \lesssim h^{s+1} \abs{\F}_{s+1,\RefDomain},
    \end{align*}
    where the hidden constant depends only on the mesh regularity of $\RefMesh$ and on $l$.
\end{assumption}
For instance, if $\F_h$ is a VEM interpolant of $\F$, then the approximation properties required by Assumption~\ref{assumption::F_hAccuracy} are guaranteed by Theorem~\ref{theorem::GlobalCurvedVEMInterp}. This case is depicted in Figure \ref{fig::mappedEDemo}, showing a physical element $E=\F(\hat{E})$ and its VEM interpolant $E_h = \F_h(\hat{E})$ of degree $l=1,2,3$. However, we prefer to keep the above as an assumption to cater for other cases, such as $\F_h$ being a sufficiently accurate VEM solution to a separate problem.

As a consequence of the bi-Lipschitz continuity condition on $\F$, it can be easily shown that a reference element $\hat{E}$ is shape regular if and only if its corresponding physical element $E = \F(\hat{E})$ is shape regular under the conditions of Assumption \ref{assumption::ShapeRegular}. The ratio between the regularity parameters is bounded above and below dependent on the bi-Lipschitz constants. It can be shown using Assumption \ref{assumption::F_hAccuracy} and a continuity argument that $\mesh$ is a virtual domain of shape regular elements; see \cite{Weier2019AnisotropicAnalysis,Antonietti2019TheDiscretizations} for a framework to estimate perturbed regularity parameters on transformed elements. 

\section{Formulation of the isoparametric VEMs}\label{sec::formulation}

Both methods presented below are based on  the virtual element map of Section~\ref{sec:mapping}. 
Letting $\RefMesh$ be a mesh of $\RefDomain$ 
satisfying the shape-regularity Assumption~\ref{assumption::ShapeRegular}, 
we consider the corresponding reference virtual element space $ \hat{V}_{h}^l$ for some discrete map order $l \in \mathbb{N}$ and let $\F_h\in [\hat{V}_{h}^l]^2$ be an approximation of $\F$ satisfying Assumption~\ref{assumption::F_hAccuracy}. 
We generally construct, over the virtual computational domain $\F_h$,  VEMs of order $k\ge 1$ with the understanding that an isoparametric method is obtained by setting $l=k$.

\subsection{Reference IsoVEM formulation}\label{sec::MethodIFormulation}
The reference-frame isoparametric approach is to consider a discretization of the problem transformed back onto the reference domain. 

We operate under the assumption that the transformation  of the problem's data from the physical to the reference frame   is readily available  or can be estimated with sufficient accuracy. Namely, we assume we are given $\hat{a}$, $\hat{\mathbf{b}}$, $\hat{c}$, and $\hat{f}$ defined over $\RefDomain$ and satisfying the same properties imposed on the original data. 

The problem written in the reference frame reads: find  $\hat{u} \in H^1_0(\RefDomain)$ such that
\begin{equation}\label{eq::EllipticModelPDE_ref}
   \hat{\mathcal{A}}(\hat{u},\hat{v}) := \hat{A}(\hat{u},\hat{v}) + \hat{B}(\hat{u},\hat{v}) + \hat{C}(\hat{u},\hat{v}) = \hat{l}(\hat{v}) \ \ \ \forall \hat{v} \in H^1_0(\RefDomain),
\end{equation}
where, for all $\hat{u},\hat{v} \in H^1_0(\RefDomain)$,
\begin{align*}
    \hat{A}(\hat{u},\hat{v}) &= \int_{\RefDomain} \hat{a}\ \Jacobian_{\F}^{-T} \RefGrad \hat{u} \cdot \Jacobian_{\F}^{-T} \RefGrad \hat{v}\ j\ d\RefCoordinate,\qquad\qquad
    \hat{C}(\hat{u},\hat{v}) = \int_{\RefDomain} \hat{\mu}\  \hat{u}\ \hat{v}\ j\ d\RefCoordinate,
    \\
    \hat{B}(\hat{u},\hat{v}) &= \frac{1}{2} \int_{\RefDomain} \hat{\mathbf{b}} \cdot (\hat{v} \Jacobian_{\F}^{-T} \RefGrad \hat{u} - \hat{u} \Jacobian_{\F}^{-T} \RefGrad \hat{v})  \ j\ d\RefCoordinate,
     \quad
    \hat{l}(\hat{v}) = \int_{\RefDomain}   \hat{f} \ \hat{v}\ j\ d\RefCoordinate.
\end{align*}

We exploit the virtual element map $\F_h$ to approximate the components of the transformed weak formulation, such as the Jacobian and its determinant.
In particular, we produce a computable approximation of the Jacobian $\Jacobian_{\F}$ via $L^2$-projection:
\begin{equation}\label{eq::ProjJacobian}
    \Jacobian_{\F,h} := \Pi^1_{l-1} \F_h.
\end{equation}
From this, an approximate determinant,
\[
j_h := \det{\Jacobian_{\F,h}},
\]
is also available as a polynomial of degree $2(l-1)$. 

\begin{remark}
The idea of using virtual element maps to transfer the solution and the Jacobian across meshes  was first introduced in~\cite{Lipnikov2019AMeshes} where the Jacobian was approximated, instead of with the $L^2$-projector, with the $H^1$-semi norm or gradient projection $\Pi^\nabla$, cf. Remark~\ref{remark:proj}. Here, we take the idea one step further by using virtual elements to solve the problem at hand, rather than simply as a solution transfer tool.
\end{remark}

We are now ready to introduce the reference IsoVEM of order  $k \ge 1$.
The solution will be sought in the reference virtual element space $\hat{V}_{h,0}^k$. On each $\hat{E}\in\RefMesh$ we define the local forms by
\begin{align*}
    A_h^{\hat{E}}(\hat{u},\hat{v}) &:= \int_{\hat{E}} \hat{a}\ \Jacobian_{\F,h}^{-T} \Pi^1_{k-1} \hat{u}_h \cdot \Jacobian_{\F,h}^{-T} \Pi^1_{k-1} \hat{v}_h\ j_h\ d\RefCoordinate,
    \\
    B_h^{\hat{E}}(\hat{u},\hat{v}) &:= \frac{1}{2} \int_{\hat{E}} \hat{\mathbf{b}} \cdot (\PiO_{k} \hat{v}_h\ \Jacobian_{\F,h}^{-T} \Pi^1_{k-1} \hat{u}_h - \PiO_{k} \hat{u}_h\ \Jacobian_{\F,h}^{-T} \Pi^1_{k-1} \hat{v}_h)\ j_h\ d\RefCoordinate,
    \\
    C_h^{\hat{E}}(\hat{u},\hat{v}) &:= \int_{\hat{E}} \hat{\mu}\  \PiO_{k} \hat{u}_h\ \PiO_{k} \hat{v}_h\ j_h\ d\RefCoordinate,  
    \\
    l_h^{\hat{E}}(\hat{v}) &:= \int_{\hat{E}} \hat{f}\ \PiO_{k}\ \hat{v}_h\ j_h\ d\RefCoordinate, 
\end{align*}
for all $\hat{u}, \hat{v}\in V_k(\hat{E})$. 
Our analysis is valid for any stabilisation choice satisfying Assumption \ref{assumption::Stab::Stab}.
In the numerical experiments of Section~\ref{sec::NumericalTests}, we use the  appropriately  weighted dofi-dofi stabilisation~\cite{ConNonConVEM}:
\begin{equation}
    S^{\hat{E}}(\hat{u} , \hat{v}) = (\bar{a} + h_{\hat{E}}^2 \bar{\mu}) \sum_{i}
 \text{dof}_i(\hat{u}_h ) \cdot \text{dof}_i(\hat{v}_h ), \label{eq::MethodIStab}
\end{equation}
where $\bar{\cdot}$ denotes the nodal average over $\hat{E}$ and $\{\text{dof}_i\}_i$ is the set of $\dofs$ operators associated with the $\dofs$ of Definition~\ref{definition::dofs}. A proof that this choice satisfies the above assumption can be found in~\cite{ConNonConVEM}.

We then define the discrete analogue of the global  forms in ~\eqref{eq::EllipticModelPDE_ref}  for $\hat{u}_h, \hat{v}_h\in  \hat{V}^k_{h}$ by
\begin{align}
    \hat{\mathcal{A}}_h(\hat{u}_h,\hat{v}_h) &:= \sum_{\hat{E} \in \RefMesh}\Big[ A_h^{\hat{E}}(\hat{u}_h,\hat{v}_h) + B_h^{\hat{E}}(\hat{u}_h,\hat{v}_h)
    + C_h^{\hat{E}}(\hat{u}_h,\hat{v}_h) \nonumber\\
    &\hspace{2cm}+ S^{\hat{E}}(\hat{u}_h - \Pi^0_k \hat{u}_h, \hat{v}_h - \Pi^0_k \hat{v}_h)\Big], 
    \label{eq::bilinear_ref}\\
    \hat{l}_h(\hat{v}_h) &:= \sum_{\hat{E} \in \RefMesh} l_h^{\hat{E}}(\hat{v}_h),
    \label{eq::linear_ref}
\end{align}
yielding the reference IsoVEM as: find $\hat{u}_h \in \hat{V}^k_{h,0}$ such that
\begin{align}\label{eq:refIsoVEM}
    \hat{\mathcal{A}}_h(\hat{u}_h,\hat{v}_h) =\hat{l}_h(\hat{v}_h) \qquad \forall \hat{v}_h \in \hat{V}^k_{h,0}.
\end{align}

\subsection{Physical IsoVEM formulation}\label{sec::MethodIIFormulation}

The physical IsoVEM is constructed over the physical approximate domain $\Omega_h$, cf.~\eqref{eq:omega_h}. Given that the mapped elements $E_h\in\mesh$ are curved in general, the curved virtual element space $V_h^k$  relative to  $\mesh$ given by~\eqref{eq::Method2GlobalVEMSpace} can be used as the solution space. With this, a VEM may be directly constructed as in~\cite{daVeigaCurvedVEM}.
However, our approach is not a direct application of the curved VEM.  Instead, we \emph{exploit the virtual nature of the curved polygonal elements} by accessing their geometrical features only through their  $L^2$-projection. 
As such, our approach only requires the availability of the degrees of freedom of $\F_h$
and hence  it \emph{does not} rely on full knowledge of the mapping $\F$ (and, in fact,  $\F_h$). Moreover, by using the projection of the element domains, we have a computable geometry on which we can more easily integrate polynomials exactly. 
It was already argued in \cite{daVeigaCurvedVEM}  that a sufficient amount of quadrature points must be used to compute integrals over curved edges and curvilinear polygons. These quadratures are required, for instance, to compute \emph{approximately} the projection operators which are core to the VEM. 
In our approach, we introduce two altogether new computable operators that serve as approximate projectors. The benefit to this is that we can compute these new operators using only the projected elements and the DoFs with minimal quadrature rules required.

In keeping with the VEM philosophy, we decompose $\F_h$ into a polynomial and non-polynomial component, namely
\begin{align*}
    \F_h = \Pi^0_l \F_h + (\F_h - \Pi^0_l \F_h).
\end{align*}
The physical virtual element space $V_h^k$ is constructed over the mesh $\mesh=\{E_h:=\F_h(\hat{E})\}_h$, cf. ~\eqref{eq:omega_h}, using~\eqref{eq::Method2GlobalVEMSpace}.
The projected element $\projE := \Pi^0_l \F_h(\hat{E})$ is used to construct the necessary projection operators as follows.  

We start with the approximate local $L^2$-projection. 
In the forthcoming definitions,  for $v_h \in V_k(E_h)$, we let $\{\text{dof}_\alpha(v_h)\}_\alpha$ denote the set of internal moment $\dofs$ of Definition~\ref{definition::dofs}. Without loss of generality, these are assumed to be ordered so that $\text{dof}_\alpha(v_h)$ corresponds to the internal moment with respect to the $\alpha^{{\rm th}}$ scaled monomial in $\mathcal{M}_{k}(E_h)$, cf. Definition~\ref{definition::dofs}. In the physical IsoVEM, we always consider the scaled monomials to be extended to $\pp_k(\rr^2)$ by evaluating the polynomial expansion of a given $m_\alpha \in \mathcal{M}_k(E_h)$ at any given point in $\rr^2$.

\begin{definition}
\label{def::PiStarValue}
The operator $\Pi^0_{h,k}: V_h^k\rightarrow \pp_k(\mesh)$ satisfies, on each  $E_h \in \VirtualMesh$,
\begin{align*}
\int_{\projE} \Pi_{h,k}^0 v_h\ m_\alpha\ d\PhysCoordinate = 
    \begin{cases}
    \abs{\projE} \text{dof}_\alpha (v_h) \qquad & \forall m_\alpha \in \mathcal{M}_{k-2}(E_h),\vspace{2mm} \\
    \displaystyle
    \int_{\projE} \Pi^*_k v_h\ m_\alpha\ d\PhysCoordinate \qquad & \forall m_\alpha \in \mathcal{M}_{k}(E_h) \backslash \mathcal{M}_{k-2}(E_h),
    \end{cases}
\end{align*}
for all $v_h \in V_k(E_h)$.
\end{definition}
Noting that, by definition $\int_{E_h} \Pi^*_k v_h q\ d\PhysCoordinate=\int_{E_h}  v_hq\ d\PhysCoordinate $, cf.~\eqref{eq::CurvedVEMSpace}, we observe that
\begin{equation}
\label{eq::pi_identity}
\int_{\projE} \Pi_{h,k}^0 v_h\ m_\alpha\ d\PhysCoordinate=\frac{\abs{\projE}}{\abs{E_h}}
\int_{E_h} v_h\ m_\alpha\ d\PhysCoordinate
\qquad\forall m_\alpha \in \mathcal{M}_{k-2}(E_h).
\end{equation}
Next we construct an approximation of the gradient projection as follows.
\begin{definition}
\label{def::PiStarValueGrad}
The operator $\Pi^1_{h,k-1}: V_h^k\rightarrow [\pp_{k-1}(\mesh)]^2$ satisfies, on each  $E_h \in \VirtualMesh$,
\begin{align*}
\int_{\projE} \Pi^1_{h,k-1} v_h \cdot \mathbf{m}_\alpha\ d\PhysCoordinate = &-  \int_{\projE} \Pi^0_{h,k} v_h \ \nabla \cdot \mathbf{m}_\alpha\ d\PhysCoordinate \\
&+ \int_{\partial E_h} v_h \mathbf{m}_{\alpha} \cdot \mathbf{n}\ dS \qquad \forall \mathbf{m}_{\alpha} \in [\mathcal{M}_{k-1}(E_h)]^2, 
\end{align*}
for all $v_h \in V_k(E_h)$.
\end{definition}

\begin{remark}
    In Definition \ref{def::PiStarValueGrad}, we exploit the full knowledge of the boundary of $E_h$ to compute edge integrals over $\partial E_h$ instead of $\partial \projE$. 
    The alternative approach of computing $\Pi^1_{h,k-1}$ using only $\partial \projE$ is also viable, at the expense of a more involved analysis, cf.\cite{WellsThesis}. Note that the latter approach would be the only option available in higher-dimensions as the elemental interfaces  are virtual in this case.
\end{remark}

\begin{remark}
We shall refer to these new operators as projections but note that they  are \emph{not} exact projections as they do not reproduce polynomials on $E_h$. 
The error estimates of Theorems~\ref{theorem::PiStarOpBound} and~\ref{theorem::PiStarValueGrad} in Section \ref{sec::MethodII} show that this inconsistency introduces errors which are appropriately controlled.
\end{remark}

We are now ready to introduce the physical IsoVEM of order  $k \in \mathbb{N}$ over the physical curved polygonal mesh $\VirtualMesh$.
On each $E_h\in\VirtualMesh$ we define the local forms
\begin{align*}
    A_h^{E_h}(u_h,v_h) &= \int_{\projE} \Exten{a}\ \Pi^1_{h,k-1} u_h \cdot \Pi^1_{h,k-1} v_h\ d\PhysCoordinate,
    \\\
    B_h^{E_h}(u_h,v_h) &= \frac{1}{2} \int_{\projE} \Exten{\mathbf{b}} \cdot (\Pi_{h,k}^0 v_h \Pi^1_{h,k-1} u_h - \Pi_{h,k}^0 u_h \Pi^1_{h,k-1} v_h) \ d\PhysCoordinate,
    \\
    C_h^{E_h}(u_h,v_h) &= \int_{\projE} \Exten{\mu}\  \Pi_{h,k}^0 u_h\ \Pi_{h,k}^0 v_h\  d\PhysCoordinate, \\
    l_h^{E_h}(v_h) &= \int_{\projE} \Exten{f}\ \Pi_{h,k}^0 v_h\ d\PhysCoordinate.
\end{align*}
Due to the fact that $\projE$ may fall outside of $\Omega$, the local forms rely on extensions of the PDE data of the problem to be available. 

We thus introduce the following forms by summing the local contributions over $\VirtualMesh$: 
\begin{align}
    \mathcal{A}_h(u_h,v_h) &= \sum_{E_h \in \VirtualMesh} \Big[ A_h^{E_h}(u_h,v_h) + B_h^{E_h}(u_h,v_h) + C_h^{E_h}(u_h,v_h) \nonumber \\
    &\hspace{2cm}+ S^{E_h}(u_h-\Pi_{h,k}^0 u_h,v_h - \Pi_{h,k}^0 v_h)\Big], \label{eq::bilinear_phy}\\
    l_h(v_h) &= \sum_{E_h \in \VirtualMesh} l_h^{E_h}(v_h).
    \label{eq::linear_phy}
\end{align}
Here, as before, $S^{E_h}$ represents any VEM stabilisation form satisfying Assumption~\ref{assumption::Stab::Stab} with $\Pi_{h,k}^0$ in place of $\Pi_k^0$. In Section \ref{sec::NumericalTests} we again choose a dofi-dofi form (cf. Equation \eqref{eq::MethodIStab})
\begin{equation}\label{eq::MethodIIStab}
    S^{E_h}(u,v) = (\bar{a} + h_{E_h}^2 \bar{\mu}) \sum_{i} \text{dof}_i(u) \cdot \text{dof}_i(v).
\end{equation}

The physical IsoVEM is then defined, recalling the definition of $V_h^k$ and $V_{h,0}^k$ in Equation \eqref{eq::Method2GlobalVEMSpace}, as follows: find $u_h \in V_{h,0}^k$ such that
\begin{align}\label{eq::physicalIsoVEM}
    \mathcal{A}_h(u_h,v_h) =l_h(v_h) \qquad \forall v_h \in V_{h,0}^k.
\end{align}

\begin{remark}
 An alternative formulation is possible defining the method directly over the physical isoparametric elements $E_h\in\mesh$. However, while exact integration of polynomials is directly available over the projected element through accurate enough quadrature rules defined over the reference elements, integration over the true physical elements would require an isoparametric subtriangulation of each $E_h$.
 Moreover, the physical element based formulation would not be computable in three space dimensions, as the faces of three-dimensional elements are virtual. 
For these reasons we prefer to construct the physical IsoVEM formulation over the projected elements and show in Section~\ref{sec::MethodII} that the  extra source of error so introduced does not reduce the order of accuracy of the method
\end{remark}

\section{Reference IsoVEM analysis}\label{sec::MethodI}
For simplicity, we constrain the regularity of the domain transformation prior to presenting our analysis. In what follows, $m \in \mathbb{N}$ is a fixed integer such that $\F \in \left[ H^{m+1}(\hat{\Omega}) \right]^2$. In all subsequent results, applicable to both the reference and physical IsoVEM, we exclusively consider VEM discretizations with degrees $k,l \leq m$.
We should note that here, $k$ denotes the degree of the solution space, while $l$ pertains to the mapping space. While an isoparametric method requires that $k=l$, our analysis encompasses the general case, illustrating the impact of mixing discretization degrees. 

We  record here some straightforward bounds which will be used repeatedly in the forthcoming analysis. These can be proved in any space dimension. However, in keeping with the rest of the paper, we present them here explicitly for the two-dimensional case only.

For a matrix $\mathbf{A}$ with entries $A_{i,j} \in W^{s}_{\infty}(\omega)$, $i,j=1,2$, and a vector $\mathbf{v} \in \left[ H^{s}(\omega) \right]^2$ for some $s \in \mathbb{N} \cup \{ 0 \}$, it holds that
\begin{align}\label{eq::SobolevNormProdRule}
    \norm{\mathbf{A} \mathbf{v}}_{s,\omega} \lesssim \norm{\mathbf{A}}_{s,\infty,\omega} \norm{\mathbf{v}}_{s,\omega},
\end{align}
with the hidden constant only depending on $s$. Further, if $\det{\mathbf{A}} > a_0 > 0$ a.e. in  $\omega$ then it holds
\begin{equation}\label{eq::DetAndInvStab}
    \norm{\det{\mathbf{A}}}_{0,\infty,\omega} \leq 2 \norm{\mathbf{A}}_{0,\infty,\omega}^2, \quad 
    \norm{\mathbf{A}^{-1}}_{0,\infty,\omega} \leq \frac{1}{a_0} \norm{\mathbf{C}_{\mathbf{A}}^\top}_{0,\infty,\omega} \leq \frac{1}{a_0} \norm{\mathbf{A}}_{0,\infty,\omega}.
\end{equation}
Here, we expressed the inverse of a matrix $\mathbf{A}$ in the form
$\mathbf{A}^{-1} = \frac{1}{\det{\mathbf{A}}} \mathbf{C}_{\mathbf{A}}^T$,
where $\mathbf{C}_{\mathbf{A}}$ is the matrix of cofactors of $\mathbf{A}$. 

In what follows, the matrix of cofactors of $\Jacobian_{\F}$ and $\Jacobian_{\F,h}$ are denoted by $\mathbf{C}_{\F}$ and $\mathbf{C}_{\F,h}$, respectively.

\subsection{Jacobian error \& stability estimates}\label{sec:Jac}
We introduce and analyse the error terms related to the approximation of the Jacobian and its cofactors, the determinant, and the inverse Jacobian of the reference to physical map. Namely,
\begin{align}\label{eq:JacErrors}
   \JacobianDiff &:= \Jacobian_{\F,h} - \Jacobian_{\F},\ \CofDiff := \mathbf{C}_{\F,h} - \mathbf{C}_\F,\ \DetDiff := j_h - j,\ \InvDiff := \Jacobian_{\F,h}^{-1} - \Jacobian_{\F}^{-1}.  
\end{align}

\begin{lemma}[Jacobian estimates]\label{lemma::ProjJacobian}
Let $\RefMesh$ satisfy Assumption~\ref{assumption::ShapeRegular}  and $\F_h$ satisfy Assumption~\ref{assumption::F_hAccuracy}. For the projected Jacobian ~\eqref{eq::ProjJacobian} and related errors in~\eqref{eq:JacErrors} there hold
\begin{align}
    \label{eq:JacError}
    \norm{\JacobianDiff}_{0,\infty,\hat{E}} &\lesssim h^s_{\hat{E}} \abs{\F}_{s+1,\infty,\hat{E}},\\
    \label{eq:JacStab}
    \norm{\Jacobian_{\F,h}}_{0,\infty,\hat{E}} &\lesssim \abs{\F}_{1,\infty,\hat{E}},\\
    \label{eq:CofError}
    \norm{\CofDiff}_{0,\infty,\hat{E}} &\lesssim h^s_{\hat{E}} \abs{\F}_{s+1,\infty,\hat{E}},\\
    \label{eq:CofStab}
    \norm{\mathbf{C}_{\F,h}}_{0,\infty,\hat{E}} &\lesssim \abs{\F}_{1,\infty,\hat{E}},\\
    \label{eq:DetError}
    \norm{\DetDiff }_{0,\infty,\hat{E}} &\lesssim h_{\hat{E}}^s  \norm{\F}_{s+1,\infty,\hat{E}}^2,\\
    \label{eq:DetStab}
    \norm{j_h}_{0,\infty,\hat{E}} &\lesssim \abs{\F}^2_{1,\infty,\hat{E}},
\end{align}
for all $\hat{E} \in \RefMesh$ and $s\in\{0,\dots, l\}$. The hidden constants depend only on the shape regularity of $\RefMesh$ and on $l$.
\begin{proof}
By definition of the approximate Jacobian \eqref{eq::ProjJacobian} and discrete mapping $\F_h$ and using the triangle inequality it holds that
\begin{align}
    \norm{\JacobianDiff}_{0,\infty,\hat{E}} \leq \norm{\Pi^1_{l-1} \F_h - \Pi^1_{l-1} \F}_{0,\infty,\hat{E}} + \norm{\Pi^1_{l-1} \F - \RefGrad \F}_{0,\infty,\hat{E}}. \label{eq::JacobianProofPt1}
\end{align}
We apply the inverse estimate found in \cite{HHOBook} (Lemma 1.25) to get
\[
    \norm{\Pi^1_{l-1} \F_h - \Pi^1_{l-1} \F}_{0,\infty,\hat{E}} \lesssim h_{\hat{E}}^{-1} \norm{\Pi^1_{l-1} \F_h - \Pi^1_{l-1} \F}_{0,\hat{E}}
    \lesssim h_{\hat{E}}^{s-1} \abs{\F}_{s+1,\hat{E}}
\]
thanks to the stability of $\Pi^1$ and Assumption \ref{assumption::F_hAccuracy}. 
Inserting this bound into Equation \eqref{eq::JacobianProofPt1} and directly applying Theorem \ref{theorem::Pi0Accuracy} provides
\[
    \norm{\JacobianDiff}_{0,\infty,\hat{E}} \lesssim h_{\hat{E}}^s \abs{\F}_{s+1,\infty,\hat{E}} + h_{\hat{E}}^{s-1} \abs{\F}_{s+1,\hat{E}} 
    \lesssim h_{\hat{E}}^s \abs{\F}_{s+1,\infty,\hat{E}},
\]
proving~\eqref{eq:JacError}, from which the stability bound~\eqref{eq:JacStab} follows by taking $s=0$ and applying the triangle inequality.

The cofactor stability bound~\eqref{eq:CofStab} is a direct consequence of~\eqref{eq::DetAndInvStab} and of~\eqref{eq:JacStab} while the error estimate~\eqref{eq:CofError} is trivial, given the definition of minors of $\Jacobian_{\F,h}$.

To show stability of the determinant of the Jacobian~\eqref{eq:DetStab}, we denote the terms of the discrete Jacobian matrix by $J_{h,i,j} = (\Jacobian_{\F,h})_{i,j}$,  $i,j=1,2$. Then, noting that $j_h$ is a linear combination of products of terms from $\Jacobian_{\F,h}$, 
applying the triangle inequality, and bounding $L^\infty$ norms leads to
\begin{align*}
    \norm{j_h}_{0,\infty,\hat{E}} &= \norm{J_{h,0,0} J_{h,1,1} - J_{h,1,0} J_{h,0,1}}_{0,\infty,\hat{E}} \leq 2 \norm{\Jacobian_{\F,h}}_{0,\infty,\hat{E}}^2,
\end{align*}
and thus~\eqref{eq:DetStab} follows applying~\eqref{eq:JacStab}. The error estimate~\eqref{eq:DetError} is proven with a similar argument as follows. By considering the entries of $\Jacobian_\F$, we have the error term of
\begin{align*}
    \DetDiff = J_{h,0,0} J_{h,1,1} - J_{0,0} J_{1,1} + J_{1,0} J_{0,1} - J_{h,1,0} J_{h,0,1}.
\end{align*}
Considering the first difference term, we can rewrite this as
\begin{align*}
   J_{h,0,0} J_{h,1,1} - J_{0,0} J_{1,1} = J_{h,0,0} (J_{h,1,1} - J_{1,1}) + J_{1,1}( J_{h,0,0}- J_{0,0}).
\end{align*}
Bounding this term using the triangle inequality and~\eqref{eq:JacError}  gives
\begin{align*}
    \norm{J_{h,0,0} J_{h,1,1} - J_{0,0} J_{1,1}}_{0,\infty,\hat{E}} &\leq \norm{\Jacobian_{\F,h}}_{0,\infty,\hat{E}} \norm{\JacobianDiff}_{0,\infty,\hat{E}} + \norm{\Jacobian_{\F}}_{0,\infty,\hat{E}} \norm{\JacobianDiff}_{0,\infty,\hat{E}} \\
    &\lesssim \abs{\F}_{1,\infty,\hat{E}} \norm{\JacobianDiff}_{0,\infty,\hat{E}} \\
    &\lesssim h_{\hat{E}}^s \abs{\F}_{1,\infty,\hat{E}} \abs{\F}_{s+1,\infty,\hat{E}}.
\end{align*}
By repeating this argument for the remaining terms in $\DetDiff$ and applying the triangle inequality we get
\begin{align*}
    \norm{\DetDiff}_{0,\infty,\hat{E}} &\lesssim h_{\hat{E}}^s \abs{\F}_{1,\infty,\hat{E}} \abs{\F}_{s+1,\infty,\hat{E}}.
\end{align*}
Bounding the semi-norms by the $W^\infty_{s+1}$ norm completes the proof.
\end{proof}
\end{lemma}

\begin{lemma}[Inverse Jacobian estimates]\label{lemma::InvError}
Let $\RefMesh$ satisfy Assumption~\ref{assumption::ShapeRegular} and $\F_h$ satisfy Assumption~\ref{assumption::F_hAccuracy} with $h$ sufficiently small.  For the inverse of the projected Jacobian~\eqref{eq::ProjJacobian} and related errors in~\eqref{eq:JacErrors}, there hold
\begin{align}
    j_h |_{\hat{E}} &\gtrsim j_0, \label{eq::LowerBoundProjDet} \\
    \norm{\Jacobian_{\F,h}^{-1}}_{0,\infty,\hat{E}} &\lesssim \abs{\F}_{1,\infty,\hat{E}}, \label{eq::InvProjStab} \\
    \norm{\InvDiff}_{0,\infty,\hat{E}} &\lesssim h^s_{\hat{E}} \norm{\F}_{s+1,\infty,\hat{E}}^3,
\end{align}
for all $\hat{E} \in \RefMesh$ and $s\in\{0,\dots, l\}$. The hidden constants only depend on $j_0$, the shape regularity of $\RefMesh$, and on $l$.
\begin{proof}
From the estimate~\eqref{eq:DetError}, there exists a $C>0$ such that 
\begin{align*}
    \abs{(j- j_h)(\RefCoordinate)} \leq C h_{\hat{E}} \norm{\F}_{2,\infty,\hat{E}}^2,
\end{align*}
observing that $j$ is continuous. 
 Then, we define
\begin{align*}
    h^* =  \frac{j_0}{2 C \norm{\F}_{2,\infty,\hat{E}}^2},
\end{align*}
such that for $h_E \leq h^*$ we have $\norm{\DetDiff}_{0,\infty,\hat{E}} \leq j_0 / 2$. Equation \eqref{eq::LowerBoundProjDet} is then deduced by using a continuity argument. 

Equation~\eqref{eq::InvProjStab} is proven by using~\eqref{eq::LowerBoundProjDet} and~\eqref{eq::DetAndInvStab}
\begin{align*}
    \norm{\Jacobian_{\F,h}^{-1}}_{0,\infty,\hat{E}} &\leq \norm{\frac{1}{j_h}}_{0,\infty,\hat{E}} \norm{\mathbf{C}_{\F,h}^\top}_{0,\infty,\hat{E}} \leq \frac{1}{j_0} \norm{\Jacobian_{\F,h}}_{0,\infty,\hat{E}}. 
\end{align*}
The bound is then concluded by applying Lemma \ref{lemma::ProjJacobian}.

Using the definition of $\JacobianDiff$ and substitution of $\mathbf{I} = \Jacobian_{\F} \Jacobian_{\F}^{-1}$ results in
$   \InvDiff = -\Jacobian_{\F,h}^{-1}\ \JacobianDiff\ \Jacobian_{\F}^{-1}$.
Expanding and bounding norms gives
\begin{align*}
    \norm{\InvDiff}_{0,\infty,\hat{E}} \leq \norm{\Jacobian_{\F,h}^{-1}}_{0,\infty,\hat{E}} \norm{\JacobianDiff}_{0,\infty,\hat{E}} \norm{\Jacobian_{\F}^{-1}}_{0,\infty,\hat{E}}.
\end{align*}
Using~\eqref{eq::DetAndInvStab}, Lemma~\ref{lemma::ProjJacobian}, and~\eqref{eq::InvProjStab} provides
\begin{align*}
    \norm{\InvDiff}_{0,\infty,\hat{E}} &\lesssim \abs{\F}_{1,\infty,\hat{E}} h^s_{\hat{E}} \abs{\F}_{s+1,\infty,\hat{E}} \frac{1}{j_0} \abs{\F}_{1,\infty,\hat{E}} \lesssim \frac{1}{j_0} h^s_{\hat{E}} \norm{\F}_{s+1,\infty,\hat{E}}^3.
\end{align*}
\end{proof}
\end{lemma}


\subsection{Well-posedness}\label{sec::MethodI::LaxMilgram}
The above estimates on the Jacobian error allow us to establish the  well-posedness of the reference IsoVEM formulation~\eqref{eq:refIsoVEM} resorting on the Lax-Milgram lemma.

\begin{theorem}[Reference IsoVEM well-posedness]\label{theorem::Method1WellPosed}
Let $\RefMesh$ satisfy Assumption~\ref{assumption::ShapeRegular}, $\F_h$ satisfy Assumption~\ref{assumption::F_hAccuracy}, and $h$ be small enough.
Then, the reference IsoVEM formulation~\eqref{eq:refIsoVEM} with stabilisation~\eqref{eq::MethodIStab} obeying Assumption~\ref{assumption::Stab::Stab} possesses a unique solution.
\begin{proof}
We start by establishing the continuity of the bilinear form $\hat{\mathcal{A}}_h$ given in~\eqref{eq::bilinear_ref}.
Let $\hat{u}_h , \hat{v}_h \in \hat{V}_{h,0}^k$.
Through applying the H\"older and  Cauchy-Schwarz 
inequalities and the operator norm inequality~\eqref{eq::SobolevNormProdRule} we have, on any $\hat{E} \in \RefMesh$,
\begin{align*}
    A_h^{\hat{E}}(\hat{u}_h,\hat{v}_h) &\leq \norm{\hat{a}}_{0,\infty,\hat{E}} \norm{\mathbf{C}_{\F,h}}_{0,\infty,\hat{E}} \norm{\Jacobian_{\F,h}^{-T}}_{0,\infty,\hat{E}}  \norm{\Pi_{k-1}^1 \hat{u}_h}_{0,\hat{E}} \norm{\Pi^1_{k-1} \hat{v}_h}_{0,\hat{E}},\\
    B_h^{\hat{E}}(\hat{u}_h,\hat{v}_h) &\leq \frac{1}{2} \norm{\hat{\mathbf{b}}}_{0,\infty,\hat{E}} \norm{\mathbf{C}_{\F,h}}_{0,\infty,\hat{E}}  \Big(\norm{\Pi^1_{k-1}\hat{u}_h}_{0,\hat{E}} \norm{\Pi^0_k\hat{v}_h}_{0,\hat{E}} \\
    &\hspace{5cm}+ \norm{\Pi^0_k \hat{u}_h}_{0,\hat{E}} \norm{\Pi^1_{k-1}\hat{v}_h}_{0,\hat{E}} \Big), \\
    C_h^{\hat{E}}(\hat{u}_h,\hat{v}_h) &\leq \norm{\hat{\mu}}_{0,\infty,\hat{E}} \norm{j_h}_{0,\infty,\hat{E}} \norm{ \Pi^0_k \hat{u}_h}_{0,\hat{E}} \norm{ \Pi^0_k \hat{v}_h}_{0,\hat{E}}.
\end{align*}
Applying the stability of the projection operators and bounding these norms by the $H^1$ norm leads to
\begin{align*}
    A_h^{\hat{E}}(\hat{u}_h,\hat{v}_h) &\leq  \norm{\hat{a}}_{0,\infty,\hat{E}} \norm{\mathbf{C}_{\F,h}}_{0,\infty,\hat{E}} \norm{\Jacobian_{\F,h}^{-T}}_{0,\infty,\hat{E}}  \norm{\hat{u}_h}_{1,\hat{E}} \norm{\hat{v}_h}_{1,\hat{E}},\\
    B_h^{\hat{E}}(\hat{u}_h,\hat{v}_h) &\leq \norm{\hat{\mathbf{b}}}_{0,\infty,\hat{E}} \norm{\mathbf{C}_{\F,h}}_{0,\infty,\hat{E}} \norm{\hat{u}_h}_{1,\hat{E}} \norm{\hat{v}_h}_{1,\hat{E}}, \\
    C_h^{\hat{E}}(\hat{u}_h,\hat{v}_h) &\leq \norm{\hat{\mu}}_{0,\infty,\hat{E}} \norm{j_h}_{0,\infty,\hat{E}} \norm{\hat{u}_h}_{1,\hat{E}} \norm{\hat{v}_h}_{1,\hat{E}}.
\end{align*}
Next we apply Lemmas \ref{lemma::ProjJacobian} and \ref{lemma::InvError} to get, respectively,
\begin{align*}
    A_h^{\hat{E}}(\hat{u}_h,\hat{v}_h) &\lesssim \norm{\hat{a}}_{0,\infty,\hat{E}}  \abs{\F}_{1,\infty,\hat{E}}^2  \norm{\hat{u}_h}_{1,\hat{E}} \norm{\hat{v}_h}_{1,\hat{E}},\\
    B_h^{\hat{E}}(\hat{u}_h,\hat{v}_h) &\lesssim \norm{\hat{\mathbf{b}}}_{0,\infty,\hat{E}}  \abs{\F}_{1,\infty,\hat{E}} \norm{\hat{u}_h}_{1,\hat{E}} \norm{\hat{v}_h}_{1,\hat{E}}, \\
    C_h^{\hat{E}}(\hat{u}_h,\hat{v}_h) &\lesssim \norm{\hat{\mu}}_{0,\infty,\hat{E}}\abs{\F}^2_{1,\infty,\hat{E}} \norm{\hat{u}_h}_{1,\hat{E}} \norm{\hat{v}_h}_{1,\hat{E}}.
\end{align*}
For the stabilization term, we employ the stability Assumption  \ref{assumption::Stab::Stab} and the continuity of $\hat{\mathcal{A}}(\cdot,\cdot)$ to get
\begin{align*}
    S^{\hat{E}}(\hat{u}_h - \Pi^0_k \hat{u}_h, \hat{v}_h - \Pi^0_k \hat{v}_h) \lesssim \norm{\hat{u}_h}_{1,\hat{E}} \norm{\hat{v}_h}_{1,\hat{E}}.
\end{align*}
Summing these bounds over $\RefMesh$ 
results in
\begin{align}
    \label{eq::continuity_ref}
    \hat{\mathcal{A}}_h(\hat{u}_h,\hat{v}_h) \lesssim \left( \abs{\F}_{1,\infty,\RefDomain} + \abs{\F}_{1,\infty,\RefDomain}^{2} +1 \right) \norm{ \hat{u}_h}_{1,\RefDomain} \norm{ \hat{v}_h}_{1,\RefDomain}\lesssim 
    \norm{ \hat{u}_h}_{1,\RefDomain} \norm{ \hat{v}_h}_{1,\RefDomain},
\end{align} 
thus establishing the continuity of $\hat{\mathcal{A}}_h$. The (hidden) continuity constant depends on shape regularity of $\RefMesh$, on $l$, on appropriate norms of the data and the $W^1_{\infty}$ semi-norm of $\F$, and  on the stability constant in Assumption~\ref{assumption::Stab::Stab}.
The continuity of the linear form $\hat{l}_h(\cdot)$ given in \eqref{eq::linear_ref} follows easily by the same argument used for the reaction term $C_h^{\hat{E}}$.

It remains to show the coercivity of the bilinear form $\hat{\mathcal{A}}_h$. Considering the $L^2$ norm of the $\Pi^1_{k-1}$ projection on an element $\hat{E} \in\RefMesh$ and assuming a sufficiently small $h$ such that $j_0 \leq j_h$ (see Lemma \ref{lemma::ProjJacobian}) gives,  for all $\hat{v}_h \in \hat{V}_{h,0}^k$,
    \begin{align*}
         a_0 j_0 \norm{\Pi^1_{k-1} \hat{v}_h}_{0,\hat{E}}^2 & \lesssim \int_{\hat{E}} \hat{a} \abs{\Pi^1_{k-1} \hat{v}_h}^2 j_h\ d\RefCoordinate 
        \lesssim \int_{\hat{E}} \hat{a} \abs{ \Jacobian_{\F,h}^T \Jacobian_{\F,h}^{-T} \Pi^1_{k-1} \hat{v}_h}^2 j_h\ d\RefCoordinate \\
        &\lesssim \norm{\Jacobian_{\F,h}^T}_{0,\infty,\hat{E}}^2 \int_{\hat{E}} \hat{a} \abs{ \Jacobian_{\F,h}^{-T} \Pi^1_{k-1} \hat{v}_h}^2 j_h\ d\RefCoordinate
        \\
        &\lesssim  \norm{\Jacobian_{\F,h}^T}_{0,\infty,\hat{E}}^2 A_h^{\hat{E}}(\hat{v}_h,\hat{v}_h)
        \lesssim A_h^{\hat{E}}(\hat{v}_h,\hat{v}_h),
    \end{align*}
    by definition of $A_h^{\hat{E}}$ and Lemma \ref{lemma::ProjJacobian}.
    Arguing similarly for the reaction term $C_h^{\hat{E}}$ and noting that  $B_h^{\hat{E}}(\hat{v}_h,\hat{v}_h) = 0$ we arrive at
    \begin{align}
          \mathcal{A}_h^{\hat{E}}(\hat{v}_h,\hat{v}_h) &\gtrsim \mu_0 j_0 \norm{\Pi^0_k \hat{v}_h}_{0,\hat{E}}^2 + a_0 j_0 \norm{\Pi^1_{k-1} \hat{v}_h}_{0,\hat{E}}^2 + S^{\hat{E}}(\hat{v}_h - \Pi^0_k \hat{v}_h, \hat{v}_h - \Pi^0_k \hat{v}_h).\label{eq::A_hCoercivityExpansionPt2}
    \end{align}
    Following standard VEM analysis \cite{ellipticVEM,ConNonConVEM}, we apply Assumption \ref{assumption::Stab::Stab} to get
    \begin{align*}
        S^{\hat{E}}(\hat{v}_h - \Pi^0_k \hat{v}_h, \hat{v}_h - \Pi^0_k \hat{v}_h) \gtrsim a_0 \norm{\RefGrad \hat{v}_h - \Pi^1_{k-1} \hat{v}_h}_{0,\hat{E}}^2 + \mu_0 \norm{\hat{v}_h - \Pi^0_{k} \hat{v}_h}_{0,\hat{E}}^2,
    \end{align*}
    from which substitution into Equation \eqref{eq::A_hCoercivityExpansionPt2} obtains
    \begin{align*}
          \mathcal{A}_h^{\hat{E}}(\hat{v}_h,\hat{v}_h) &\gtrsim j_0 \left( \mu_0 \norm{\hat{v}_h}_{0,\hat{E}}^2 + a_0 \norm{\RefGrad \hat{v}_h}_{0,\hat{E}}^2 \right) \gtrsim \norm{\hat{v}_h}_{1,\hat{E}}^2.
    \end{align*}

    The coercivity of $\hat{\mathcal{A}}_h$ now follows summing over elements and bounding norms yielding,  for all $\hat{v}_h \in \hat{V}_{h,0}^k$, 
    \begin{align}
    \label{eq::coercive_ref}
    \norm{\hat{v}_h}_{1,\RefDomain}^2 \lesssim \hat{\mathcal{A}}_h(\hat{v}_h,\hat{v}_h).
    \end{align}
    The hidden constant depends on the shape regularity of $\RefMesh$, on $l$, $j_0$, $a_0$, and $\mu_0$, the stability constant of Assumption \ref{assumption::Stab::Stab} and the $W^1_{\infty}$-seminorm of $\F$.
   The well-posedness of the reference IsoVEM formulation is thus given by the Lax-Milgram lemma.
\end{proof}
\end{theorem}

\subsection{Reference IsoVEM error estimate}\label{sec::MethodI::H1Error}
We now prove an $H^1$-norm error estimate for the reference IsoVEM. It is straightforward to check that the Strang-type bound of Theorem~\ref{theorem::StrangBound1} applies to the IsoVEM formulation~\eqref{eq:refIsoVEM}. 

The difficulty lays in analysing the  inconsistency terms in~\eqref{eq::StrangBound}.
These include new inconsistencies produced by the virtual map and hence we cannot simply resort on the analysis in~\cite{ConNonConVEM,ellipticVEM}. Instead, the analysis relies on the bounds of Section~\ref{sec:Jac} as follows.

\begin{lemma}[Reference IsoVEM consistency error] \label{lemma::MethodIConsistencyError} 
    Let $\RefMesh$ satisfy Assumption~\ref{assumption::ShapeRegular}, $\F_h$ satisfy Assumption~\ref{assumption::F_hAccuracy}, and $h$ be small enough.
    Assume that $\hat{a},\hat{\mu}, \hat{\mathbf{b}}  \in W^{s+1}_{\infty}(\RefDomain)$ and $\hat{f} \in  H^{s-1}(\RefDomain)$, for some $s\in\{1,\dots, \min\{k,l\}\}$. For all $\hat{v}_h\in \hat{V}^k_{h}$ and $\hat{p}\in\pp_k(\RefMesh)$, there hold
    \begin{align}
        \sum_{\hat{E}\in\RefMesh}\abs{\hat{\mathcal{A}}^{\hat{E}}_h(\hat{p},\hat{v}_h) - \hat{\mathcal{A}}^{\hat{E}}(\hat{p},\hat{v}_h)} &\lesssim h^s \norm{\hat{p}}_{s+1,\RefMesh} \norm{\hat{v}_h}_{1,\RefDomain},
        \label{eq::a_hat_cons}\\
        \abs{\hat{l}_h(\hat{v}_h) - \hat{l}(\hat{v}_h)} &\lesssim h^s \norm{\hat{f}}_{s-1,\RefDomain} \norm{\hat{v}_h}_{1,\RefDomain}.\label{eq::l_hat_cons}
    \end{align}
    The hidden constants depend on the shape regularity of $\RefMesh$, on $s$ and $j_0$, and on the $W^{s+1}_{\infty}$-norm of  $\F$. The first constant depends also on the $W^{s+1}_{\infty}(\RefDomain)$-norm of  $\hat{a}$, $\hat{\mathbf{b}}$, and $\hat{\mu}$.
    \begin{proof}
    Let $\hat{E}\in\RefMesh$.
        For brevity, we only detail the estimation of the diffusion term $A_h^{\hat{E}}(\hat{p},\hat{v}_h) - A^{\hat{E}}(\hat{p},\hat{v}_h)$; all other terms can be treated similarly as in~\cite{ConNonConVEM}. Adding and subtracting terms gives
        \begin{align*}
             A^{\hat{E}}_h(\hat{p},\hat{v}_h) - A^{\hat{E}}(\hat{p},\hat{v}_h) 
            =& \int_{\hat{E}} \hat{a} \CofDiff \RefGrad \hat{p} \cdot \Jacobian^{-T}_{\F} \RefGrad \hat{v}_h\ d\RefCoordinate 
            + \int_{\hat{E}} \hat{a} \mathbf{C}_{\F,h} \RefGrad \hat{p} \cdot \InvDiff^T \RefGrad \hat{v}_h\ d\RefCoordinate\\
            &+ \int_{\hat{E}} \hat{a} \Jacobian_{\F,h}^{-T} \RefGrad \hat{p} \cdot \mathbf{C}_{\F,h} (\RefGrad - \Pi^1_{k-1})  \hat{v}_h\ d\RefCoordinate\\
            =:& \; T_1 + T_2 + T_3.
        \end{align*}
        Combining H\"{o}lder and Cauchy-Schwarz inequalities with~\eqref{eq::DetAndInvStab} we have
        \begin{align*}
            |T_1| &\leq \norm{\hat{a}}_{0,\infty,\hat{E}} \norm{\CofDiff}_{0,\infty,\hat{E}} \norm{\Jacobian^{-T}_{\F}}_{0,\infty,\hat{E}} \abs{\hat{p}}_{1,\hat{E}} \abs{\hat{v}_h}_{1,\hat{E}} \\
            &\leq \norm{\hat{a}}_{0,\infty,\hat{E}} \norm{\CofDiff}_{0,\infty,\hat{E}} \frac{1}{j_0} \norm{\mathbf{C}_{\F}}_{0,\infty,\hat{E}} \abs{\hat{p}}_{1,\hat{E}} \abs{\hat{v}_h}_{1,\hat{E}},
        \end{align*}
        from which, applying  Lemma \ref{lemma::ProjJacobian}, yields
        \begin{align*}
            |T_1| &\lesssim \norm{\hat{a}}_{0,\infty,\hat{E}} \frac{1}{j_0}  h_{\hat{E}}^s \norm{\F}_{s+1,\infty,\hat{E}}^2 \abs{\hat{p}}_{1,\hat{E}} \abs{\hat{v}_h}_{1,\hat{E}}.
        \end{align*}
        A similar argument using, instead, Lemma \ref{lemma::InvError} provides
        \begin{align*}    
            |T_2| &\lesssim \left\{ \norm{\hat{a}}_{0,\infty,\hat{E}} \frac{1}{j_0}  \norm{\F}_{s+1,\infty,\hat{E}}^4 \right\} h_{\hat{E}}^s\abs{\hat{p}}_{1,\hat{E}} \abs{\hat{v}_h}_{1,\hat{E}}
        \end{align*}
        Combining these bounds and hiding constants leads to
        \begin{align*}
            |T_1| + |T_2| &\lesssim h_{\hat{E}}^s \norm{\hat{p}}_{1,\hat{E}} \norm{\hat{v}_h}_{1,\hat{E}}.
        \end{align*}
        To bound $T_3$, we emulate the steps taken in \cite{ConNonConVEM}, writing
        \begin{align*}
            T_3 &= \int_{\hat  {E}} \hat{a} \Jacobian_{\F,h}^{-T} \RefGrad \hat{p} \cdot \mathbf{C}_{\F,h} (\RefGrad \hat{v}_h - \Pi^1_{k-1} \hat{v}_h)\  d\RefCoordinate\\ 
            &= \int_{\hat{E}} \boldsymbol{\gamma} \cdot (\RefGrad  \hat{v}_h- \Pi^1_{k-1} \hat{v}_h)\   d\RefCoordinate\\
            &=\int_{\hat{E}} (\boldsymbol{\gamma} - \Pi^0_{k-1} \boldsymbol{\gamma} ) \cdot (\RefGrad \hat{v}_h- \Pi^1_{k-1}\hat{v}_h)  \  d\RefCoordinate.
        \end{align*}
        noting that $\mathbf{C}_{\F,h}$ has polynomial entries of degree $(l-1)$, introducing $\boldsymbol{\gamma} := ( \hat{a} \Jacobian_{\F,h}^{-T} \RefGrad \hat{p} )^{T} \mathbf{C}_{\F,h}$. 
        Then, applying Cauchy-Schwarz and Theorem \ref{theorem::Pi0Accuracy}, gives
        \begin{align}
            |T_3| \leq \norm{\boldsymbol{\gamma} - \Pi^0_{k-1} \boldsymbol{\gamma}}_{0,\hat{E}} \norm{(\RefGrad - \Pi^1_{k-1}) \hat{v}_h}_{0,\hat{E}} 
            \lesssim h_{\hat{E}}^{s} \norm{\boldsymbol{\gamma}}_{s,\hat{E}} \norm{ \hat{v}_h}_{1,\hat{E}}.
            \label{eq::ConsistencyPt6}
        \end{align}
        Note that $\boldsymbol{\gamma}$ can be written as
        $\boldsymbol{\gamma} = \hat{a} \frac{1}{j_h} (\mathbf{C}_{\F,h} \RefGrad \hat{p})^\top C_{\F,h}$,
        and thus a bound of norms of ${1}/{j_h}$ is required. To this end, we use the chain rule and quotient rule yielding, for any $|\alpha|$ multi-index with $|\alpha|\le s$ and noting that $j_0 \in (0,1)$,
        \begin{align*}
            \left|D^{\alpha} \frac{1}{j_h}\right|  = \left|(-1)^{\abs{\alpha}} \frac{\alpha !}{j_h^{\abs{\alpha} +1}} D^{\alpha} j_h\right| &\leq \frac{s!}{\abs{j_h}^{\abs{\alpha}+1}} \abs{D^{\alpha} j_h} \lesssim s!\left( \frac{1}{j_0}\right)^{s+1}  \norm{j_h}_{s,\infty,\hat{E}},
        \end{align*}
        and we conclude that 
         $   \norm{{1}/{j_h}}_{s,\infty,\hat{E}} \lesssim \norm{j_h}_{s,\infty,\hat{E}}$
         holds with the hidden constant depending on $s$ and $j_0$.
        Hence,  by repeated applications of~\eqref{eq::SobolevNormProdRule} we arrive at
        \begin{align}
            \norm{\boldsymbol{\gamma}}_{s,\hat{E}} &\lesssim \norm{\hat{a}}_{s,\infty,\hat{E}} \norm{j_h}_{s,\infty,\hat{E}} \norm{\mathbf{C}_{\F,h}}_{s,\infty,\hat{E}}^2 \norm{\RefGrad \hat{p}}_{s,\hat{E}}. \label{eq::ConsistencyPt0}
        \end{align}
        We further bound $\norm{j_h}_{s,\infty,\hat{E}}$  adding and subtracting $\Pi^0_{2l-2} j$ and applying the triangle inequality, yielding
        \begin{align*}
            \norm{j_h}_{s,\infty,\hat{E}} &\leq \norm{j_h - \Pi^0_{2l-2} j}_{s,\infty,\hat{E}} + \norm{\Pi^0_{2l-2} j}_{s,\infty,\hat{E}} \\
            &= \norm{\Pi^0_{2l-2} \DetDiff}_{s,\infty,\hat{E}} + \norm{\Pi^0_{2l-2} j}_{s,\infty,\hat{E}},
        \end{align*}
        from which, an inverse inequality
         \cite{HHOBook} (Corollary 1.29) gives
        \begin{align*}
            \norm{j_h}_{s,\infty,\hat{E}} &\lesssim h_{\hat{E}}^{-s} \norm{\Pi^0_{2l-2} \DetDiff}_{0,\infty,\hat{E}} + \norm{\Pi^0_{2l-2} j}_{s,\infty,\hat{E}}\\
            &\lesssim h_{\hat{E}}^{-s} \norm{\DetDiff}_{0,\infty,\hat{E}} + \norm{j}_{s,\infty,\hat{E}}.
        \end{align*}
        Lemma \ref{lemma::ProjJacobian} is then applied to the first term and the second term is bounded using~\eqref{eq::SobolevNormProdRule}, giving
        \begin{align*}
            \norm{j_h}_{s,\infty,\hat{E}} \lesssim \norm{\F}_{s+1,\infty,\hat{E}}^2. 
        \end{align*}
        A similar argument yields
        \begin{align*}
            \norm{\mathbf{C}_{\F,h}}_{s,\infty,\hat{E}} \lesssim \norm{\F}_{s+1,\infty,\hat{E}},  \label{eq::ConsistencyPt4}
        \end{align*}
        and, inserting these bounds in \eqref{eq::ConsistencyPt0}, we obtain
        \begin{align*}
            \norm{\boldsymbol{\gamma}}_{s,\hat{E}} &\lesssim \norm{\hat{a}}_{s,\infty,\hat{E}} \norm{\F}_{s+1,\infty,\hat{E}}^4 \norm{\RefGrad \hat{p}}_{s,\hat{E}} 
            \lesssim \norm{\hat{p}}_{s+1,\hat{E}}.
        \end{align*}
        Using this in \eqref{eq::ConsistencyPt6} completes the bound of $T_3$ which, together with the  bound of $T_1+T_2$ gives
        \begin{align*}
            \abs{A^{\hat{E}}(\hat{p},\hat{v}_h) - A_h^{\hat{E}}(\hat{p},\hat{v}_h)} &\lesssim h_{\hat{E}}^s \norm{\hat{p}}_{1,\hat{E}} \norm{\hat{v}_h}_{1,\hat{E}} + h_{\hat{E}}^{s} \norm{\boldsymbol{\gamma}}_{s,\hat{E}} \norm{ \hat{v}_h}_{1,\hat{E}} 
            \\
            &\lesssim h_{\hat{E}}^{s} \norm{\hat{p}}_{s+1,\hat{E}} \norm{ \hat{v}_h}_{1,\hat{E}}
        \end{align*}
        A similar line of reasoning provides bounds for the lower order terms, while we note that the stabilization term is zero in this case. 
        Summing up the contributions from all elements   already gives~\eqref{eq::a_hat_cons}. The proof of~\eqref{eq::l_hat_cons} is similar and omitted for brevity.
    \end{proof}
\end{lemma}

As for the analogous Theorem~\ref{theorem::CurvedH1Error},  we can now conclude the analysis of the reference IsoVEM  by inserting in the Strang-type bound~\eqref{eq::StrangBound} the results from Theorem~\ref{theorem::Pi0Accuracy}, Theorem~\ref{theorem::GlobalCurvedVEMInterp}, and the above Lemma~\ref{lemma::MethodIConsistencyError}; see eg.~\cite{ConNonConVEM} for the details.

\begin{theorem}[Reference IsoVEM error estimate]\label{theorem::MethodIH1Error} 
    Let $\RefMesh$ satisfy Assumption~\ref{assumption::ShapeRegular}, $\F_h$ satisfy Assumption~\ref{assumption::F_hAccuracy}, and $h$ be small enough.
    Assume that $\hat{a},\hat{\mu}, \hat{\mathbf{b}}  \in W^{s+1}_{\infty}(\RefDomain)$ and $\hat{f} \in  H^{s-1}(\RefDomain)$, for some $s\in\{1,\dots, \min\{k,l\}\}$.
    Suppose additionally that  $\hat{u} \in H^{s+1}(\RefDomain)$. Then, the reference IsoVEM solution $\hat{u}_h \in \hat{V}^k_{h,0}$ of~\eqref{eq:refIsoVEM} satisfies
    \begin{equation*}
        \norm{\hat{u} - \hat{u}_h}_{1,\RefDomain} 
        + \norm{\nabla\hat{u} - \Pi^1_{k-1}\hat{u}_h}_{0,\RefDomain}
        \lesssim h^s \left( \norm{\hat{u}}_{s+1,\RefDomain} + \norm{\hat{f}}_{s-1,\RefDomain} \right).
    \end{equation*}
     The hidden constant depends on the shape regularity of $\RefMesh$, on $s$, $l$, and $k$, the lower bounds $j_0$, $a_0$, $\mu_0$, and the $W^{s+1}_\infty(\RefDomain)$ norms of $\hat{a}$, $\hat{\mathbf{b}}$, $\hat{\mu}$, and $\F$.
\end{theorem}

\section{Physical IsoVEM analysis} \label{sec::MethodII}
The main difficulty in the analysis of the physical IsoVEM is in assessing the error due to the use of the projected elements in both the definition of the projection operators and the method itself. We thus proceed with a series of results aimed at assessing these sources of errors. One important step is given in Lemma~\ref{lemma::MeshScalingV2} below, showing that the mesh sizes of the physical elements scale like the respective reference and projected elements, and thus can be used interchangeably.

We recall from Section~\ref{sec::MethodIIFormulation} that the physical IsoVEM~\eqref{eq::physicalIsoVEM} is defined over the virtual domain $\Omega_h$ by fixing the physical mesh $\mesh$ and corresponding virtual element space $V_h^k$ via a reference mesh $\RefMesh$ and corresponding  virtual map $\F_h$.

\subsection{Projected elements}
We quantify the effects of using the projected elements for integration in the following technical results.

\begin{lemma}[Projected element area approximation]\label{lemma::AreaError} 
Let $\RefMesh$ satisfy Assumption~\ref{assumption::ShapeRegular}, $\F_h \in \left[ \hat{V}_h^l \right]^2$ satisfy Assumption~\ref{assumption::F_hAccuracy}, and $h$ be small enough.
For each reference element $\hat{E} \in \RefMesh$ with corresponding physical and projected elements $E_h \in \mesh$ and  $\projE$, respectively, 
there holds 
    \begin{align*}
        |\abs{E_h} - \abs{E_h^\pi} | &\lesssim h_{\hat{E}}^{s+2} \norm{\F}_{s+1,\infty,\hat{E}}^2, \\
        \abs{E_h \triangle \projE} = \abs{E_h \backslash \projE}+\abs{\projE \backslash E_h} &\lesssim h_{\hat{E}}^{s+2} \norm{\F}^{3}_{s+1,\infty,\hat{E}}, \\
        \abs{E_h \triangle E} = \abs{E_h \backslash E}+\abs{E \backslash E_h} &\lesssim h_{\hat{E}}^{s+2} \norm{\F}^{3}_{s+1,\infty,\hat{E}},
    \end{align*}
    with $s \in \{ 0,1,...,l\}$.
    The symbol $\triangle$ is used to denote the symmetric difference between the two sets. So, for instance, $E_h \triangle \projE = (E_h \backslash \projE) \cup (\projE \backslash E_h)$. 
    The hidden constants depend on the shape regularity of $\RefMesh$ and on $l$.
\begin{proof}
     By similar arguments  as in  Lemma~\ref{lemma::InvError}, for $h$ small enough, the determinants 
    $\det{\RefGrad \F_h}$ and $ \det{\RefGrad \Pi^0_l \F_h}$
    are strictly positive. Hence, \holder's inequality gives
    \begin{align}
        \abs{\abs{E_h} - \abs{E_h^\pi}} &= \abs{ \int_{\hat{E}} \left(\det{\RefGrad \F_h} -  \det{\RefGrad \Pi^0_l \F_h}\right)\ d\RefCoordinate
        }\nonumber \\
        &\leq h^2_{\hat{E}} \norm{\det{\RefGrad \F_h} -  \det{\RefGrad \Pi^0_l \F_h}}_{0,\infty,\hat{E}}. \label{eq::measureLemmaPt1}
    \end{align}
    Adding and subtracting the determinants of ${\RefGrad \F}$ and ${\RefGrad \Pi^0_l \F}$ as in the proof of Lemma \ref{lemma::ProjJacobian} we easily get
    \begin{align*}
        \norm{\det{\RefGrad \F_h} -  \det{\RefGrad \Pi^0_l \F_h}}_{0,\infty,\hat{E}} &\lesssim h_{\hat{E}}^s \norm{\F}_{s+1,\infty,\hat{E}}^2,
    \end{align*}
    and substitution of this into Equation \eqref{eq::measureLemmaPt1} completes the first estimate.

In view of  bounding the set difference $ \abs{\projE \backslash E_h}$,  we reason that this set must be contained within a band around $\projE \backslash E_h$ of thickness $ r= \norm{\F_h - \Pi^0_l \F_h}_{0,\infty,\hat{E}}$. In turn, this can be estimated from above considering the bands constructed over each edge of $E_h$ plus the area of circles of the same radius $r$ centred on each vertex. Thus,
\begin{align*}
        \abs{\projE \backslash E_h} &\le  \sum_{e_h \subset \partial E_h} \left(
        \norm{\F_h - \Pi^0_l \F_h}_{0,\infty,\hat{E}} 
        \abs{e_h} +
        \frac{\pi}{4} \norm{\F_h - \Pi^0_l \F_h}_{0,\infty,\hat{E}}^2 
        \right).
    \end{align*}
From this, using once more the arguments of Lemma \ref{lemma::ProjJacobian}, we obtain
\begin{align*}
        \abs{\projE \backslash E_h} &\lesssim
        h_{\hat{E}}^{s+1} \norm{\F}_{s+1,\infty,\hat{E}}\left( h_{\hat{E}}\abs{\F_h}^2_{1,\infty,\hat{E}}+\frac{\pi}{4} h_{\hat{E}}^{s+1} \norm{\F}_{s+1,\infty,\hat{E}} \right) \lesssim
        h_{\hat{E}}^s \norm{\F}_{s+1,\infty,\hat{E}}^3.
\end{align*}
A similar argument follows to provide bounds for $\abs{E_h \backslash \projE}$. The bounds of $\abs{E \backslash E_h}$ and $\abs{E_h \backslash E}$ also follow similarly making use of Assumption~\ref{assumption::F_hAccuracy}.
\end{proof}
\end{lemma}

When quantifying the error of the physical IsoVEM, we have to consider the scaling of the mesh size terms $h_\omega$ for a choice of $\omega = \hat{E}, E_h, \projE$. A first consequence of the above lemma is that reference, physical, and projected elements mesh sizes can, in fact, be used interchangeably.
\begin{lemma}[equivalence of mesh sizes]\label{lemma::MeshScalingV2} 
Let $\RefMesh$ satisfy Assumption~\ref{assumption::ShapeRegular}, $\F_h$ satisfy Assumption~\ref{assumption::F_hAccuracy}, and $h$ be small enough. For each $\hat{E}\in \RefMesh$, we have
\begin{align}
    h_{\hat{E}} \sim h_{E_h} \sim h_{\projE}
\end{align}
where the hidden constants depend on $j_0$ and the shape regularity of $\hat{E}$. The second  hidden constant depends also on the $W^1_\infty(\hat{E})$-norm of $\F$.
\begin{proof}
We first show that $h_{\hat{E}} \sim h_{E_h}$. 
Using 
the shape regularity of $\hat{E}$ and the arguments of Lemma \ref{lemma::InvError} we have
\begin{align}
    j_0 h_{\hat{E}}^2 &\lesssim j_0 \abs{\hat{E}}  \lesssim \abs{E_h} \lesssim h_{E_h}^2.
\end{align}
Taking square roots completes the lower bound. The upper bound is given using the shape regularity of $E_h$,  computing the area of $E_h$ using $\hat{E}$, applying \holder's inequality and the stability of the VEM interpolation
\begin{align}
    h_{E_h}^2 \lesssim \abs{E_h} &\leq \norm{\det{\RefGrad \F_h}}_{0,\infty,\hat{E}} \abs{\hat{E}} \leq 2 \abs{\F_h}^2_{1,\infty,\hat{E}} \abs{\hat{E}} \lesssim \abs{\F}_{1,\infty,\hat{E}} h_{\hat{E}}^2.
\end{align}
The second equivalence relation, namely $h_{E_h} \sim h_{\projE}$, follows easily from the first and Lemma~\ref{lemma::AreaError}.
\end{proof}
\end{lemma}

\begin{remark}
    We note that if the computational reference domain $\RefDomain$ is defined as being \emph{close to} $\Omega$, then the hidden constants in Lemma \ref{lemma::MeshScalingV2} are close to $1$.
\end{remark}

\subsection{Extension of polynomials}
For polynomials $p \in \pp_k(E_h)$, $E_h\in\mesh$, we consider their natural extensions to the whole of $\rr^2$, still denoted by $p$, in view of their use over the projected element $\projE$.
We emphasise that this extension does not satisfy Stein's extension Theorem~\ref{theorem::SteinExtension}. Instead we have the following local polynomial extension bound.

\begin{lemma}[Polynomial extension]\label{lemma::PolyExten} 
Let $\RefMesh$ satisfy Assumption~\ref{assumption::ShapeRegular}, $\F_h$ satisfy Assumption~\ref{assumption::F_hAccuracy}, and $h$ be small enough.
For each physical element $E_h \in \mesh$  and corresponding projected element $\projE$, 
there holds 
\begin{align*}
    \norm{p}_{0,\infty,\projE} \lesssim \norm{p}_{0,\infty,E_h},
\end{align*}
for all $p \in \pp_k(E_h)$.
The hidden constant depends on the shape regularity of $\RefMesh$, on $k$, and on the $W^1_\infty(\hat{E})$-norm of $\F$.
\begin{proof}
We proceed by bounding the ratio between the two norms $R:={\norm{p}_{0,\infty,\projE}}/{\norm{p}_{0,\infty,E_h}}$.
Applying an inverse estimate \cite[Lemma 1.25]{HHOBook} to the numerator and a Sobolev inequality \cite{Adams2003} to the denominator of $R$ gives
\begin{align*}
    R &\lesssim \frac{h_{\projE}^{-1} \norm{p}_{0,\projE}}{h_{E_h}^{-1} \norm{p}_{0,E_h}} \lesssim \frac{ \norm{p}_{0,\projE}}{\norm{p}_{0,E_h}},
\end{align*}
by the equivalence $h_{\projE}\sim h_{E_h}$ from Lemma \ref{lemma::MeshScalingV2}.
Passing to squares and expanding the numerator gives
\begin{align}
    R^2 &\lesssim \frac{ \norm{p}^2_{0,\projE}}{\norm{p}^2_{0,E_h}} =  \left( \frac{\norm{p}^2_{0,E_h} - \norm{p}^2_{0,E_h \backslash \projE} + \norm{p}^2_{0,\projE \backslash E_h}}{\norm{p}^2_{0,E_h}} \right) \lesssim \left( 1 + \frac{\norm{p}^2_{0,\projE \backslash E_h}}{\norm{p}^2_{0,E_h}} \right). \label{eq::Rpt1}
\end{align}
The ratio of norms in the bound above can be further bounded by applying \holder's inequality on the numerator and, once again, an inverse estimate on the denominator to give
\begin{align*}
    \frac{\norm{p}^2_{0,\projE \backslash E_h}}{\norm{p}^2_{0,E_h}} \lesssim \frac{\abs{\projE \backslash E_h} }{\abs{E_h}} R^2
    \lesssim h_{\hat{E}}^s  R^2,
\end{align*}
having used Lemmas \ref{lemma::AreaError} and \ref{lemma::MeshScalingV2} to obtain the last inequality.
Substitution into \eqref{eq::Rpt1} implies that there exist constants $C_1, C_2 >0$ such that
\begin{align*}
    R^2 &\leq C_1 + C_2 R^2 h_{\hat{E}}^s,
\end{align*}
therefore assuming $h_{\hat{E}}$ is such that $h_{\hat{E}}^s \leq  C_2^{-1} \slash 2 $
results in $R^2 \leq 2 C_1$, implying that the ratio between norms $R$ is uniformly bounded for $h$ small enough, as required.
\end{proof}
\end{lemma}

\subsection{Quadrature error}
Herein we analyse the error due to the use of projected elements within the physical IsoVEM. 
 
\begin{lemma}[Projected element integration error]\label{lemma::QuadratureError} 
Let $\RefMesh$ satisfy Assumption~\ref{assumption::ShapeRegular}, $\F_h$ satisfy Assumption~\ref{assumption::F_hAccuracy}, and $h$ be small enough.
For all $E_h\in\mesh$, $g\in L^\infty(E_h\cup \projE)$, and $p \in \pp_k(E_h)$, there holds 
\begin{align}
    \abs{\int_{E_h} {g} p\ d\PhysCoordinate - \int_{\projE} {g} p\ d\mathbf{x}} &\lesssim h^{s+1}_{\hat{E}} \norm{{g}}_{0,\infty,E_h \triangle \projE}\norm{p}_{0,E_h}, \label{eq:proj_ele_err} 
\end{align}
with $ s \in \{0,1,...,l\}$.
The hidden constant depends on the shape regularity of $\RefMesh$, on $k,l$, on 
the $W^{s+1}_\infty(\hat{E})$-norm $\F$, and on $j_0$. 
\begin{proof}
We decompose the integral difference into

\begin{align}
    \abs{\int_{E_h} {g} p\ d\PhysCoordinate - \int_{\projE} {g} p\ d\mathbf{x}} &= \abs{\int_{E_h \backslash \projE} {g} p\ d\PhysCoordinate - \int_{\projE \backslash E_h} {g} p\ d\PhysCoordinate}\nonumber\\
    &\leq  \abs{E_h \triangle \projE} \norm{{g}}_{0,\infty,E_h \triangle \projE} \norm{p}_{0,\infty,E_h \triangle \projE}. \label{eq::QuadPt1}
\end{align}
Next, we  bound  the $L^\infty$-norm of $p$ using Lemma \ref{lemma::PolyExten}, yielding 
\begin{align*}
    \norm{p}_{0,\infty,E_h \triangle \projE} &\leq \norm{p}_{0,\infty,E_h \backslash \projE} + \norm{p}_{0,\infty,\projE \backslash E_h} \leq \norm{p}_{0,\infty,E_h} + \norm{p}_{0,\infty,\projE}\lesssim \norm{p}_{0,\infty,E_h}.
\end{align*}
Substituting this bound into Equation \eqref{eq::QuadPt1}, applying Lemma \ref{lemma::AreaError} to bound the measure of $E_h \triangle \projE$, hiding norms of $\F$ 
within the hidden constant, and using once again an inverse estimate, results in
\begin{align*}
    \abs{\int_{E_h} {g} p\ d\PhysCoordinate - \int_{\projE} \tilde{g} p\ d\mathbf{x}} &\lesssim h_{\hat{E}}^{s+2} 
    \norm{{g}}_{0,\infty,E_h \triangle \projE}
    \norm{p}_{0,\infty,E_h}
    \\
    &\lesssim \abs{E_h}^{-1/2} h_{\hat{E}}^{s+2} \norm{{g}}_{0,\infty,E_h \triangle \projE}\norm{p}_{0,E_h},
\end{align*}
from which applying Lemma \ref{lemma::MeshScalingV2} concludes the proof.
\end{proof}
\end{lemma}

\subsection{Projection error estimates}
We are now in a position to assess the accuracy with respect to the true $L^2$-projection of the approximate projection operators $\Pi^0_h$ and $\Pi^1_h$ given by Definitions~\ref{def::PiStarValue} and~\ref{def::PiStarValueGrad}, respectively.

\newcommand{\Cstar}[1]{C^*_{#1}}

\begin{theorem}\label{theorem::PiStarOpBound} 
Let $\RefMesh$ satisfy Assumption~\ref{assumption::ShapeRegular}, $\F_h$ satisfy Assumption~\ref{assumption::F_hAccuracy}, and $h$ be small enough.
The operator $\Pi_{h,k}^0$ given by Definition~\ref{def::PiStarValue} 
satisfies, for each $E_h \in \mesh$ and $v_h \in V_k(E_h)$,
\begin{align*}
    \norm{\Pi^0_k v_h - \Pi_{h,k}^0 v_h}_{0,E_h} &\lesssim h_{\hat{E}}^{s}  \norm{v_h}_{0,E_h},\\
     \norm{\Pi_{h,k}^0 v_h}_{0,E_h} &\lesssim \norm{v_h}_{0,E_h},
\end{align*}
for  $s \in \{0,1,..., l \}$.
The hidden constants depend on the shape regularity of $\RefMesh$, on $k$, $l$, $j_0$, and on the $W^{s+1}_\infty(\hat{E})$-norm of $\F$.

\begin{proof}
\renewcommand{\a}{\mathfrak{a}}
\renewcommand{\b}{\mathbf{b}}
\newcommand{\M}{\mathbf{M}}
\newcommand{\load}{\mathbf{b}}
\newcommand{\basis}{\boldsymbol{\phi}}
\renewcommand{\d}{\mathbf{d}}
\renewcommand{\t}{\mathbf{t}}
\renewcommand{\c}{\mathbf{c}}

In this proof we adopt the notation used in \cite{DaVeiga2014} in defining the matrix equations to compute the $L^2$ projection. We consider the error for some $v_h \in V_h$ as a vector-matrix product
\begin{equation*}
    \norm{\Pi^0_k v_h - \Pi_{h,k}^0 v_h}_{0,E_h}^2 = \left(\t - \t_h\right)^{T} \M \left(\t - \t_h \right), 
\end{equation*}
where we define the mass matrix $\mathbf{M}$ on $E_h$ by
\begin{align*}
    M_{\alpha,\beta} &= \int_{E_h} m_\alpha m_\beta\ d\PhysCoordinate,
\end{align*}
for $m_\alpha, m_\beta \in \mathcal{M}_k(E_h)$ (see Definition \ref{definition::dofs}) and the vectors of coefficients of $\Pi^0$ and $\Pi_h^0$ are given by $\t$ and $\t_h$ respectively. By using the $l_2$ norm $\norm{\cdot}_2$ and the matrix norm induced by $\norm{\cdot}_2$ and analysing individual matrix terms we bound the error noting that the mass matrix terms scale like $\abs{E_h}$
\begin{align}
    \norm{\Pi^0_k v_h - \Pi_{h,k}^0 v_h}_{0,E_h}^2 &\leq \norm{\M}_2  \norm{\t - \t_h}_2^2 \lesssim h_{E_h}^2 \norm{\t - \t_h}_2^2, \label{eq::OriginalError}
\end{align}
and it remains to bound the error in the projection coefficients. The $\Pi^0_k$ projection is determined by $ \t = \M^{-1} \c$ where
\begin{equation}
    c_\alpha = 
    \begin{cases}
    \abs{E_h} \text{dof}_\alpha (v_h) \qquad &m_\alpha \in \mathcal{M}_{k-2}(E_h), \nonumber \\
    \int_{E_h} \Pi^*_k\ v_h\ m_\alpha\ d\PhysCoordinate \qquad &m_\alpha \in \mathcal{M}_{k}(E_h) \backslash \mathcal{M}_{k-2}(E_h),
    \end{cases}
\end{equation}
while the $\Pi_h^0$ projection is determined by $\t_h = \M_h^{-1} \c_h$ where
\begin{equation*}
    c_{\alpha,h} = \begin{cases}
    \abs{\projE} \text{dof}_\alpha (v_h) \qquad &m_\alpha \in \mathcal{M}_{k-2}(E_h), \\
    \int_{\projE} \Pi^*_k v_h\ m_\alpha\ d\PhysCoordinate \qquad &m_\alpha \in \mathcal{M}_{k}(E_h) \backslash \mathcal{M}_{k-2}(E_h).
    \end{cases}
\end{equation*}
We then consider the difference of these coefficients, noting that
\begin{equation*}
    \t - \t_h = \left(\M^{-1} - \M_h^{-1}  \right) \c\ +\ \M_h^{-1} \left( \c - \c_h \right).
\end{equation*}
Passing to  norms and applying the triangle and matrix norm inequalities results in
\begin{align}\label{eq::projproof1}
    \norm{\t - \t_h}_2 \leq \norm{\M^{-1} - \M_h^{-1}}_2 \norm{\c}_2 + \norm{\M_h^{-1}}_2 \norm{ \c - \c_h}_2.
\end{align}
The first term can be expanded by
\begin{equation*}
    \M^{-1} - \M_h^{-1} = \M_h^{-1} \left( \M_h - \M \right) \M^{-1},
\end{equation*}
and hence, noting that the matrix inverses norms  $\norm{\M^{-1}}_2$ and $\norm{\M_h^{-1}}_2$ scale like  $h_{E_h}^{-2}$ and applying Lemma \ref{lemma::QuadratureError} with $g\equiv 1$ and $p=m_\alpha m_\beta$, and inverse estimates we obtain
\begin{align*}
\norm{\M^{-1} - \M_h^{-1}}_2 &\leq \norm{\M_h^{-1}}_2 \norm{ \M_h - \M }_2 \norm{\M^{-1}}_2 \lesssim h_{E_h}^{-4} \norm{ \M_h - \M }_2 \\
&\lesssim h_{E_h}^{-4} \max_{\alpha,\beta} \abs{\int_{E_h} m_\alpha m_\beta\ d\PhysCoordinate - \int_{\projE} m_\alpha m_\beta\ d\PhysCoordinate} \\
&\lesssim h_{E_h}^{-4} h_{\hat{E}}^s\norm{m_\alpha}_{0,E_h} \norm{m_\beta}_{0,E_h} 
\lesssim h_{E_h}^{-2} h^s_{\hat{E}},
\end{align*}
given that $\norm{m_\alpha}_{0,E_h} \sim h_{E_h}$ for all $m_\alpha$. Substituting this back  in~\eqref{eq::projproof1} gives
\begin{align}
    \norm{\t - \t_h}_2 \lesssim h_{E_h}^{-2} h^s_{\hat{E}} \norm{\c}_2 + h_{E_h}^{-2} \norm{ \c - \c_h}_2.\label{eq::PiStarValueBoundPt1}
\end{align}
The difference of $\c - \c_h$ requires a case-by-case analysis.

\paragraph{Case I: $m_\alpha \in \mathcal{M}_{k-2}(E_h)$.}
Taking the difference and applying Lemma \ref{lemma::AreaError} gives
\begin{align}
    \abs{c_{\alpha} - c_{\alpha,h}} &\leq \abs{ \text{dof}_{\alpha}(v_h) } \abs{ \abs{E_h} - \abs{\projE} } \nonumber \\ 
    &\lesssim \abs{ \text{dof}_{\alpha}(v_h) }  h_{\hat{E}}^{s+2} \norm{\F}^2_{s+1,\infty,\hat{E}}
    \lesssim h_{E_h} h_{\hat{E}}^{s} \norm{v_h}_{0,E_h},\label{eq::calphah}
\end{align}
having used in the last bound the fact that  the scaled monomials scale like one, that is
$\text{dof}_\alpha(v_h) \sim h_{E_h}^{-1} \norm{v_h}_{0,E_h}$, cf. Definition~\ref{definition::dofs} and \cite{DaVeiga2014}.

\paragraph{Case II: $m_\alpha \in \mathcal{M}_{k}(E_h) \backslash \mathcal{M}_{k-2}(E_h)$.} Using  the definition of the VEM space \eqref{eq::CurvedVEMSpace}, and then applying Lemma~\ref{lemma::QuadratureError} with $g\equiv 1$ and $p=\Pi^*_k v_h\, m_\alpha$, inverse estimates, and Lemma \ref{lemma::MeshScalingV2}, we get
\begin{align*}
    \abs{c_{\alpha} - c_{\alpha,h}} &\lesssim h_{\hat{E}}^s \norm{\Pi^*_k v_h}_{0,E_h} \norm{m_\alpha}_{0,E_h}\lesssim h_{E_h} h_{\hat{E}}^{s} \norm{v_h}_{0,E_h}.
\end{align*}

Having bounded both cases, we can now bound $\c-\c_h$ via
\begin{align*}
    \norm{\c-\c_h}_2 &\lesssim h_{\hat{E}}^{s+1} \norm{v_h}_{0,E_h}.\label{eq::PiStarValueBoundPt2}
\end{align*}
The norm of $\c$ can also be bounded by a similar argument, yielding
\begin{align*}
    \norm{\c}_2 \lesssim h_{E_h} \norm{v_h}_{0,E_h}, 
\end{align*}
and, using these bounds and the scaling of $\M_h^{-1}$ in~\eqref{eq::PiStarValueBoundPt1}, we arrive to
\begin{align*}
    \norm{\t-\t_h}_2 &\lesssim h_{\hat{E}}^s h_{E_h}^{-1} \norm{v_h}_{0,E_h}.
\end{align*}
To conclude the proof, we substitute this bound back into the original error, Equation \eqref{eq::OriginalError}, to get
\begin{align}
    \norm{\Pi^0_k v_h - \Pi_{h,k}^0 v_h}_{0,E_h}^2 &\lesssim h_{E_h}^2 \norm{\t - \t_h}_2^2 \lesssim h_{E_h}^2  \left( h_{\hat{E}}^s h_{E_h}^{-1} \norm{v_h}_{0,E_h} \right)^2 \lesssim h_{\hat{E}}^{2s} \norm{v_h}_{0,E_h}^2.
\end{align}
Taking the square root completes the proof. The stability bound is a consequence of setting $s=0$ and applying the triangle inequality.
\end{proof}
\end{theorem}

\begin{theorem}\label{theorem::PiStarValueGrad} 
Let $\RefMesh$ satisfy Assumption~\ref{assumption::ShapeRegular}, $\F_h$ satisfy Assumption~\ref{assumption::F_hAccuracy}, and $h$ be small enough.
The operator $\Pi_{h,k-1}^1$ given by Definition~\ref{def::PiStarValueGrad} satisfies, for each $E_h \in \mesh$ and $v_h \in V_k(E_h)$,

\begin{align*}
    \norm{\Pi^1_{k-1} v_h - \Pi_{h,k-1}^1 v_h}_{0,E_h} &\lesssim h_{\hat{E}}^{s}  \norm{v_h}_{1,E_h},\\
    \norm{\Pi_{h,k-1}^1 v_h}_{0,E_h} &\lesssim \norm{v_h}_{1,E_h},
\end{align*}
with $s \in \{0,1,..., l \}$.
The hidden constants depend on the shape regularity of $\RefMesh$, on $k$, $l$, $j_0$, and the $W^{s+1}_\infty(\hat{E})$-norm of $\F$.

\begin{proof}
\renewcommand{\a}{\mathfrak{a}}
\renewcommand{\b}{\mathbf{b}}
\newcommand{\M}{\mathbf{M}}
\newcommand{\load}{\mathbf{b}}
\newcommand{\basis}{\boldsymbol{\phi}}
\renewcommand{\d}{\mathbf{d}}
\renewcommand{\t}{\mathbf{t}}
\renewcommand{\c}{\mathbf{c}}

The proof of the estimate follows the same arguments as Theorem \ref{theorem::PiStarOpBound}. We have, 

\begin{align}
    \norm{\Pi^1_{k-1} v_h - \Pi_{h,k-1}^1 v_h}_{0,E_h}^2 &\leq \norm{\M}_2 \norm{\t-\t_h}_2^2 \lesssim h_{E_h}^2 \norm{\t - \t_h}_2^2, \label{eq::GradProjPt1}
\end{align}
where $\M$ and the corresponding approximation $\M_h$ are the vector equivalents of $\M$ and $\M_h$ defined in the proof of Theorem \ref{theorem::PiStarOpBound}.
We recall that $\Pi^1_{k-1} v_h$  is determined by $\t = \M^{-1} \c$
where
\begin{align}
    c_\alpha = -\int_{E_h} v_h \nabla \cdot \mathbf{m}_\alpha\ d\PhysCoordinate + \int_{\partial E_h} v_h \mathbf{m}_\alpha \cdot \mathbf{n}\ dS, \qquad \mathbf{m}_\alpha \in \left[ \mathcal{M}_{k-1}(E_h) \right]^2. \label{eq::GradProjC} 
\end{align}
Similarly, $\Pi_{h,k-1}^1 v_h$ is determined by $\t_h = \M_h^{-1} \c_h$
where 
\begin{align*}
    c_{\alpha,h} = -  \int_{\projE} \Pi^0_{h,k} v_h \ \nabla \cdot \mathbf{m}_\alpha\ d\PhysCoordinate + \int_{\partial E_h} v_h \mathbf{m}_\alpha \cdot \mathbf{n}\ dS, \qquad \mathbf{m}_\alpha \in \left[ \mathcal{M}_{k-1}(E_h) \right]^2.
\end{align*}
Following the proof of Theorem \ref{theorem::PiStarOpBound} up to Equation \eqref{eq::PiStarValueBoundPt1} gives
\begin{align}
    \norm{\t-\t_h}_2 \lesssim h_{E_h}^{-2} h_{\hat{E}}^s \norm{\c}_2 + h_{E_h}^{-2} \norm{\c-\c_h}_2. \label{eq::GradProjPt2} 
\end{align}
Further, applying integration by parts to Equation \eqref{eq::GradProjC} gives
\begin{align*}
    c_\alpha = \int_{E_h} \mathbf{m}_\alpha \cdot \nabla v_h \ d\PhysCoordinate,
\end{align*}
from which we bound $\c$ using the Cauchy-Schwarz inequality and the scaling of $\mathbf{m}_\alpha$, yielding
\begin{align}
    \norm{\c}_2 \lesssim h_{E_h} \norm{\nabla v_h}_{0,E_h}. \label{eq::GradProjPt2.5}
\end{align}
It remains to estimate $\norm{\mathbf{c} - \mathbf{c}_h}_2$. 
This error term is given as
\begin{align*}
     c_{\alpha} - c_{\alpha,h} = \int_{E_h} v_h \nabla \cdot \mathbf{m}_\alpha\ d\PhysCoordinate - \int_{\projE} \Pi^0_{h,k} v_h\ \nabla \cdot \mathbf{m}_\alpha \ d\PhysCoordinate.
\end{align*}
Since $\nabla \cdot \mathbf{m}_\alpha \in [\pp_{k-2}(E_h)]^2$, we can evaluate this error term using only the internal dofs (see Definition \ref{def::PiStarValue}): there exists a set $\{ A_\beta \}_{\abs{\beta} \leq k-2} \subset \rr$ such that
\begin{align}
    c_{\alpha,h} - c_{\alpha} = \sum_{\abs{\beta} \leq k-2} A_\beta \text{dof}_\beta (v_h) (\abs{E_h} - \abs{E_h^\pi}). \label{eq::GradProjPt3}
\end{align}
The difference in volume is bounded as done in~\eqref{eq::calphah} by applying Lemma \ref{lemma::AreaError}, resulting in
\begin{align*}
    \abs{\abs{\projE} - \abs{E_h}} \lesssim  h_{\hat{E}}^{s} \abs{E_h}.
\end{align*}
We have, for all $|\beta|\le k-2$, that
\begin{align*}
    \text{dof}_\beta (v_h) =
    \frac{1}{\abs{E_h}} \int_{E_h} v_h \nabla \cdot \mathbf{m}_\beta\ d\PhysCoordinate \lesssim h_{E_h}^{-1} \norm{v_h}_{0,E_h}.
\end{align*}
Inserting these bounds into Equation \eqref{eq::GradProjPt3} and bounding terms gives
\begin{align*}
    \norm{\c - \c_{h}}_2 &\lesssim h_{\hat{E}}^{s} h_{E_h} \norm{v_h}_{0,E_h}.
\end{align*}
Next we insert this bound into Equation \eqref{eq::GradProjPt2} along with Equation \eqref{eq::GradProjPt2.5}, to obtain
\begin{align*}
    \norm{\t-\t_h}_2 &\lesssim h_{\hat{E}}^s h_{E_h}^{-1} \left( \norm{\nabla v_h}_{0,E_h} +  \norm{v_h}_{0,E_h} \right) \lesssim h_{\hat{E}}^s h_{E_h}^{-1} \norm{v_h}_{1,E_h}.
\end{align*}
Finally, we have from Equation \eqref{eq::GradProjPt1} that
\begin{align}
    \norm{\Pi^1_{h,k-1} v_h - \Pi^1_{k-1} v_h}_{0,E_h}^2 &\lesssim h_{E_h}^2 \left( h_{\hat{E}}^s h_{E_h}^{-1} \norm{v_h}_{1,E_h} \right)^2 \lesssim h_{\hat{E}}^{2s}   \norm{v_h}_{1,E_h}^2.
\end{align}
Taking the square root completes the proof. The stability bound is a consequence of setting $s=0$ and applying the triangle inequality.
\end{proof}
\end{theorem}

\subsection{Well-posedness}\label{sec::MethodII::LaxMilgram}
The use of \emph{overlapping} projected elements in the assembly of the physical IsoVEM introduces a new form of inconsistency. As a consequence, the proof of coercivity of the physical IsoVEM bilinear form $\mathcal{A}_h$ given in~\eqref{eq::bilinear_phy} is not-trivial and requires a careful analysis. 
To this end, we introduce the discrete bilinear form 
$\pertB_h:V_h^k\times V_h^k\rightarrow\rr$ defined by summing up the elemental contributions
\begin{align}
    \pertB_h^{E_h}(u,v) = 
    &
    \int_{E_h} \Big(\Exten{a} \Pi^1_{k-1} u \cdot \Pi^1_{k-1} v + \frac{1}{2} \Exten{\mathbf{b}} \cdot ( \Pi^0_k v \Pi^1_{k-1} u - \Pi^0_k u \Pi^1_{k-1} v) + \Exten{\mu} \Pi^0_k u\ \Pi^0_k v\Big) \ d\mathbf{x} \nonumber \\
    +& S^{E_h}(u-\Pi^0_k u,v - \Pi^0_k v),
    \label{eq::B_hDef}
\end{align}
where the stabilisation term is assumed to be identical to the choice of stabilisation in Section \ref{sec::MethodIIFormulation}.

\begin{lemma}\label{lemma::IDKYet}
Let $s=\min \{k,l\}$. Then, for all $u_h,v_h \in V_h^k$,  it holds that
\begin{equation*}
    \sum_{E_h \in \VirtualMesh}\abs{\pertB_h^{E_h}(u_h,v_h) - \Ah^{E_h}(u_h,v_h)} \lesssim  h^s \norm{u_h}_{1,\Omega_h} \norm{v_h}_{1,\Omega_h},
\end{equation*}
where the hidden constant depends on  the shape regularity of $\RefMesh$, on $k$, $l$,  $j_0$, the $W^{s+1}_\infty(\RefDomain)$-norm of $\F$, and the $L^\infty(\Omega)$-norm of $a$, $\mathbf{b}$, and $\mu$.
\begin{proof}
Noting that there is no contribution to $\pertB_h- \Ah$ from the stabilisation term,  we split the difference via
\begin{align*}
    \pertB_h^{E_h}(u_h,v_h) - \mathcal{A}_h^{E_h}(u_h,v_h) =(T_q^{E_h} + T_a^{E_h} + T_b^{E_h} + T_c^{E_h}),
\end{align*}
where through adding and subtracting terms we have,
\begin{align*}
    T_q^{E_h} &= \int_{E_h} \Exten{a} \Pi^1_{h,k-1} u_h \cdot \Pi^1_{h,k-1} v_h\ d\PhysCoordinate - \int_{\projE} \Exten{a} \Pi^1_{h,k-1} u_h \cdot \Pi^1_{h,k-1} v_h\ d\PhysCoordinate \\
    &+\frac{1}{2} \Bigg( \int_{E_h}  \Exten{\mathbf{b}} \cdot (\Pi_{h,k}^0 v_h \Pi^1_{h,k-1} u_h - \Pi_{h,k}^0 u_h \Pi^1_{h,k-1} v_h)\ d\PhysCoordinate 
    \\
    &\hspace{1cm}-  \int_{\projE}  \Exten{\mathbf{b}} \cdot (\Pi_{h,k}^0 v_h \Pi^1_{h,k-1} u_h - \Pi_{h,k}^0 u_h \Pi^1_{h,k-1} v_h)\ d\PhysCoordinate\Bigg) \\
    &+ \int_{E_h} \Exten{\mu} \Pi_{h,k}^0 u_h \Pi_{h,k}^0 v_h\ d\PhysCoordinate - \int_{\projE} \Exten{\mu} \Pi_{h,k}^0 u_h \Pi_{h,k}^0 v_h\ d\PhysCoordinate, \\
    T_a^{E_h} &= \int_{E_h}  \Exten{a} (\Pi^1_{k-1} - \Pi^1_{h,k-1}) u_h \cdot \Pi^1_{k-1} v_h + \Exten{a} \Pi^1_{h,k-1} u_h \cdot (\Pi^1_{k-1} - \Pi^1_{h,k-1}) v_h\ d\PhysCoordinate, \\
    T_b^{E_h} &= \frac{1}{2} \int_{E_h} \Exten{\mathbf{b}} \cdot \big( \Pi_{k}^0 v_h \Pi^1_{k-1} u_h - \Pi_{k}^0 u_h \Pi^1_{k-1} v_h \\
    &\hspace{4cm}+ \Pi_{h,k}^0 v_h \Pi^1_{h,k-1} u_h - \Pi_{h,k}^0 u_h \Pi^1_{h,k-1} v_h \big) d\PhysCoordinate, \\
    T_c^{E_h} &= \int_{E_h} \left( \Exten{\mu} (\Pi^0_k-\Pi_{h,k}^0)u_h \Pi^0_k v_h + \Exten{\mu} \Pi^0_{h,k} u_h (\Pi^0_k-\Pi_{h,k}^0) v_h\ \right) d\PhysCoordinate.
\end{align*}
The local quadrature error term $T_q^{E_h}$ is bounded by applying Lemma \ref{lemma::QuadratureError} to give
\begin{align*}
    \abs{T_q^{E_h}} &\lesssim h^s_{\hat{E}} \Big( \norm{\Exten{a}}_{0,\infty,E_h\triangle \projE}\norm{\Pi^1_{h,k-1} u_h}_{0,E_h} \norm{\Pi^1_{h,k-1} v_h}_{0,E_h}\\
    &\hspace{2cm}+\frac{1}{2}\norm{\mathbf{\Exten{b}}}_{0,\infty,E_h\triangle \projE} \norm{\Pi^1_{h,k-1} u_h}_{0,E_h} \norm{\Pi_{h,k}^0 v_h}_{0,E_h}\\
    &\hspace{2cm}+\frac{1}{2}\norm{\mathbf{\Exten{b}}}_{0,\infty,E_h\triangle \projE}  \norm{\Pi^1_{h,k-1} v_h}_{0,E_h} \norm{\Pi_{h,k}^0 u_h}_{0,E_h} \\
    &\hspace{2cm}+ \norm{\Exten{\mu}}_{0,\infty,E_h\triangle \projE} \norm{\Pi_{h,k}^0 u_h}_{0,E_h} \norm{\Pi_{h,k}^0 v_h}_{0,E_h}\nonumber \Big).
\end{align*}
Then, using the stability bounds in Theorems \ref{theorem::PiStarOpBound} and \ref{theorem::PiStarValueGrad}, provides
\begin{align*}
    \abs{T_q^{E_h}}  &\lesssim h^s_{\hat{E}} \max \left( \norm{\Exten{a}}_{0,\infty,E_h\triangle \projE}, \norm{\mathbf{\Exten{b}}}_{0,\infty,E_h\triangle \projE}, \norm{\Exten{\mu}}_{0,\infty,E_h\triangle \projE} \right) \norm{ u_h}_{1,E_h} \norm{ v_h}_{1,E_h} 
\end{align*}
Considering $T_a^{E_h}$, 
we use Theorem \ref{theorem::PiStarValueGrad} to get
\begin{align*}
    \abs{T_a^{E_h}} &\leq \norm{\Exten{a}}_{0,\infty,E_h} \left(\norm{\Pi^1_{k-1} u_h - \Pi^1_{h,k-1} u_h}_{0,E_h} \norm{\Pi^1_{k-1} v_h}_{0,E_h}\right.\\
    &
    \qquad \qquad\qquad + \left. \norm{\Pi^1_{h,k-1} u_h}_{0,E_h} \norm{\Pi^1_{k-1} v_h - \Pi^1_{h,k-1} v_h}_{0,E_h}\right) \\
    &\lesssim h_{\hat{E}}^s \norm{\Exten{a}}_{0,\infty,E_h}\norm{ u_h}_{1,E_h} \norm{ v_h}_{1,E_h}.
\end{align*}
 The terms $T_b^{E_h}$ and $T_c^{E_h}$ can be similarly bounded. 
We conclude the proof by summing up the above bounds over all elements in the mesh and applying Stein's extension Theorem~\ref{theorem::SteinExtension} to bound the infinity norms of $\Exten{a}$, $\mathbf{\Exten{b}}$, and $\Exten{\mu}$ from $\Omega_h$ back to $\Omega$.
\end{proof}
\end{lemma}

\begin{theorem}[Physical IsoVEM well-posedness]\label{theorem::Method2WellPosed} 
Let $\RefMesh$ satisfy Assumption~\ref{assumption::ShapeRegular}, $\F_h$ satisfy Assumption~\ref{assumption::F_hAccuracy}, and $h$ be small enough. Then, the physical IsoVEM formulation~\eqref{eq::physicalIsoVEM} with stabilisation~\eqref{eq::MethodIIStab} obeying Assumption~\ref{assumption::Stab::Stab} possesses a unique solution.
\begin{proof}
 We verify once more the assumptions of the Lax-Milgram lemma, starting with the continuity of $\mathcal{A}_h$. Application of the triangle, \holder \ and Cauchy-Schwarz inequality we have
\begin{align*}
    \abs{\mathcal{A}_h(u_h,v_h)} &\leq \sum_{E_h \in \mesh} \Big( \norm{\Exten{a}}_{0,\infty,E_h}\norm{\Pi^1_{h,k-1} u_h}_{0,E_h} \norm{\Pi^1_{h,k-1} v_h}_{0,E_h} \\
    &\hspace{2cm} +\frac{1}{2} \norm{\Exten{\mathbf{b}}}_{0,\infty,E_h} \norm{\Pi^1_{h,k-1} u_h}_{0,E_h} \norm{\Pi_{h,k}^0 v_h}_{0,E_h} \\
    &\hspace{2cm}+ \frac{1}{2} \norm{\Exten{\mathbf{b}}}_{0,\infty,E_h} \norm{\Pi^1_{h,k-1} v_h}_{0,E_h} \norm{\Pi_{h,k}^0 u_h}_{0,E_h} \\ 
    &\hspace{2cm}+ \norm{\Exten{\mu}}_{0,\infty,E_h} \norm{\Pi_{h,k}^0 u_h}_{0,E_h} \norm{\Pi_{h,k}^0 v_h}_{0,E_h} + S^{E_h}(u_h,v_h)\Big).
\end{align*}
Then we apply the stability estimates from Theorems \ref{theorem::PiStarOpBound} and \ref{theorem::PiStarValueGrad}, Assumption~\ref{assumption::Stab::Stab} on the VEM stabilisation term, and finally Cauchy-Schwarz to get
\begin{align*}
    \abs{\mathcal{A}_h(u_h,v_h)} &\lesssim \max\left\{ \norm{\Exten{a}}_{0,\infty,\Omega_h}, \norm{\Exten{\mathbf{b}}}_{0,\infty,\Omega_h}, \norm{\Exten{\mu}}_{0,\infty,\Omega_h} \right\} \norm{u_h}_{1,\Omega_h} \norm{ v_h}_{1,\Omega_h},
\end{align*}
which can be further bounded in terms of norms of the data over $\Omega$ thanks to Stein's extension Theorem~\ref{theorem::SteinExtension}, thus completing the continuity proof.
The continuity of $l_h$ can be similarly proven.

To prove the coercivity of $\mathcal{A}_h$, we consider the bilinear form $\pertB_h^{E_h}$ defined in Equation \eqref{eq::B_hDef} and consider the coercivity constant $C_{1} >0$ given in \cite{ConNonConVEM} and the continuity constant $C_2$ of Lemma \ref{lemma::IDKYet} to get
\begin{align*}
    C_{1} \norm{v_h}_{1,\Omega_h}^2 &\leq \pertB_h(v_h,v_h) \leq \Ah(v_h,v_h) + \sum_{E_h \in \VirtualMesh} \abs{\Ah^{E_h}(v_h,v_h) - \pertB_h^{E_h}(v_h,v_h)} \\
    &\leq \Ah(v_h,v_h) + C_{2} h^s \norm{v_h}^2_{1,\Omega_h}.
\end{align*}
Coercivity now follows  for sufficiently small $h$ such that $C_1-C_2h^s > 0$.
\end{proof}

\end{theorem}

\subsection{Error estimate}\label{sec::MethodII::H1Error}
Before embarking on the error analysis for the physical IsoVEM~\eqref{eq::physicalIsoVEM}, we note that the discrete solution $u_h \in V_{h,0}^k$ is defined over the virtual domain $\Omega_h$ while the exact solution $u\in H^1_0(\Omega)$ is defined over the true physical domain $\Omega$, hence they are \emph{not} directly comparable. Instead, we shall compare $u_h$ with the restriction on $\Omega_h$ of the extension $\Exten{u}$ given by Stein's Theorem~\ref{theorem::SteinExtension}.
Similarly, we introduce a new bilinear form $\tilde{\mathcal{A}}$  to replace the  form $\mathcal{A}$ used in the continuum problem~\eqref{eq::EllipticModelPDE}. This is defined over  $H^1(\Omega_h)$ by replacing   $\Omega, {a}, {\mathbf{b}},$ and ${\mu}$ with, respectively,  $\Omega_h, \tilde{a}, \tilde{\mathbf{b}},$ and $\tilde{\mu}$ in~\eqref{eq::modelAB}-\eqref{eq::modelCl}. Finally, following Ciarlet~\cite{CiarletElliptic}, we introduce the extension $f_*$ of the forcing function $f$ given a.e. in $\Omega_h$ by
\begin{equation}
    f_* := -\nabla \cdot (\Exten{a} \nabla \Exten{u}) + \Exten{\mathbf{b}} \cdot \nabla \Exten{u} + \Exten{c} \Exten{u}, \label{eq::fStarDef}
\end{equation}
assuming that $f_*\in L^2(\Omega_h)$. We thus notice that the extended solution $\tilde{u}$ satisfies
\begin{align}\label{eq:starform}
\tilde{\mathcal{A}}(\Exten{u},v) = l_*(v) \qquad\forall v\in H^1_0(\Omega_h),
\end{align}
where
the linear operator $l_*$ is given by
\begin{align*}
    l_*(v) = \int_{\Omega_h} f_* v\ d\PhysCoordinate
    \qquad \forall v\in H^1(\Omega_h).
\end{align*}
These extensions allow derivation of the  Strang-type bound analogue to that of Theorem~\ref{theorem::StrangBound1}, cf. also~\cite{CiarletElliptic},  and thus prove the following result.

\begin{theorem}[Physical IsoVEM Strang-type bound]\label{theorem::Method2Strang}
Let $\RefMesh$ satisfy Assumption~\ref{assumption::ShapeRegular}, $\F_h$ satisfy Assumption~\ref{assumption::F_hAccuracy}, and $h$ be small enough. 
For the solution  $u_h \in V_{h,0}^k$  of the physical IsoVEM~\eqref{eq::physicalIsoVEM} it holds that
\begin{align}
    \norm{\Exten{u}-u_h}_{1,\Omega_h} &\lesssim \inf_{v_h \in V_{h,0}^k} \norm{\Exten{u} - v_h}_{1,\Omega_h}   + \inf_{p \in \mathbb{P}_k(\VirtualMesh)} \norm{\Exten{u} - p}_{1,\Omega_h} \nonumber \\
    &+ \sup_{w_h \in V_{h,0}^k \backslash \{ 0 \}}\frac{\left| l_h(w_h) - l_*(w_h) \right|}{\norm{w_h}_{1,\Omega_h}} \nonumber \\
    &+ \inf_{p \in \mathbb{P}_k(\VirtualMesh)}
    \sup_{w_h \in V_{h,0}^k\backslash \{ 0 \}}  \frac{\sum_{E_h \in \VirtualMesh} \left| \tilde{\mathcal{A}}^{E_h}(p,w_h) - \pertB_h^{E_h}(p,w_h) \right|}{\norm{w_h}_{1,\Omega_h}} \nonumber \\&
    + \inf_{p \in \mathbb{P}_k(\VirtualMesh)}   \sup_{w_h \in V_{h,0}^k \backslash \{ 0 \}} \frac{\sum_{E_h \in \VirtualMesh}\left| \pertB_h^{E_h}(p,w_h) - \mathcal{A}_h^{E_h}(p,w_h) \right|}{\norm{w_h}_{1,\Omega_h}},
    \label{eq::StrangBound2}
\end{align}
where the hidden constant depends on the shape regularity of $\RefMesh$, on $k$ and $l$, on the lower bounds $a_0$, $\mu_0$, $j_0$ and the $L^\infty(\Omega)$-norms of $a$, $\mathbf{b}$, and $\mu$, and on the $W^{s+1}_\infty(\hat{\Omega})$ norm of $\F$.
\begin{proof}
We have once again following \cite{ConNonConVEM}  the basic Strang-type bound  analogue to~\eqref{theorem::StrangBound1}, which in the present setting reads
\begin{align*}
    \norm{\Exten{u}-u_h}_{1,\Omega_h} &\lesssim \inf_{v_h \in V_{h,0}^k} \norm{\Exten{u} - v_h}_{1,\Omega_h}   + \inf_{p \in \mathbb{P}_k(\VirtualMesh)} \norm{\Exten{u} - p}_{1,\Omega_h} \\
     &\hspace{-.4cm}+ \sup_{w_h \in V_{h,0}^k \backslash \{ 0 \}}\frac{\left| l_h(w_h) - \tilde{\mathcal{A}}(\Exten{u},w_h) \right|}{\norm{w_h}_{1,\Omega_h}}\\
    &\hspace{-.4cm}+ \inf_{p \in \mathbb{P}_k(\VirtualMesh)}  \sup_{w_h \in V_k(E_h) \backslash \{ 0 \}} \frac{\sum_{E_h \in \VirtualMesh} \left| \tilde{\mathcal{A}}^{E_h}(p,w_h) - \mathcal{A}_h^{E_h}(p,w_h) \right|}{\norm{w_h}_{1,\Omega_h}}. 
\end{align*}
Inserting~\eqref{eq:starform} in the third term and adding and subtracting $\pertB_h^{E_h}(p,w_h)$ in the fourth we conclude.
\end{proof}
\end{theorem}

We assess the consistency error terms in the above Strang-type bound. The fifth term in~\eqref{eq::StrangBound2}, measuring the bilinear form inconsistency due to projection, was already treated in Lemma \ref{lemma::IDKYet}. The fourth term, measuring the consistency between the exact and virtual bilinear forms defined on $\Omega_h$ is standard and can be bounded by applying verbatim the results of \cite{ConNonConVEM}. We record this bound in the following lemma.

\begin{lemma}\label{lemma::Method2BoundPt1} 
Under the assumptions of Theorem~\ref{theorem::Method2Strang}, if, moreover,  $a,\mathbf{b},\mu \in W_\infty^{s+1}(\Omega)$ and $u \in H^{s+1}(\Omega)$, for some $s\in\{1,\dots, \min\{k,l\}\}$, there holds 
\begin{equation*}
 \sup_{w_h \in V_{h,0}^k\backslash \{ 0 \}}  \frac{\sum_{E_h \in \VirtualMesh} \left| {\tilde{\mathcal{A}}^{E_h}}(p,w_h) - \pertB_h^{E_h}(p,w_h) \right|}{\norm{w_h}_{1,\Omega_h}}
\lesssim h^s \norm{\Exten{u}}_{s+1,\Omega_h},
\end{equation*}
for all $p\in \mathbb{P}_k(\VirtualMesh)$. The hidden constant depends on the shape regularity of $\mesh$, on $k$ and $l$, and on the $W^{s+1}_\infty(\Omega)$-norms of ${a}$, ${\mathbf{b}}$, and ${\mu}$.
\end{lemma}

We proceed with bounding the remaining term three in~\eqref{eq::StrangBound2} relating to the linear functional consistency.

\begin{lemma}\label{lemma::Method2BoundPt3}
Under the assumptions of Theorem~\ref{theorem::Method2Strang}, if, moreover,  $a \in W_\infty^{s}(\Omega)$, $\mathbf{b},\mu \in W_\infty^{s-1}(\Omega)$, $f \in H^{s-1}(\Omega) \cap L^{\infty}(\Omega)$, $u \in H^{s+1}(\Omega) \cap W_\infty^2(\Omega)$ for some $s\in\{1,\dots, \min\{k,l\}\}$, then
\begin{equation*}
    \frac{\left| l_h(w_h) - l_*(w_h) \right|}{\norm{w_h}_{1,\Omega_h}} \lesssim  h^{s}\left( \norm{\Exten{f}}_{s-1,\Omega_h} +  \norm{\Exten{u}}_{s+1,\Omega_h} + \norm{\Exten{f}}_{0,\infty,\Omega_h} + \norm{\Exten{u}}_{2,\infty,\Omega_h} \right),
\end{equation*}
holds true with hidden constant depending on the shape regularity of $\mesh$ and $\RefMesh$, on $k$, $l$,  and $j_0$, and on the $W^{s+1}_\infty(\RefDomain)$-norm of $\F$. 
\begin{proof}
Adding and subtracting $\Pi^0_{k-2} \Exten{f}$ and using the orthogonality of $\Pi^0$ we rewrite the difference term as
\begin{align}
    l_h(w_h) - l_*(w_h) =& \sum_{E_h \in \VirtualMesh}\left[ \int_{\projE} \Exten{f}\ \Pi_{h,k}^0 w_h\ d\PhysCoordinate- \int_{E_h} f_* w_h\ d\PhysCoordinate\right]\nonumber\\
    =& \sum_{E_h \in \VirtualMesh} 
    \Bigg[\left(\int_{\projE} \Pi^0_{k-2} \Exten{f} \ \Pi^0_{h,k} w_h\ d\PhysCoordinate -  \int_{E_h} \Pi^0_{k-2} \Exten{f} \ w_h\ d\PhysCoordinate \right) \nonumber \\
    & \qquad + \int_{\projE} (\Exten{f} - \Pi^0_{k-2} \Exten{f}) \Pi^0_{h,k} w_h\ d\PhysCoordinate
    \nonumber\\
    & \qquad + \int_{E_h} (\Pi^0_{k-2} f_* - f_*) (w_h - \Pi^0_{k-2} w_h)\ d\PhysCoordinate \nonumber \\
    & \qquad + \int_{E_h} (\Exten{f} - f_*) \Pi^0_{k-2} w_h\ d\PhysCoordinate
    \Bigg]
    \nonumber\\
    =& \displaystyle\sum_{E_h \in \VirtualMesh} (T_1^{E_h} + T_2^{E_h}+T_3^{E_h}+T_4^{E_h}).\label{eq:lstar1}
\end{align}
The term $T_1^{E_h}$ can be rewritten using~\eqref{eq::pi_identity} 
as
\begin{equation*}
    T_1^{E_h} = \frac{\abs{E_h^\pi} - \abs{E_h}}{\abs{E_h}} \int_{E_h} \Pi^0_{k-2} \Exten{f}\ w_h\ d\PhysCoordinate.
\end{equation*}
Applying the Cauchy-Schwarz inequality, the stability of $\Pi^0$, and Lemmas \ref{lemma::AreaError} and \ref{lemma::MeshScalingV2} results in
\begin{align*}
    \abs{T_1^{E_h}} \lesssim h_{\hat{E}}^s \norm{\Exten{f}}_{0,E_h} \norm{w_h}_{0,E_h}.
\end{align*}
To bound $T_2^{E_h}$, we introduce the $L^2$-projection onto $\mathbb{P}_{k-2}(\projE)$, denoted by $\mathcal{P}_{k-2}$, and apply \holder's inequality, Lemma \ref{lemma::PolyExten}, and inverse estimates over $E_h$, yielding
\begin{align*}
     \abs{ \int_{\projE} (\mathcal{P}_{k-2} \Exten{f} - \Pi^0_{k-2} \Exten{f}) \Pi^0_{h,k} w_h\ d\PhysCoordinate} &\leq \abs{\projE} \norm{\mathcal{P}_{k-2} \Exten{f} - \Pi^0_{k-2} \Exten{f}}_{0,\infty,\projE} \norm{ \Pi^0_{h,k} w_h}_{0,\infty,\projE} \\
    &\lesssim \abs{\projE} \norm{\mathcal{P}_{k-2} \Exten{f} - \Pi^0_{k-2} \Exten{f}}_{0,\infty,E_h} \norm{ \Pi^0_{h,k} w_h}_{0,\infty,E_h} \\
    &\lesssim \norm{\mathcal{P}_{k-2} \Exten{f} - \Pi^0_{k-2} \Exten{f}}_{0,E_h} \norm{ \Pi^0_{h,k} w_h}_{0,E_h}.
\end{align*}
Following the same arguments as in the proof of Theorem \ref{theorem::PiStarOpBound} results in
\begin{align*}
    \abs{T_2^{E_h}} \lesssim h_{\hat{E}}^s \norm{\Exten{f}}_{0,E_h} \norm{w_h}_{0,E_h}.
\end{align*}

$T_3^{E_h}$ is bounded via Cauchy-Schwarz and Theorem \ref{theorem::Pi0Accuracy}, resulting in
\begin{align*}
    \abs{T_3^{E_h}} \lesssim h_{\hat{E}}^s \norm{f_*}_{s-1,E_h} \norm{w_h}_{1,E_h}.
\end{align*}

We rewrite $T_4^{E_h}$ by splitting the integral into $E_h \backslash E$ and $E_h \cap E$ and noting that $f_* = \Exten{f}$ a.e. on $E_h \cap E$
\begin{align*}
    \int_{E_h} (\Exten{f} - f_*) \Pi^0_{k-2} w_h\ d\PhysCoordinate&= \int_{E_h\setminus E} (\Exten{f} - f_*)\ \Pi^0_{k-2} w_h\ d\PhysCoordinate .
\end{align*}
$T_4^{E_h}$ is bounded by applying \holder's inequality, Lemma \ref{lemma::AreaError}, inverse estimates over $E_h$, and the stability of $\Pi^0$
\begin{align*}
    \abs{\int_{E_h \backslash E}  ( \Exten{f} - f_*)\ \Pi^0_{k-2} w_h\ d\PhysCoordinate} &\leq \abs{E_h \backslash E}  \norm{ \Exten{f} - f_* }_{0,\infty,E_h}  \norm{\Pi^0_{k-2} w_h}_{0,\infty,E_h}\\
    &\lesssim h_{\hat{E}}^{s+1}  \norm{\Exten{f}- f_*}_{0,\infty,E_h} \norm{\Pi^0_{k-2} w_h}_{0,E_h} \\
    &\lesssim h_{\hat{E}}^{s+1}  \norm{\Exten{f}- f_*}_{0,\infty,E_h} \norm{w_h}_{0,E_h}.
\end{align*}

Substituting the bound of each term back into~\eqref{eq:lstar1} and bounding norms gives
\begin{align*}
\abs{ l_h(w_h) - l_*(w_h)}&\lesssim
    \sum_{E_h\in\mesh} h_{\hat{E}}^{s} \left( \norm{\Exten{f}}_{s-1,E_h} + \norm{f_*}_{s-1,E_h} + \norm{\Exten{f}- f_*}_{0,\infty,E_h} \right) \norm{w_h}_{1,E_h}\\
    & \lesssim h^{s}\left(  \norm{ \Exten{f}}_{s-1,\Omega_h } +  \norm{f_*}_{s-1,\Omega_h} + \norm{\Exten{f}- f_*}_{0,\infty,\Omega_h} \right)\norm{w_h}_{1,\Omega_h}.
\end{align*}
Bounding norms of $f_*$ with norms of $\Exten{u}$ using~\eqref{eq::fStarDef}
and dividing by $\norm{w_h}_{1,\Omega_h}$ the required estimate readily follows.
\end{proof}
\end{lemma}

\begin{remark}
    Assuming, as in the classical analysis of isoparametric FEM by Ciarlet~\cite{CiarletElliptic}, that $\Exten{f} \equiv f_*$, we obtain an $O(h^s)$ estimate without requiring the additional $L^\infty$ regularity on $\tilde{u}$ and $\tilde{f}$, cf.~\cite{WellsThesis}. 
    Also, for $\Exten{f} \neq f_*$ as in the above lemma, following  a slightly different line of proof, we can obtain an estimate of the relevant terms of $O(h^{3s/2-1})$, without the additional $L^\infty$ regularity. Thus, the additional regularity assumption can be dropped from the statement of Lemma~\ref{lemma::Method2BoundPt3} when $k,l \geq 2$.
\end{remark}

We are finally in a position to conclude the a priori analysis of the physical IsoVEM with the following error bound.
\begin{theorem}[Physical IsoVEM error estimate]\label{theorem::MethodIIH1Estimate}
Let $\RefMesh$ satisfy Assumption~\ref{assumption::ShapeRegular}, $\F_h$ satisfy Assumption~\ref{assumption::F_hAccuracy}, and $h$ be small enough.  Suppose additionally that the solution to \eqref{eq::EllipticModelPDE} satisfies $u\in H^{s+1}(\Omega) \cap W^2_\infty(\Omega)$ with coefficients $a,\mathbf{b},\mu \in W_\infty^{s+1}(\Omega)$, and forcing term $f \in H^{s-1}(\Omega) \cap L^\infty(\Omega)$, for some $s\in\{1,\dots, \min\{k,l\}\}$. Then, for the solution $u_h\in V_{h,0}^k$ of  the physical IsoVEM~\eqref{eq::physicalIsoVEM}, there holds 
 \begin{align*}
     \norm{\Exten{u} - u_h}_{1,\Omega_h}+ \norm{\nabla\Exten{u} - \Pi^1_{h,k-1}{u}_h}_{0,\Omega_h} \lesssim h^s \Big( &\norm{u}_{s+1,\Omega} + \norm{f}_{s-1,\Omega} \\
      &\hspace{.4cm}+\norm{f}_{0,\infty,\Omega} + \norm{u}_{2,\infty,\Omega} \Big).
 \end{align*}
 The hidden constant is dependent on the shape regularity of $\RefMesh$, on $s$, $k$, $l$, and $j_0$, on the $W^{s+1}_\infty(\RefDomain)$-norm of $\F$, and on the $W^{s+1}_\infty(\Omega)$-norms of $a$, $\mathbf{b}$, and $\mu$.
 
 \begin{proof}
 The proof of the first term's bound is standard, cf. Section~\ref{sec::VEM} and~\cite{ConNonConVEM}. It follows by bounding each term in the Strang-type bound~\eqref{eq::StrangBound2}. The first two terms are bounded by  considering $v_h = \Exten{u}_I$ and $p = \Pi^0_k \Exten{u}$ and then using the polynomial projection and virtual interpolation error bounds from Theorem~\ref{theorem::Pi0Accuracy} and Theorem~\ref{theorem::GlobalCurvedVEMInterp}. The remaining, consistency error, terms are bounded using the results from Lemmata~\ref{lemma::IDKYet}, \ref{lemma::Method2BoundPt1}, and~\ref{lemma::Method2BoundPt3}. The proof is concluded by bounding  norms of the extensions $\Exten{u}$ and $\Exten{f}$ back to $\Omega$ using Stein's extension Theorem~\ref{theorem::SteinExtension}.

To bound the second term we add and subtract the $\Pi^1_{k-1}$ projection of both $\Exten{u}$ and $u_h$ to obtain
\begin{align*}
\norm{\nabla\Exten{u} - \Pi^1_{h,k-1}{u}_h}_{0,\Omega_h}&\le 
\norm{\nabla\Exten{u} - \Pi^1_{k-1}\Exten{u}}_{0,\Omega_h}+
\norm{\Pi^1_{k-1}(\Exten{u}-{u}_h)}_{0,\Omega_h} \\
&+
\norm{\Pi^1_{k-1}{u}_h - \Pi^1_{h,k-1}{u}_h}_{0,\Omega_h}.
\end{align*}
The result easily follows from Theorems~\ref{theorem::Pi0Accuracy} and~\ref{theorem::PiStarValueGrad} and the bound of $ \Exten{u} - u_h$ already derived. 
 \end{proof}
\end{theorem}

\section{Numerical Results}\label{sec::NumericalTests}
\begin{figure}[htbp]

\begin{minipage}[t]{.3\textwidth}
  \centering
  \includegraphics[width=1\linewidth]{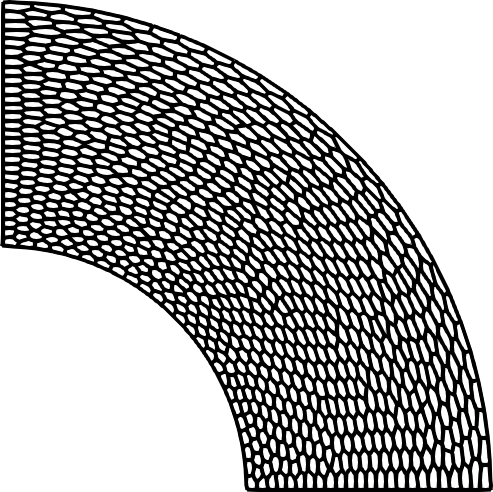}
\end{minipage}
~
\centering\hspace{2cm}
~
\begin{minipage}[t]{.3\textwidth}
  \centering
  \includegraphics[width=1\linewidth]{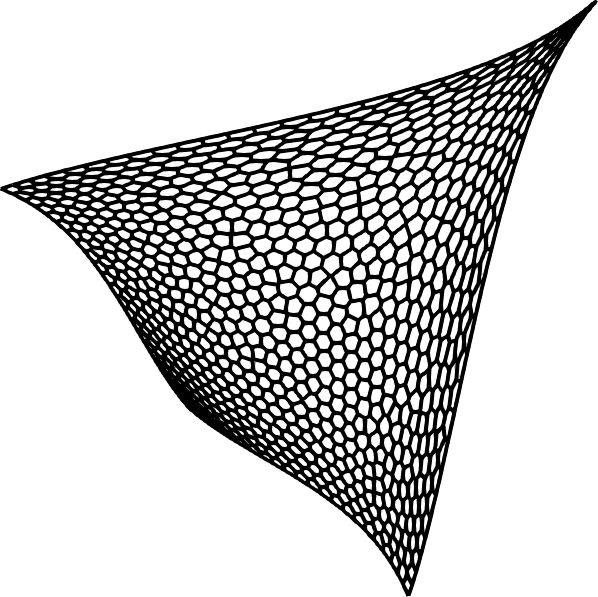}
\end{minipage}
\captionof{figure}{The annulus (left) and plane (right) domains approximated by a VEM mapping on a discretisation of $\RefDomain = [0,1]^2$ using 800 polygonal elements.}\label{fig::TestCases}
\end{figure}

We present numerical experiments in support of the analysis of the reference and physical IsoVEMs. Here, VEM spaces of degree $k,l=1,2,3$ are employed across two distinct domain transformations $\F$ which we refer to as the \textit{annulus} mapping and the \textit{plane} mapping. We refer to \cite{WellsThesis} for a more comprehensive set of numerical experiments for both IsoVEMs. In both experiments, the reference domain is chosen to be the unit square $\RefDomain = [0,1]^2$ with $\Omega$ being given in each test case by $\F(\RefDomain)$. The transformed domains and sample polygonal meshes are shown in Figure~\ref{fig::TestCases}. In the tests we also report results on the $L^2$-norm error, not considered in our a priori analysis.

Noting that the numerical solutions are defined in two different frames, reference and physical, to obtain a quantitatively meaningful comparison, we compute all errors on the physical frame. Thus, for the reference IsoVEM we consider the error norms
\begin{align*}
    \abs{u- u_h}_{h,1}^2 &:= \sum_{\hat{E} \in \RefMesh} \int_{\hat{E}} \abs{ \Jacobian_{\F,h}^{-\top} \ \left( \RefGrad \hat{u}-\Pi^1_{k-1} \hat{u}_h \right) }^2\ j_{h}\ d\RefCoordinate,  \\
    \norm{u - u_h}_{h,0}^2 &:= \sum_{\hat{E} \in \RefMesh} \int_{\hat{E}} \abs{ \hat{u}-\ \Pi^0_{k} \hat{u}_h }^2\ j_{h}\ d\RefCoordinate,
\end{align*}
and compare them with the following directly computable error norms for the physical IsoVEM:
\begin{align*}
    \abs{u - u_h}_{h,1}^2 &:= \sum_{E_h \in \mesh} \int_{\projE} \abs{ \nabla \Exten{u} -\Pi^1_{h,k-1}  u_h }^2\  d\PhysCoordinate, \\
    \norm{u - u_h}_{h,0}^2 &:= \sum_{E_h \in \mesh} \int_{\projE} \abs{\Exten{u}- \Pi^0_{h,k} u_h   }^2\ d\PhysCoordinate.
\end{align*}

The numerical experiments are performed within the DUNE software environment \cite{Bastian2021TheDevelopments, Dedner2012Dune-Fem:Computing, Alkamper2016TheModule} using the Python bindings presented in \cite{Dedner2020PythonModule}. We consider a sequence of 5 CVT Voronoi meshes \cite{Senechal1995SpatialDiagrams} which roughly halve in mesh size with each new generated mesh. The mesh files are generated using PolyMesher within MATLAB \cite{Talischi2012PolyMesher:Matlab} and imported into DUNE.

\paragraph{The Annulus Map.} 
We consider the use of IsoVEM in approximating a annulus shape domain. The polar coordinate system is given by a linear transformation of $(r,\theta) = (\xi_1 + 1, \pi/2\ \xi_2)$ and the domain transformation is then defined as
\[
\F(r, \theta) = \left[ r \cos{\theta},\ r \sin{\theta} \right] = \left[ (\xi_1 + 1) \cos \left( \frac{\pi \xi_2}{2} \right),\ (\xi_1+1) \sin \left( \frac{\pi \xi_2}{2} \right) \right].
\]
Over $\Omega=\F([0,1]^2)$, we consider problem~\eqref{eq::prob1}-\eqref{eq::prob2} with $a=1, \mathbf{b} = (-y, x)$, and $c=0$, whose solution in the physical coordinate system can be shown to be given by $u(r,\theta) = {\ln r}/{\ln 2}$.

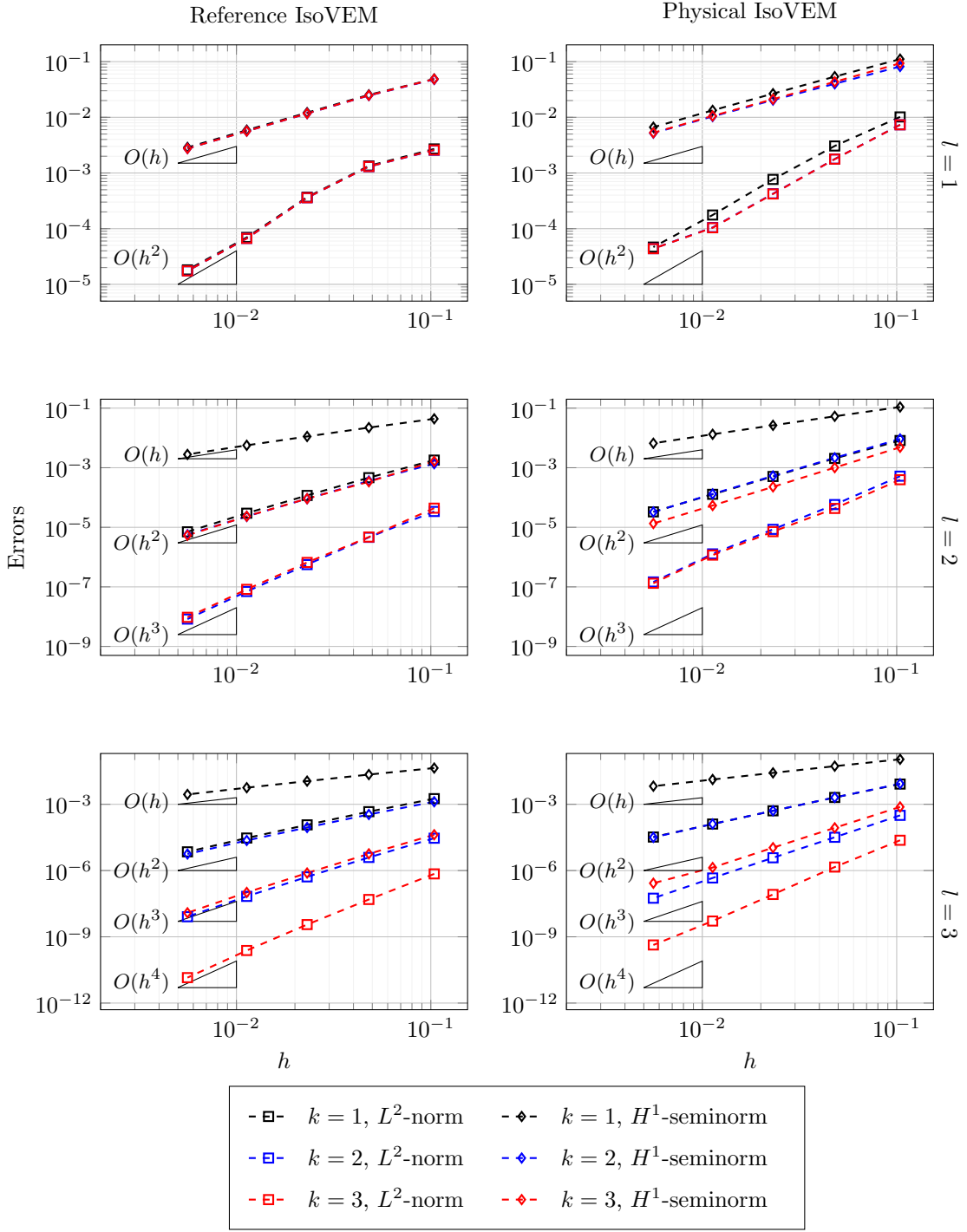
\begin{figure}[htbp]
	\centering
	\begin{tikzpicture}
    \begin{groupplot}[group style={
                      group name= myplot,
                      group size= 2 by 4,
                      horizontal sep=1.5cm,
                      vertical sep=1.5cm},
                      height=5.5cm,
                      width=7.15cm,
                      xmode=log,
                      ymode=log, xmin=0.002,
                      axis background/.style={fill=gray!0}, 
					  legend pos=south east,
					  grid=both,
					  grid style={line width=.1pt, draw=gray!10},
   					  major grid style={line width=.2pt,draw=gray!50}]
                      	
        \nextgroupplot[title = Reference IsoVEM,ymax=0.2, ymin=5e-6]
                \addplot+[mark=square, thick, dashed, black, mark options={black, solid}] table [x=h, y=L2_error, col sep=comma] {data/new_tests/ANNULUS_jacobian_elliptic_mesh_1_sol_1.csv}; 
		    	\addplot+[mark=diamond, thick, dashed, black, mark options={black, solid}] table [x=h, y=H1_error, col sep=comma] {data/new_tests/ANNULUS_jacobian_elliptic_mesh_1_sol_1.csv}; 
		    	\addplot+[mark=square, thick, dashed, blue, mark options={blue, solid}] table [x=h, y=L2_error, col sep=comma] {data/new_tests/ANNULUS_jacobian_elliptic_mesh_1_sol_2.csv}; 
		    	\addplot+[mark=diamond, thick, dashed, blue, mark options={blue, solid}] table [x=h, y=H1_error, col sep=comma] {data/new_tests/ANNULUS_jacobian_elliptic_mesh_1_sol_2.csv}; 
		    	\addplot+[mark=square, thick, dashed, red, mark options={red, solid}] table [x=h, y=L2_error, col sep=comma] {data/new_tests/ANNULUS_jacobian_elliptic_mesh_1_sol_3.csv}; 
	    		\addplot+[mark=diamond, thick, dashed, red, mark options={red, solid}] table [x=h, y=H1_error, col sep=comma] {data/new_tests/ANNULUS_jacobian_elliptic_mesh_1_sol_3.csv}; 
	    		
	    		\addplot[mark=none, solid, black] coordinates {(0.005,1e-5) (0.01,4e-5) (0.01,1e-5) (0.005,1e-5)};
	    		\plot[mark=none] (0.005,3e-5) node[anchor=east] {$O(h^2)$};
	    		
	    		\addplot[mark=none, solid, black] coordinates {(0.005,3*5e-4) (0.01,3*5*2e-4) (0.01,5*3*1e-4) (0.005,5*3e-4)};
	    		\plot[mark=none] (0.005,3*6e-4) node[anchor=east] {$O(h)$};
	    		
        \nextgroupplot[title = Physical IsoVEM,ymax=0.2, ymin=5e-6]
                \addplot+[mark=square, thick, dashed, black, mark options={black, solid}] table [x=h, y=L2_error, col sep=comma] {data/new_tests/ANNULUS_mapped_elliptic_mesh_1_sol_1.csv}; 
		    	\addplot+[mark=diamond, thick, dashed, black, mark options={black, solid}] table [x=h, y=H1_error, col sep=comma] {data/new_tests/ANNULUS_mapped_elliptic_mesh_1_sol_1.csv}; 
		    	\addplot+[mark=square, thick, dashed, blue, mark options={blue, solid}] table [x=h, y=L2_error, col sep=comma] {data/new_tests/ANNULUS_mapped_elliptic_mesh_1_sol_2.csv}; 
		    	\addplot+[mark=diamond, thick, dashed, blue, mark options={blue, solid}] table [x=h, y=H1_error, col sep=comma] {data/new_tests/ANNULUS_mapped_elliptic_mesh_1_sol_2.csv}; 
		    	\addplot+[mark=square, thick, dashed, red, mark options={red, solid}] table [x=h, y=L2_error, col sep=comma] {data/new_tests/ANNULUS_mapped_elliptic_mesh_1_sol_3.csv}; 
	    		\addplot+[mark=diamond, thick, dashed, red, mark options={red, solid}] table [x=h, y=H1_error, col sep=comma] {data/new_tests/ANNULUS_mapped_elliptic_mesh_1_sol_3.csv}; 
	    		
	    		\addplot[mark=none, solid, black] coordinates {(0.005,1e-5) (0.01,4e-5) (0.01,1e-5) (0.005,1e-5)};
	    		\plot[mark=none] (0.005,3e-5) node[anchor=east] {$O(h^2)$};
	    		
	    		\addplot[mark=none, solid, black] coordinates {(0.005,3*5e-4) (0.01,3*5*2e-4) (0.01,5*3*1e-4) (0.005,5*3e-4)};
	    		\plot[mark=none] (0.005,3*6e-4) node[anchor=east] {$O(h)$};

       \nextgroupplot[ymax=0.2, ymin=5e-10]
                \addplot+[mark=square, thick, dashed, black, mark options={black, solid}] table [x=h, y=L2_error, col sep=comma] {data/new_tests/ANNULUS_jacobian_elliptic_mesh_2_sol_1.csv}; 
		    	\addplot+[mark=diamond, thick, dashed, black, mark options={black, solid}] table [x=h, y=H1_error, col sep=comma] {data/new_tests/ANNULUS_jacobian_elliptic_mesh_2_sol_1.csv}; 
		    	\addplot+[mark=square, thick, dashed, blue, mark options={blue, solid}] table [x=h, y=L2_error, col sep=comma] {data/new_tests/ANNULUS_jacobian_elliptic_mesh_2_sol_2.csv}; 
		    	\addplot+[mark=diamond, thick, dashed, blue, mark options={blue, solid}] table [x=h, y=H1_error, col sep=comma] {data/new_tests/ANNULUS_jacobian_elliptic_mesh_2_sol_2.csv}; 
		    	\addplot+[mark=square, thick, dashed, red, mark options={red, solid}] table [x=h, y=L2_error, col sep=comma] {data/new_tests/ANNULUS_jacobian_elliptic_mesh_2_sol_3.csv}; 
	    		\addplot+[mark=diamond, thick, dashed, red, mark options={red, solid}] table [x=h, y=H1_error, col sep=comma] {data/new_tests/ANNULUS_jacobian_elliptic_mesh_2_sol_3.csv}; 
	    		
		        \addplot[mark=none, solid, black] coordinates {(0.005,0.25e-8) (0.01,2e-8) (0.01,0.25e-8) (0.005,0.25e-8)};
	    		\plot[mark=none] (0.005,0.25e-8) node[anchor=east] {$O(h^3)$};
	    		
	    		\addplot[mark=none, solid, black] coordinates {(0.005,0.3*1e-5) (0.01,0.3*4e-5) (0.01,0.3*1e-5) (0.005,0.3*1e-5)};
	    		\plot[mark=none] (0.005,0.3*1e-5) node[anchor=east] {$O(h^2)$};
	    		
	    		\addplot[mark=none, solid, black] coordinates {(0.005,0.4*5e-3) (0.01,0.4*5*2e-3) (0.01,0.4*5*1e-3) (0.005,0.4*5e-3)};
	    		\plot[mark=none] (0.005,0.4*6e-3) node[anchor=east] {$O(h)$};
	    		
        \nextgroupplot[ymax=0.2, ymin=5e-10]
                \addplot+[mark=square, thick, dashed, black, mark options={black, solid}] table [x=h, y=L2_error, col sep=comma] {data/new_tests/ANNULUS_mapped_elliptic_mesh_2_sol_1.csv}; 
		    	\addplot+[mark=diamond, thick, dashed, black, mark options={black, solid}] table [x=h, y=H1_error, col sep=comma] {data/new_tests/ANNULUS_mapped_elliptic_mesh_2_sol_1.csv}; 
		    	\addplot+[mark=square, thick, dashed, blue, mark options={blue, solid}] table [x=h, y=L2_error, col sep=comma] {data/new_tests/ANNULUS_mapped_elliptic_mesh_2_sol_2.csv}; 
		    	\addplot+[mark=diamond, thick, dashed, blue, mark options={blue, solid}] table [x=h, y=H1_error, col sep=comma] {data/new_tests/ANNULUS_mapped_elliptic_mesh_2_sol_2.csv}; 
		    	\addplot+[mark=square, thick, dashed, red, mark options={red, solid}] table [x=h, y=L2_error, col sep=comma] {data/new_tests/ANNULUS_mapped_elliptic_mesh_2_sol_3.csv}; 
	    		\addplot+[mark=diamond, thick, dashed, red, mark options={red, solid}] table [x=h, y=H1_error, col sep=comma] {data/new_tests/ANNULUS_mapped_elliptic_mesh_2_sol_3.csv}; 
	    		
		        \addplot[mark=none, solid, black] coordinates {(0.005,0.25e-8) (0.01,2e-8) (0.01,0.25e-8) (0.005,0.25e-8)};
	    		\plot[mark=none] (0.005,0.25e-8) node[anchor=east] {$O(h^3)$};
	    		
	    		\addplot[mark=none, solid, black] coordinates {(0.005,0.3*1e-5) (0.01,0.3*4e-5) (0.01,0.3*1e-5) (0.005,0.3*1e-5)};
	    		\plot[mark=none] (0.005,0.3*1e-5) node[anchor=east] {$O(h^2)$};
	    		
	    		\addplot[mark=none, solid, black] coordinates {(0.005,0.4*5e-3) (0.01,0.4*5*2e-3) (0.01,0.4*5*1e-3) (0.005,0.4*5e-3)};
	    		\plot[mark=none] (0.005,0.4*6e-3) node[anchor=east] {$O(h)$};

       \nextgroupplot[xlabel={$h$},ymax=0.2, ymin=5e-13]
                \addplot+[mark=square, thick, dashed, black, mark options={black, solid}] table [x=h, y=L2_error, col sep=comma] {data/new_tests/ANNULUS_jacobian_elliptic_mesh_3_sol_1.csv}; 
		    	\addplot+[mark=diamond, thick, dashed, black, mark options={black, solid}] table [x=h, y=H1_error, col sep=comma] {data/new_tests/ANNULUS_jacobian_elliptic_mesh_3_sol_1.csv}; 
		    	\addplot+[mark=square, thick, dashed, blue, mark options={blue, solid}] table [x=h, y=L2_error, col sep=comma] {data/new_tests/ANNULUS_jacobian_elliptic_mesh_3_sol_2.csv}; 
		    	\addplot+[mark=diamond, thick, dashed, blue, mark options={blue, solid}] table [x=h, y=H1_error, col sep=comma] {data/new_tests/ANNULUS_jacobian_elliptic_mesh_3_sol_2.csv}; 
		    	\addplot+[mark=square, thick, dashed, red, mark options={red, solid}] table [x=h, y=L2_error, col sep=comma] {data/new_tests/ANNULUS_jacobian_elliptic_mesh_3_sol_3.csv}; 
	    		\addplot+[mark=diamond, thick, dashed, red, mark options={red, solid}] table [x=h, y=H1_error, col sep=comma] {data/new_tests/ANNULUS_jacobian_elliptic_mesh_3_sol_3.csv}; 
	    		
	    		\addplot[mark=none, solid, black] coordinates {(0.005,0.5e-11) (0.01,0.8e-10) (0.01,0.5e-11) (0.005,0.5e-11)};
	    		\plot[mark=none] (0.005,1e-11) node[anchor=east] {$O(h^4)$};
		
		        \addplot[mark=none, solid, black] coordinates {(0.005,0.5*1e-8) (0.01,0.5*4*2e-8) (0.01,0.5*1e-8) (0.005,0.5*1e-8)};
	    		\plot[mark=none] (0.005,0.5*1e-8) node[anchor=east] {$O(h^3)$};
	    		
	    		\addplot[mark=none, solid, black] coordinates {(0.005,0.1*1e-5) (0.01,0.1*4e-5) (0.01,0.1*1e-5) (0.005,0.1*1e-5)};
	    		\plot[mark=none] (0.005,0.1*1e-5) node[anchor=east] {$O(h^2)$};
	    		
	    		\addplot[mark=none, solid, black] coordinates {(0.005,0.2*5e-3) (0.01,0.2*5*2e-3) (0.01,0.2*5*1e-3) (0.005,0.2*5e-3)};
	    		\plot[mark=none] (0.005,0.2*6e-3) node[anchor=east] {$O(h)$};
	    		
        \nextgroupplot[xlabel={$h$},ymax=0.2, ymin=5e-13]
                \addplot+[mark=square, thick, dashed, black, mark options={black, solid}] table [x=h, y=L2_error, col sep=comma] {data/new_tests/ANNULUS_mapped_elliptic_mesh_3_sol_1.csv}; 
		    	\addplot+[mark=diamond, thick, dashed, black, mark options={black, solid}] table [x=h, y=H1_error, col sep=comma] {data/new_tests/ANNULUS_mapped_elliptic_mesh_3_sol_1.csv}; 
		    	\addplot+[mark=square, thick, dashed, blue, mark options={blue, solid}] table [x=h, y=L2_error, col sep=comma] {data/new_tests/ANNULUS_mapped_elliptic_mesh_3_sol_2.csv}; 
		    	\addplot+[mark=diamond, thick, dashed, blue, mark options={blue, solid}] table [x=h, y=H1_error, col sep=comma] {data/new_tests/ANNULUS_mapped_elliptic_mesh_3_sol_2.csv}; 
		    	\addplot+[mark=square, thick, dashed, red, mark options={red, solid}] table [x=h, y=L2_error, col sep=comma] {data/new_tests/ANNULUS_mapped_elliptic_mesh_3_sol_3.csv}; 
	    		\addplot+[mark=diamond, thick, dashed, red, mark options={red, solid}] table [x=h, y=H1_error, col sep=comma] {data/new_tests/ANNULUS_mapped_elliptic_mesh_3_sol_3.csv}; 
	    		
	    		\addplot[mark=none, solid, black] coordinates {(0.005,0.5e-11) (0.01,0.8e-10) (0.01,0.5e-11) (0.005,0.5e-11)};
	    		\plot[mark=none] (0.005,1e-11) node[anchor=east] {$O(h^4)$};
		
		        \addplot[mark=none, solid, black] coordinates {(0.005,0.5*1e-8) (0.01,0.5*4*2e-8) (0.01,0.5*1e-8) (0.005,0.5*1e-8)};
	    		\plot[mark=none] (0.005,0.5*1e-8) node[anchor=east] {$O(h^3)$};
	    		
	    		\addplot[mark=none, solid, black] coordinates {(0.005,0.1*1e-5) (0.01,0.1*4e-5) (0.01,0.1*1e-5) (0.005,0.1*1e-5)};
	    		\plot[mark=none] (0.005,0.1*1e-5) node[anchor=east] {$O(h^2)$};
	    		
	    		\addplot[mark=none, solid, black] coordinates {(0.005,0.2*5e-3) (0.01,0.2*5*2e-3) (0.01,0.2*5*1e-3) (0.005,0.2*5e-3)};
	    		\plot[mark=none] (0.005,0.2*6e-3) node[anchor=east] {$O(h)$};
       
    \end{groupplot}
    \path (myplot c1r1.outer north west)
          -- node[anchor=south,rotate=90] {Errors}
          (myplot c1r3.outer south west)
    ;
    \path (myplot c2r1.outer north east)
          -- node[anchor=south,rotate=-90] {$l=1$}
          (myplot c2r1.outer south east)
    ;
    \path (myplot c2r2.outer north east)
          -- node[anchor=south,rotate=-90] {$l=2$}
          (myplot c2r2.outer south east)
    ;
    \path (myplot c2r3.outer north east)
          -- node[anchor=south,rotate=-90] {$l=3$}
          (myplot c2r3.outer south east)
    ;
\path (myplot c1r3.south west|-current bounding box.south)--
      coordinate(legendpos)
      (myplot c2r3.south east|-current bounding box.south);
\matrix[
    matrix of nodes,
    anchor=north,
    draw,
    inner sep=0.4em,
    draw
  ]at([yshift=-1ex]legendpos)
  {
    \ref{plots:K1L2}& $k=1$, $L^2$-norm &[8pt]
    \ref{plots:K1H1}& $k=1$, $H^1$-seminorm &[8pt] \\
    \ref{plots:K2L2}& $k=2$, $L^2$-norm &[8pt]
    \ref{plots:K2H1}& $k=2$, $H^1$-seminorm &[8pt] \\
    \ref{plots:K3L2}& $k=3$, $L^2$-norm &[8pt]
    \ref{plots:K3H1}& $k=3$, $H^1$-seminorm &[8pt]
    \\};
\end{tikzpicture}
\caption{Error plots for the annulus mapping problem with isoparametric elements of degree $l,k=1,2,3$.}\label{fig::Annulus}
\end{figure}

To demonstrate the dependency on the convergence rates on both $k,l$, as outlined in the $H^1$-error estimates of Theorems~\ref{theorem::MethodIH1Error} and~\ref{theorem::MethodIIH1Estimate}, we consider all possible combinations of VEM discretisation degrees of $l,k=1,2,3$. 
In Figure \ref{fig::Annulus} the numerical errors are presented for each possible combination of $l$ and $k$.  For both IsoVEMs, we observe that the $H^1$-seminorm error converges with order $\min \{ l, k \}$ in agreement with theory. In addition, we observe that the $L^2$-norm errors converge in each case with order $1+\min \{ l, k \}$.
 These results demonstrate that both IsoVEMs are optimal in the $H^1$- and $L^2$-norms with respect to $\min \{ l, k \}$.  We observe that in this experiment, where the domain approximation is expected to dominate the solution error, as long as $k\ge l$, the choice of solution degree $k$ has little impact on the $H^1$-norm and $L^2$-norm errors. Additionally, the reference IsoVEM produces slightly better $H^1$ and $L^2$ errors in the higher-order cases.

\paragraph{The Plane Map.} Our second test applies the IsoVEM to a more complex PDE and domain transformation. For ease of implementation we construct our test case by defining the solution in the physical coordinates first and computing the forcing data and boundary conditions accordingly.   
The domain translation is defined by,
\begin{equation*}
    \mathbf{F}(\RefCoordinate) = \left[ \tanh (2\xi_1-1)\exp (\xi_2),\ \tanh (2 \xi_2-1) \exp (\xi_1) \right].
\end{equation*}
Over $\Omega=\F([0,1]^2)$, we consider the general elliptic equation~\eqref{eq::prob1} with coefficients
\begin{align*}
    a = 
    \begin{pmatrix}
y^2+1 & -xy \\
-xy & x^2+1
\end{pmatrix}, \qquad \mathbf{b} = (x,y), \qquad c = x^2 + y^3 + 2,
\end{align*}
The forcing term is chosen such that the true solution in the physical coordinates is given by $u(\mathbf{x}) = \cos (xy) \sin (x)$.

\begin{figure}[htbp]
	\centering
	\begin{tikzpicture}
    \begin{groupplot}[group style={
                      group name= myplot,
                      group size= 2 by 2,
                      horizontal sep=1.5cm,
                      vertical sep=1.5cm},
                      height=5.5cm,
                      width=7.15cm,
                      xmode=log,
                      ymode=log, xmin=0.002,
                      axis background/.style={fill=gray!0}, 
					  legend pos=south east,
					  grid=both,
					  grid style={line width=.1pt, draw=gray!10},
   					  major grid style={line width=.2pt,draw=gray!50}]
                      	
        \nextgroupplot[xlabel={$h$},title = Reference IsoVEM,ymax=3, ymin=5e-9]
                \addplot+[mark=square, thick, dashed, black, mark options={black, solid}] table [x=h, y=L2_error, col sep=comma] {data/new_tests/PLANE_jacobian_elliptic_mesh_1_sol_1.csv};
                \label{plots:K1L2}
		    	\addplot+[mark=diamond, thick, dashed, black, mark options={black, solid}] table [x=h, y=H1_error, col sep=comma] {data/new_tests/PLANE_jacobian_elliptic_mesh_1_sol_1.csv};
                \label{plots:K1H1}
		    	\addplot+[mark=square, thick, dashed, blue, mark options={blue, solid}] table [x=h, y=L2_error, col sep=comma] {data/new_tests/PLANE_jacobian_elliptic_mesh_2_sol_2.csv}; 
                \label{plots:K2L2}
                \addplot+[mark=diamond, thick, dashed, blue, mark options={blue, solid}] table [x=h, y=H1_error, col sep=comma] {data/new_tests/PLANE_jacobian_elliptic_mesh_2_sol_2.csv}; 
                \label{plots:K2H1}
                \addplot+[mark=square, thick, dashed, red, mark options={red, solid}] table [x=h, y=L2_error, col sep=comma] {data/new_tests/PLANE_jacobian_elliptic_mesh_3_sol_3.csv}; 
                \label{plots:K3L2}
                \addplot+[mark=diamond, thick, dashed, red, mark options={red, solid}] table [x=h, y=H1_error, col sep=comma] {data/new_tests/PLANE_jacobian_elliptic_mesh_3_sol_3.csv}; 
                \label{plots:K3H1}

		        \addplot[mark=none, solid, black] coordinates {(0.005,250*2*5*0.5e-11) (0.01,250*2*5*0.8e-10) (0.01,250*2*5*0.5e-11) (0.005,250*2*5*0.5e-11)};
	    		\plot[mark=none] (0.005,250*2*5*0.5e-11) node[anchor=east] {$O(h^4)$};
		
		        \addplot[mark=none, solid, black] coordinates {(0.005,5*6*2*0.25e-7) (0.01,5*6*2*2e-7) (0.01,5*6*2*0.25e-7) (0.005,5*6*2*0.25e-7)};
	    		\plot[mark=none] (0.005,5*6*2*0.25e-7) node[anchor=east] {$O(h^3)$};
	    		
	    		\addplot[mark=none, solid, black] coordinates {(0.005,20*1e-5) (0.01,20*4e-5) (0.01,20*1e-5) (0.005,20*1e-5)};
	    		\plot[mark=none] (0.005,20*1e-5) node[anchor=east] {$O(h^2)$};
	    		
	    		\addplot[mark=none, solid, black] coordinates {(0.005,2*2*5e-3) (0.01,2*2*5*2e-3) (0.01,2*2*5*1e-3) (0.005,2*2*5e-3)};
	    		\plot[mark=none] (0.005,2*2*5e-3) node[anchor=east] {$O(h)$};
	    		
        \nextgroupplot[xlabel={$h$}, title = Physical IsoVEM,ymax=3, ymin=5e-9]
                \addplot+[mark=square, thick, dashed, black, mark options={black, solid}] table [x=h, y=L2_error, col sep=comma] {data/new_tests/PLANE_mapped_elliptic_mesh_1_sol_1.csv}; 
		    	\addplot+[mark=diamond, thick, dashed, black, mark options={black, solid}] table [x=h, y=H1_error, col sep=comma] {data/new_tests/PLANE_mapped_elliptic_mesh_1_sol_1.csv}; 
		    	\addplot+[mark=square, thick, dashed, blue, mark options={blue, solid}] table [x=h, y=L2_error, col sep=comma] {data/new_tests/PLANE_mapped_elliptic_mesh_2_sol_2.csv}; 
		    	\addplot+[mark=diamond, thick, dashed, blue, mark options={blue, solid}] table [x=h, y=H1_error, col sep=comma] {data/new_tests/PLANE_mapped_elliptic_mesh_2_sol_2.csv}; 
		    	\addplot+[mark=square, thick, dashed, red, mark options={red, solid}] table [x=h, y=L2_error, col sep=comma] {data/new_tests/PLANE_mapped_elliptic_mesh_3_sol_3.csv}; 
	    		\addplot+[mark=diamond, thick, dashed, red, mark options={red, solid}] table [x=h, y=H1_error, col sep=comma] {data/new_tests/PLANE_mapped_elliptic_mesh_3_sol_3.csv}; 
	    		
	    		\addplot[mark=none, solid, black] coordinates {(0.005,250*2*5*0.5e-11) (0.01,250*2*5*0.8e-10) (0.01,250*2*5*0.5e-11) (0.005,250*2*5*0.5e-11)};
	    		\plot[mark=none] (0.005,250*2*5*0.5e-11) node[anchor=east] {$O(h^4)$};
		
		        \addplot[mark=none, solid, black] coordinates {(0.005,5*6*2*0.25e-7) (0.01,5*6*2*2e-7) (0.01,5*6*2*0.25e-7) (0.005,5*6*2*0.25e-7)};
	    		\plot[mark=none] (0.005,5*6*2*0.25e-7) node[anchor=east] {$O(h^3)$};
	    		
	    		\addplot[mark=none, solid, black] coordinates {(0.005,20*1e-5) (0.01,20*4e-5) (0.01,20*1e-5) (0.005,20*1e-5)};
	    		\plot[mark=none] (0.005,20*1e-5) node[anchor=east] {$O(h^2)$};
	    		
	    		\addplot[mark=none, solid, black] coordinates {(0.005,2*2*5e-3) (0.01,2*2*5*2e-3) (0.01,2*2*5*1e-3) (0.005,2*2*5e-3)};
	    		\plot[mark=none] (0.005,2*2*5e-3) node[anchor=east] {$O(h)$};

    \end{groupplot}
\path (myplot c1r1.south west|-current bounding box.south)--
      coordinate(legendpos)
      (myplot c2r1.south east|-current bounding box.south);
\matrix[
    matrix of nodes,
    anchor=north,
    draw,
    inner sep=0.4em,
    draw
  ]at([yshift=-1ex]legendpos)
  {
    \ref{plots:K1L2}& $k=1$, $L^2$-norm &[8pt]
    \ref{plots:K1H1}& $k=1$, $H^1$-seminorm &[8pt] \\
    \ref{plots:K2L2}& $k=2$, $L^2$-norm &[8pt]
    \ref{plots:K2H1}& $k=2$, $H^1$-seminorm &[8pt] \\
    \ref{plots:K3L2}& $k=3$, $L^2$-norm &[8pt]
    \ref{plots:K3H1}& $k=3$, $H^1$-seminorm &[8pt]
    \\};
\end{tikzpicture}
\caption{Error plots for the plane mapping problem with isoparametric elements ($l=k$) of degree $k=1,2,3$.} \label{fig::DiamondResultl=1}
\end{figure}
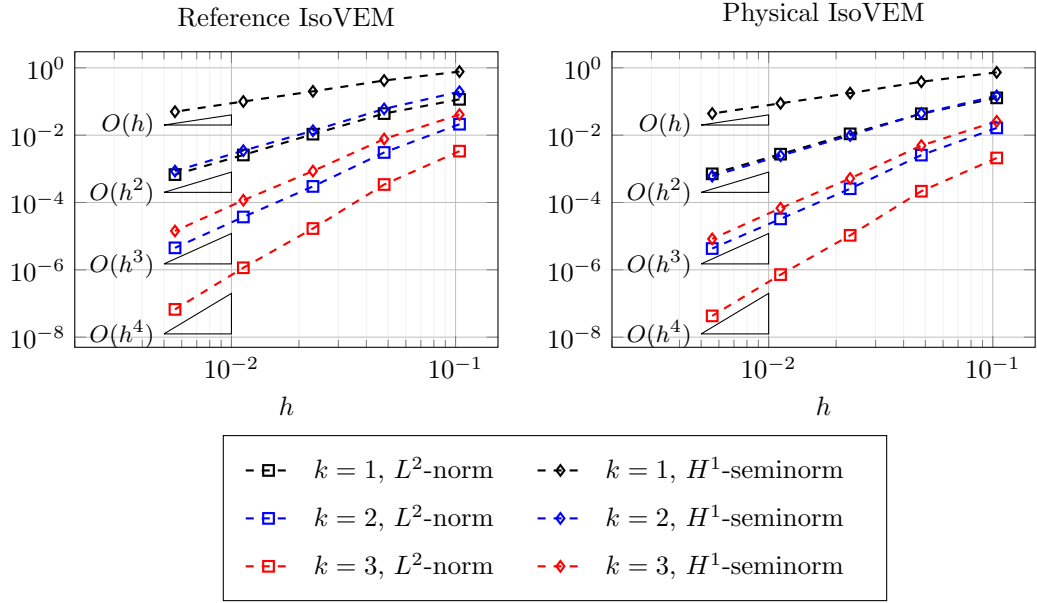

In Figure \ref{fig::DiamondResultl=1}, the numerical errors are presented only for the truly isoparametric case, namely for $l=k$, with $k=1,2,3$. For both IsoVEMs, we observe that the $H^1$-seminorm error converges with order $k$ in agreement with the estimates of Theorems~\ref{theorem::MethodIH1Error} and~\ref{theorem::MethodIIH1Estimate}. In addition to the $H^1$-norm errors we observe that the $L^2$-norm errors converge in each case with order $k+1$. 
 Unlike the annulus test, the numerical errors produced by the two methods are very similar.
These results demonstrate that both IsoVEMs are optimal in the $H^1$- and $L^2$-norm with respect to the isoparametric discretisation degree. 

\section{Conclusion}\label{sec::Conclusion}
We presented two  isoparametric VEM: a method that approximates the  problem on a reference computational domain and a method that constructs a discretization on a physical \emph{virtual} domain. In both cases, optimal $H^1$-norm error estimates are proven and numerical results are shown supporting the theoretical results. In addition, we have verified experimentally that  the methods achieve optimal  orders of convergence in the $L^2$-norm as well.
Extensions of the analysis to $L^2$-norm estimates and three-dimensional problems will be the subject of future work.

In contrast to other works in curved virtual element methods, this paper considers the mapping between domains to be subject to an approximation, not known explicitly \textit{a priori} and only computable through a set of degrees of freedom. This perspective naturally allows for the consideration of more complex problems upon which the mapping is the subject of further numerical approximation. As such, our analysis provides the basis for the development of accurate VEMs for time-dependent problems requiring discretization of both the solution and domain evolution, such as fluid-structure interaction problems.

\

{\bf Acknowledgements.}
This work benefited from the research visit of AD at SISSA funded by the visiting professors fund of the Gruppo Nazionale Calcolo Scientifico-Istituto Nazionale di Alta Matematica (GNCS-INdAM), of which AC is a member. AC was partially supported by the iNEST - Interconnected Nord-Est Innovation Ecosystem funded by the PNRR-NextGenerationEU programme. HW was supported by EPSRC doctoral training grants EP/N50970X/1 and EP/R513283/1. These support is gratefully acknowledged.

\bibliography{references}


\begin{thebibliography}{43}
\ifx \bisbn   \undefined \def \bisbn  #1{ISBN #1}\fi
\ifx \binits  \undefined \def \binits#1{#1}\fi
\ifx \bauthor  \undefined \def \bauthor#1{#1}\fi
\ifx \batitle  \undefined \def \batitle#1{#1}\fi
\ifx \bjtitle  \undefined \def \bjtitle#1{#1}\fi
\ifx \bvolume  \undefined \def \bvolume#1{\textbf{#1}}\fi
\ifx \byear  \undefined \def \byear#1{#1}\fi
\ifx \bissue  \undefined \def \bissue#1{#1}\fi
\ifx \bfpage  \undefined \def \bfpage#1{#1}\fi
\ifx \blpage  \undefined \def \blpage #1{#1}\fi
\ifx \burl  \undefined \def \burl#1{\textsf{#1}}\fi
\ifx \doiurl  \undefined \def \doiurl#1{\url{https://doi.org/#1}}\fi
\ifx \betal  \undefined \def \betal{\textit{et al.}}\fi
\ifx \binstitute  \undefined \def \binstitute#1{#1}\fi
\ifx \binstitutionaled  \undefined \def \binstitutionaled#1{#1}\fi
\ifx \bctitle  \undefined \def \bctitle#1{#1}\fi
\ifx \beditor  \undefined \def \beditor#1{#1}\fi
\ifx \bpublisher  \undefined \def \bpublisher#1{#1}\fi
\ifx \bbtitle  \undefined \def \bbtitle#1{#1}\fi
\ifx \bedition  \undefined \def \bedition#1{#1}\fi
\ifx \bseriesno  \undefined \def \bseriesno#1{#1}\fi
\ifx \blocation  \undefined \def \blocation#1{#1}\fi
\ifx \bsertitle  \undefined \def \bsertitle#1{#1}\fi
\ifx \bsnm \undefined \def \bsnm#1{#1}\fi
\ifx \bsuffix \undefined \def \bsuffix#1{#1}\fi
\ifx \bparticle \undefined \def \bparticle#1{#1}\fi
\ifx \barticle \undefined \def \barticle#1{#1}\fi
\bibcommenthead
\ifx \bconfdate \undefined \def \bconfdate #1{#1}\fi
\ifx \botherref \undefined \def \botherref #1{#1}\fi
\ifx \url \undefined \def \url#1{\textsf{#1}}\fi
\ifx \bchapter \undefined \def \bchapter#1{#1}\fi
\ifx \bbook \undefined \def \bbook#1{#1}\fi
\ifx \bcomment \undefined \def \bcomment#1{#1}\fi
\ifx \oauthor \undefined \def \oauthor#1{#1}\fi
\ifx \citeauthoryear \undefined \def \citeauthoryear#1{#1}\fi
\ifx \endbibitem  \undefined \def \endbibitem {}\fi
\ifx \bconflocation  \undefined \def \bconflocation#1{#1}\fi
\ifx \arxivurl  \undefined \def \arxivurl#1{\textsf{#1}}\fi
\csname PreBibitemsHook\endcsname

\bibitem[\protect\citeauthoryear{Beir{\~{a}}o~da Veiga
  et~al.}{2013}]{basicprinciples}
\begin{barticle}
\bauthor{\bsnm{Veiga}, \binits{L.}},
\bauthor{\bsnm{Brezzi}, \binits{F.}},
\bauthor{\bsnm{Cangiani}, \binits{A.}},
\bauthor{\bsnm{Manzini}, \binits{G.}},
\bauthor{\bsnm{Marini}, \binits{L.D.}},
\bauthor{\bsnm{Russo}, \binits{A.}}:
\batitle{{Basic principles of the virtual element method}}.
\bjtitle{Mathematical Models and Methods in Applied Sciences}
\bvolume{23}(\bissue{01}),
\bfpage{199}--\blpage{214}
(\byear{2013})
\doiurl{10.1142/S0218202512500492}
\end{barticle}
\endbibitem

\bibitem[\protect\citeauthoryear{Lipnikov and
  Morgan}{2019}]{Lipnikov2019AMeshes}
\begin{barticle}
\bauthor{\bsnm{Lipnikov}, \binits{K.}},
\bauthor{\bsnm{Morgan}, \binits{N.}}:
\batitle{{A high-order discontinuous Galerkin method for level set problems on
  polygonal meshes}}.
\bjtitle{Journal of Computational Physics}
\bvolume{397},
\bfpage{108834}
(\byear{2019})
\doiurl{10.1016/j.jcp.2019.07.033}
\end{barticle}
\endbibitem

\bibitem[\protect\citeauthoryear{Lipnikov and
  Morgan}{2020}]{Lipnikov2020ConservativeMeshes}
\begin{barticle}
\bauthor{\bsnm{Lipnikov}, \binits{K.}},
\bauthor{\bsnm{Morgan}, \binits{N.}}:
\batitle{{Conservative high-order discontinuous Galerkin remap scheme on
  curvilinear polyhedral meshes}}.
\bjtitle{Journal of Computational Physics}
\bvolume{420},
\bfpage{109712}
(\byear{2020})
\doiurl{10.1016/j.jcp.2020.109712}
\end{barticle}
\endbibitem

\bibitem[\protect\citeauthoryear{Mazzia
  et~al.}{2020}]{Mazzia2020VirtualEnvironments}
\begin{botherref}
\oauthor{\bsnm{Mazzia}, \binits{A.}},
\oauthor{\bsnm{Ferronato}, \binits{M.}},
\oauthor{\bsnm{Teatini}, \binits{P.}},
\oauthor{\bsnm{Zoccarato}, \binits{C.}}:
{Virtual element method for the numerical simulation of long-term dynamics of
  transitional environments}.
Journal of Computational Physics
\textbf{407}
(2020)
\doiurl{10.1016/j.jcp.2020.109235}
\end{botherref}
\endbibitem

\bibitem[\protect\citeauthoryear{Wells et~al.}{2024}]{Wells2024AMethod}
\begin{barticle}
\bauthor{\bsnm{Wells}, \binits{H.}},
\bauthor{\bsnm{Hubbard}, \binits{M.E.}},
\bauthor{\bsnm{Cangiani}, \binits{A.}}:
\batitle{{A velocity-based moving mesh virtual element method}}.
\bjtitle{Computers {\&} Mathematics with Applications}
\bvolume{155},
\bfpage{110}--\blpage{125}
(\byear{2024})
\doiurl{10.1016/j.camwa.2023.12.005}
\end{barticle}
\endbibitem

\bibitem[\protect\citeauthoryear{Antonietti
  et~al.}{2019}]{Antonietti2019NumericalGrids}
\begin{barticle}
\bauthor{\bsnm{Antonietti}, \binits{P.}},
\bauthor{\bsnm{Verani}, \binits{M.}},
\bauthor{\bsnm{Vergara}, \binits{C.}},
\bauthor{\bsnm{Zonca}, \binits{S.}}:
\batitle{{Numerical solution of fluid-structure interaction problems by means
  of a high order Discontinuous Galerkin method on polygonal grids}}.
\bjtitle{Finite Elements in Analysis and Design}
\bvolume{159}(\bissue{February}),
\bfpage{1}--\blpage{14}
(\byear{2019})
\doiurl{10.1016/j.finel.2019.02.002}
\end{barticle}
\endbibitem

\bibitem[\protect\citeauthoryear{Gaburro et~al.}{2020}]{Gaburro2020HighChanges}
\begin{botherref}
\oauthor{\bsnm{Gaburro}, \binits{E.}},
\oauthor{\bsnm{Boscheri}, \binits{W.}},
\oauthor{\bsnm{Chiocchetti}, \binits{S.}},
\oauthor{\bsnm{Klingenberg}, \binits{C.}},
\oauthor{\bsnm{Springel}, \binits{V.}},
\oauthor{\bsnm{Dumbser}, \binits{M.}}:
{High order direct Arbitrary-Lagrangian-Eulerian schemes on moving Voronoi
  meshes with topology changes}.
Journal of Computational Physics
\textbf{407}
(2020)
\doiurl{10.1016/j.jcp.2019.109167}
\end{botherref}
\endbibitem

\bibitem[\protect\citeauthoryear{Gaburro}{2021}]{Gaburro2021AChange}
\begin{barticle}
\bauthor{\bsnm{Gaburro}, \binits{E.}}:
\batitle{{A Unified Framework for the Solution of Hyperbolic PDE Systems Using
  High Order Direct Arbitrary-Lagrangian–Eulerian Schemes on Moving
  Unstructured Meshes with Topology Change}}.
\bjtitle{Archives of Computational Methods in Engineering}
\bvolume{28}(\bissue{3}),
\bfpage{1249}--\blpage{1321}
(\byear{2021})
\doiurl{10.1007/s11831-020-09411-7}
\end{barticle}
\endbibitem

\bibitem[\protect\citeauthoryear{Huang and
  Wang}{2020}]{Huang2020AnisotropicMeshes}
\begin{barticle}
\bauthor{\bsnm{Huang}, \binits{W.}},
\bauthor{\bsnm{Wang}, \binits{Y.}}:
\batitle{{Anisotropic mesh quality measures and adaptation for polygonal
  meshes}}.
\bjtitle{Journal of Computational Physics}
\bvolume{410},
\bfpage{109368}
(\byear{2020})
\doiurl{10.1016/j.jcp.2020.109368}
\end{barticle}
\endbibitem

\bibitem[\protect\citeauthoryear{Ciarlet}{2002}]{CiarletElliptic}
\begin{bbook}
\bauthor{\bsnm{Ciarlet}, \binits{P.G.}}:
\bbtitle{{The Finite Element Method for Elliptic Problems}}.
\bpublisher{Society for Industrial and Applied Mathematics}, \blocation{???}
(\byear{2002}).
\doiurl{10.1137/1.9780898719208} .
\burl{http://epubs.siam.org/doi/book/10.1137/1.9780898719208}
\end{bbook}
\endbibitem

\bibitem[\protect\citeauthoryear{Ciarlet and
  Raviart}{1972}]{Ciarlet1972InterpolationMethods}
\begin{barticle}
\bauthor{\bsnm{Ciarlet}, \binits{P.G.}},
\bauthor{\bsnm{Raviart}, \binits{P.A.}}:
\batitle{{Interpolation theory over curved elements, with applications to
  finite element methods}}.
\bjtitle{Computer Methods in Applied Mechanics and Engineering}
\bvolume{1}(\bissue{2}),
\bfpage{217}--\blpage{249}
(\byear{1972})
\doiurl{10.1016/0045-7825(72)90006-0}
\end{barticle}
\endbibitem

\bibitem[\protect\citeauthoryear{Lenoir}{1986}]{Lenoir1986OptimalBoundaries}
\begin{barticle}
\bauthor{\bsnm{Lenoir}, \binits{M.}}:
\batitle{{Optimal Isoparametric Finite Elements and Error Estimates for Domains
  Involving Curved Boundaries}}.
\bjtitle{SIAM Journal on Numerical Analysis}
\bvolume{23}(\bissue{3}),
\bfpage{562}--\blpage{580}
(\byear{1986})
\doiurl{10.1137/0723036}
\end{barticle}
\endbibitem

\bibitem[\protect\citeauthoryear{Donea et~al.}{2004}]{Donea2004}
\begin{bchapter}
\bauthor{\bsnm{Donea}, \binits{J.}},
\bauthor{\bsnm{Huerta}, \binits{A.}},
\bauthor{\bsnm{Ponthot}, \binits{J.-P.}},
\bauthor{\bsnm{Rodr{\'{i}}guez-Ferran}, \binits{A.}}:
\bctitle{{Arbitrary Lagrangian-Eulerian Methods}}.
In: \bbtitle{Encyclopedia of Computational Mechanics},
pp. \bfpage{1}--\blpage{25}.
\bpublisher{John Wiley {\&} Sons, Ltd},
\blocation{Chichester, UK}
(\byear{2004}).
\doiurl{10.1002/0470091355.ecm009}
\end{bchapter}
\endbibitem

\bibitem[\protect\citeauthoryear{Richter}{2017}]{Richter2017}
\begin{bbook}
\bauthor{\bsnm{Richter}, \binits{T.}}:
\bbtitle{{Fluid-structure Interactions}}.
\bsertitle{Lecture Notes in Computational Science and Engineering},
vol. \bseriesno{118}.
\bpublisher{Springer},
\blocation{Cham}
(\byear{2017}).
\doiurl{10.1007/978-3-319-63970-3} .
\burl{http://link.springer.com/10.1007/978-3-319-63970-3}
\end{bbook}
\endbibitem

\bibitem[\protect\citeauthoryear{Beir{\~{a}}o~da Veiga
  et~al.}{2019}]{daVeigaCurvedVEM}
\begin{barticle}
\bauthor{\bsnm{Veiga}, \binits{L.}},
\bauthor{\bsnm{Russo}, \binits{A.}},
\bauthor{\bsnm{Vacca}, \binits{G.}}:
\batitle{{The Virtual Element Method with curved edges}}.
\bjtitle{ESAIM: Mathematical Modelling and Numerical Analysis}
\bvolume{53}(\bissue{2}),
\bfpage{375}--\blpage{404}
(\byear{2019})
\doiurl{10.1051/m2an/2018052}
\end{barticle}
\endbibitem

\bibitem[\protect\citeauthoryear{Dedner and Hodson}{2024}]{Dedner2022ASpaces}
\begin{barticle}
\bauthor{\bsnm{Dedner}, \binits{A.}},
\bauthor{\bsnm{Hodson}, \binits{A.}}:
\batitle{{A framework for implementing general virtual element spaces}}.
\bjtitle{SIAM Journal of Scientific Computing (to be published)}
(\byear{2024})
\doiurl{10.48550/arXiv.2208.08978}
\end{barticle}
\endbibitem

\bibitem[\protect\citeauthoryear{Bastian
  et~al.}{2021}]{Bastian2021TheDevelopments}
\begin{barticle}
\bauthor{\bsnm{Bastian}, \binits{P.}},
\bauthor{\bsnm{Blatt}, \binits{M.}},
\bauthor{\bsnm{Dedner}, \binits{A.}},
\bauthor{\bsnm{Dreier}, \binits{N.-A.}},
\bauthor{\bsnm{Engwer}, \binits{C.}},
\bauthor{\bsnm{Fritze}, \binits{R.}},
\bauthor{\bsnm{Gr{\"{a}}ser}, \binits{C.}},
\bauthor{\bsnm{Gr{\"{u}}ninger}, \binits{C.}},
\bauthor{\bsnm{Kempf}, \binits{D.}},
\bauthor{\bsnm{Kl{\"{o}}fkorn}, \binits{R.}},
\bauthor{\bsnm{Ohlberger}, \binits{M.}},
\bauthor{\bsnm{Sander}, \binits{O.}}:
\batitle{{The Dune framework: Basic concepts and recent developments}}.
\bjtitle{Computers {\&} Mathematics with Applications}
\bvolume{81},
\bfpage{75}--\blpage{112}
(\byear{2021})
\doiurl{10.1016/j.camwa.2020.06.007}
\end{barticle}
\endbibitem

\bibitem[\protect\citeauthoryear{Dedner
  et~al.}{2012}]{Dedner2012Dune-Fem:Computing}
\begin{bchapter}
\bauthor{\bsnm{Dedner}, \binits{A.}},
\bauthor{\bsnm{Kl{\"{o}}fkorn}, \binits{R.}},
\bauthor{\bsnm{Nolte}, \binits{M.}},
\bauthor{\bsnm{Ohlberger}, \binits{M.}}:
\bctitle{{Dune-Fem: A General Purpose Discretization Toolbox for Parallel and
  Adaptive Scientific Computing}}.
In: \bbtitle{Advances in DUNE},
pp. \bfpage{17}--\blpage{31}.
\bpublisher{Springer},
\blocation{Berlin, Heidelberg}
(\byear{2012}).
\doiurl{10.1007/978-3-642-28589-9}
\end{bchapter}
\endbibitem

\bibitem[\protect\citeauthoryear{Frittelli
  et~al.}{2021}]{Frittelli2021Bulk-surfaceDimensions}
\begin{barticle}
\bauthor{\bsnm{Frittelli}, \binits{M.}},
\bauthor{\bsnm{Madzvamuse}, \binits{A.}},
\bauthor{\bsnm{Sgura}, \binits{I.}}:
\batitle{{Bulk-surface virtual element method for systems of PDEs in two-space
  dimensions}}.
\bjtitle{Numerische Mathematik}
\bvolume{147}(\bissue{2}),
\bfpage{305}--\blpage{348}
(\byear{2021})
\doiurl{10.1007/s00211-020-01167-3}
\end{barticle}
\endbibitem

\bibitem[\protect\citeauthoryear{Frittelli and
  Sgura}{2018}]{Frittelli2018VirtualSurfaces}
\begin{barticle}
\bauthor{\bsnm{Frittelli}, \binits{M.}},
\bauthor{\bsnm{Sgura}, \binits{I.}}:
\batitle{{Virtual Element Method for the Laplace-Beltrami equation on
  surfaces}}.
\bjtitle{ESAIM: Mathematical Modelling and Numerical Analysis}
\bvolume{52}(\bissue{3}),
\bfpage{965}--\blpage{993}
(\byear{2018})
\doiurl{10.1051/m2an/2017040}
\end{barticle}
\endbibitem

\bibitem[\protect\citeauthoryear{Bachini
  et~al.}{2021}]{Bachini2021Arbitrary-orderSurfaces}
\begin{barticle}
\bauthor{\bsnm{Bachini}, \binits{E.}},
\bauthor{\bsnm{Manzini}, \binits{G.}},
\bauthor{\bsnm{Putti}, \binits{M.}}:
\batitle{{Arbitrary-order intrinsic virtual element method for elliptic
  equations on surfaces}}.
\bjtitle{Calcolo}
\bvolume{58}(\bissue{3}),
\bfpage{30}
(\byear{2021})
\doiurl{10.1007/s10092-021-00418-5}
\end{barticle}
\endbibitem

\bibitem[\protect\citeauthoryear{Li et~al.}{2022}]{Li2022LocalSurface}
\begin{barticle}
\bauthor{\bsnm{Li}, \binits{J.}},
\bauthor{\bsnm{Feng}, \binits{X.}},
\bauthor{\bsnm{He}, \binits{Y.}}:
\batitle{{Local tangential lifting virtual element method for the
  diffusion–reaction equation on the non-flat Voronoi discretized surface}}.
\bjtitle{Engineering with Computers}
(\byear{2022})
\doiurl{10.1007/s00366-021-01595-1}
\end{barticle}
\endbibitem

\bibitem[\protect\citeauthoryear{Beir{\~{a}}o Da~Veiga
  et~al.}{2016}]{ellipticVEM}
\begin{barticle}
\bauthor{\bsnm{Beir{\~{a}}o Da~Veiga}, \binits{L.}},
\bauthor{\bsnm{Brezzi}, \binits{F.}},
\bauthor{\bsnm{Marini}, \binits{L.D.}},
\bauthor{\bsnm{Russo}, \binits{A.}}:
\batitle{{Virtual Element Method for general second-order elliptic problems on
  polygonal meshes}}.
\bjtitle{Mathematical Models and Methods in Applied Sciences}
\bvolume{26}(\bissue{4}),
\bfpage{729}--\blpage{750}
(\byear{2016})
\doiurl{10.1142/S0218202516500160}
\end{barticle}
\endbibitem

\bibitem[\protect\citeauthoryear{Bonito
  et~al.}{2013a}]{Bonito2013Time-DiscreteStability}
\begin{barticle}
\bauthor{\bsnm{Bonito}, \binits{A.}},
\bauthor{\bsnm{Kyza}, \binits{I.}},
\bauthor{\bsnm{Nochetto}, \binits{R.H.}}:
\batitle{{Time-Discrete Higher-Order ALE Formulations: Stability}}.
\bjtitle{SIAM Journal on Numerical Analysis}
\bvolume{51}(\bissue{1}),
\bfpage{577}--\blpage{604}
(\byear{2013})
\doiurl{10.1137/120862715}
\end{barticle}
\endbibitem

\bibitem[\protect\citeauthoryear{Bonito
  et~al.}{2013b}]{Bonito2013Time-discreteAnalysis}
\begin{barticle}
\bauthor{\bsnm{Bonito}, \binits{A.}},
\bauthor{\bsnm{Kyza}, \binits{I.}},
\bauthor{\bsnm{Nochetto}, \binits{R.H.}}:
\batitle{{Time-discrete higher order ALE formulations: A priori error
  analysis}}.
\bjtitle{Numerische Mathematik}
\bvolume{125}(\bissue{2}),
\bfpage{225}--\blpage{257}
(\byear{2013})
\doiurl{10.1007/s00211-013-0539-3}
\end{barticle}
\endbibitem

\bibitem[\protect\citeauthoryear{Formaggia and
  Nobile}{2004}]{Formaggia2004StabilityALEFEM}
\begin{barticle}
\bauthor{\bsnm{Formaggia}, \binits{L.}},
\bauthor{\bsnm{Nobile}, \binits{F.}}:
\batitle{{Stability analysis of second-order time accurate schemes for
  ALE–FEM}}.
\bjtitle{Computer Methods in Applied Mechanics and Engineering}
\bvolume{193}(\bissue{39-41}),
\bfpage{4097}--\blpage{4116}
(\byear{2004})
\doiurl{10.1016/j.cma.2003.09.028}
\end{barticle}
\endbibitem

\bibitem[\protect\citeauthoryear{Gastaldi}{2001}]{Gastaldi2001AElements}
\begin{barticle}
\bauthor{\bsnm{Gastaldi}, \binits{L.}}:
\batitle{{A priori error estimates for the Arbitrary Lagrangian Eulerian
  formulation with finite elements}}.
\bjtitle{Journal of Numerical Mathematics}
\bvolume{9}(\bissue{2}),
\bfpage{123}--\blpage{156}
(\byear{2001})
\doiurl{10.1515/JNMA.2001.123}
\end{barticle}
\endbibitem

\bibitem[\protect\citeauthoryear{Stein}{1971}]{Stein1971SingularFunctions}
\begin{bbook}
\bauthor{\bsnm{Stein}, \binits{E.M.}}:
\bbtitle{{Singular Integrals and Differentiability Properties of Functions}}.
\bpublisher{Princeton University Press}, \blocation{???}
(\byear{1971}).
\doiurl{10.1515/9781400883882}
\end{bbook}
\endbibitem

\bibitem[\protect\citeauthoryear{Adams and Fournier}{2003}]{Adams2003}
\begin{bbook}
\bauthor{\bsnm{Adams}, \binits{R.A.}},
\bauthor{\bsnm{Fournier}, \binits{J.J.F.}}:
\bbtitle{{Sobolev Spaces}}.
\bpublisher{Elsevier}, \blocation{???}
(\byear{2003})
\end{bbook}
\endbibitem

\bibitem[\protect\citeauthoryear{Ahmad
  et~al.}{2013}]{Ahmad2013EquivalentMethods}
\begin{barticle}
\bauthor{\bsnm{Ahmad}, \binits{B.}},
\bauthor{\bsnm{Alsaedi}, \binits{A.}},
\bauthor{\bsnm{Brezzi}, \binits{F.}},
\bauthor{\bsnm{Marini}, \binits{L.D.}},
\bauthor{\bsnm{Russo}, \binits{A.}}:
\batitle{{Equivalent projectors for virtual element methods}}.
\bjtitle{Computers {\&} Mathematics with Applications}
\bvolume{66}(\bissue{3}),
\bfpage{376}--\blpage{391}
(\byear{2013})
\doiurl{10.1016/j.camwa.2013.05.015}
\end{barticle}
\endbibitem

\bibitem[\protect\citeauthoryear{Beir{\~{a}}o~da Veiga
  et~al.}{2014}]{DaVeiga2014}
\begin{barticle}
\bauthor{\bsnm{Veiga}, \binits{L.}},
\bauthor{\bsnm{Brezzi}, \binits{F.}},
\bauthor{\bsnm{Marini}, \binits{L.D.}},
\bauthor{\bsnm{Russo}, \binits{A.}}:
\batitle{{The Hitchhiker's Guide to the Virtual Element Method}}.
\bjtitle{Mathematical Models and Methods in Applied Sciences}
\bvolume{24}(\bissue{08}),
\bfpage{1541}--\blpage{1573}
(\byear{2014})
\doiurl{10.1142/S021820251440003X}
\end{barticle}
\endbibitem

\bibitem[\protect\citeauthoryear{Sutton}{2017}]{Sutton2017TheMATLAB}
\begin{barticle}
\bauthor{\bsnm{Sutton}, \binits{O.J.}}:
\batitle{{The virtual element method in 50 lines of MATLAB}}.
\bjtitle{Numerical Algorithms}
\bvolume{75}(\bissue{4}),
\bfpage{1141}--\blpage{1159}
(\byear{2017})
\doiurl{10.1007/s11075-016-0235-3}
\end{barticle}
\endbibitem

\bibitem[\protect\citeauthoryear{Di~Pietro and Droniou}{2020}]{HHOBook}
\begin{bbook}
\bauthor{\bsnm{Di~Pietro}, \binits{D.A.}},
\bauthor{\bsnm{Droniou}, \binits{J.}}:
\bbtitle{{The Hybrid High-Order Method for Polytopal Meshes}}.
\bsertitle{MS{\&}A},
vol. \bseriesno{19}.
\bpublisher{Springer},
\blocation{Cham}
(\byear{2020}).
\doiurl{10.1007/978-3-030-37203-3} .
\burl{http://link.springer.com/10.1007/978-3-030-37203-3}
\end{bbook}
\endbibitem

\bibitem[\protect\citeauthoryear{Cangiani et~al.}{2016}]{ConNonConVEM}
\begin{botherref}
\oauthor{\bsnm{Cangiani}, \binits{A.}},
\oauthor{\bsnm{Manzini}, \binits{G.}},
\oauthor{\bsnm{Sutton}, \binits{O.J.}}:
{Conforming and nonconforming virtual element methods for elliptic problems}.
IMA Journal of Numerical Analysis,
0--36
(2016)
\doiurl{10.1093/imanum/drw036}
\end{botherref}
\endbibitem

\bibitem[\protect\citeauthoryear{Dedner and
  Hodson}{2022}]{Dedner2022RobustCoefficients}
\begin{barticle}
\bauthor{\bsnm{Dedner}, \binits{A.}},
\bauthor{\bsnm{Hodson}, \binits{A.}}:
\batitle{{Robust nonconforming virtual element methods for general fourth-order
  problems with varying coefficients}}.
\bjtitle{IMA Journal of Numerical Analysis}
\bvolume{42}(\bissue{2}),
\bfpage{1364}--\blpage{1399}
(\byear{2022})
\doiurl{10.1093/imanum/drab003}
\end{barticle}
\endbibitem

\bibitem[\protect\citeauthoryear{Beir{\~{a}}o~da Veiga
  et~al.}{2020}]{BeiraodaVeiga2020PolynomialEdges}
\begin{barticle}
\bauthor{\bsnm{Veiga}, \binits{L.}},
\bauthor{\bsnm{Brezzi}, \binits{F.}},
\bauthor{\bsnm{Marini}, \binits{L.D.}},
\bauthor{\bsnm{Russo}, \binits{A.}}:
\batitle{{Polynomial preserving virtual elements with curved edges}}.
\bjtitle{Mathematical Models and Methods in Applied Sciences}
\bvolume{30}(\bissue{08}),
\bfpage{1555}--\blpage{1590}
(\byear{2020})
\doiurl{10.1142/S0218202520500311}
\end{barticle}
\endbibitem

\bibitem[\protect\citeauthoryear{Wei{\ss}er}{2019}]{Weier2019AnisotropicAnalysis}
\begin{barticle}
\bauthor{\bsnm{Wei{\ss}er}, \binits{S.}}:
\batitle{{Anisotropic polygonal and polyhedral discretizations in finite
  element analysis}}.
\bjtitle{ESAIM: Mathematical Modelling and Numerical Analysis}
\bvolume{53}(\bissue{2}),
\bfpage{475}--\blpage{501}
(\byear{2019})
\doiurl{10.1051/m2an/2018066}
\end{barticle}
\endbibitem

\bibitem[\protect\citeauthoryear{Antonietti
  et~al.}{2019}]{Antonietti2019TheDiscretizations}
\begin{bchapter}
\bauthor{\bsnm{Antonietti}, \binits{P.F.}},
\bauthor{\bsnm{Berrone}, \binits{S.}},
\bauthor{\bsnm{Verani}, \binits{M.}},
\bauthor{\bsnm{Wei{\ss}er}, \binits{S.}}:
\bctitle{{The Virtual Element Method on Anisotropic Polygonal
  Discretizations}}.
In: \bbtitle{Lecture Notes in Computational Science and Engineering}
vol. \bseriesno{126},
pp. \bfpage{725}--\blpage{733}
(\byear{2019})
\end{bchapter}
\endbibitem

\bibitem[\protect\citeauthoryear{Wells}{2023}]{WellsThesis}
\begin{botherref}
\oauthor{\bsnm{Wells}, \binits{H.}}:
{Moving Mesh Virtual Element Methods}.
PhD thesis,
University of Nottingham
(2023)
\end{botherref}
\endbibitem

\bibitem[\protect\citeauthoryear{Alk{\"{a}}mper
  et~al.}{2016}]{Alkamper2016TheModule}
\begin{barticle}
\bauthor{\bsnm{Alk{\"{a}}mper}, \binits{M.}},
\bauthor{\bsnm{Dedner}, \binits{A.}},
\bauthor{\bsnm{Kl{\"{o}}fkorn}, \binits{R.}},
\bauthor{\bsnm{Nolte}, \binits{M.}}:
\batitle{{The DUNE-ALUGRID Module}}.
\bjtitle{Archive of Numerical Software}
\bvolume{4}(\bissue{1}),
\bfpage{1}--\blpage{28}
(\byear{2016})
\doiurl{10.11588/ans.2016.1.23252}
\end{barticle}
\endbibitem

\bibitem[\protect\citeauthoryear{Dedner et~al.}{2020}]{Dedner2020PythonModule}
\begin{botherref}
\oauthor{\bsnm{Dedner}, \binits{A.}},
\oauthor{\bsnm{Kl{\"{o}}fkorn}, \binits{R.}},
\oauthor{\bsnm{Kl{\"{o}}fkorn}, \binits{R.}}:
{Python Bindings for the DUNE-FEM module}
(2020).
\doiurl{10.5281/zenodo.3706993}
\end{botherref}
\endbibitem

\bibitem[\protect\citeauthoryear{Senechal
  et~al.}{1995}]{Senechal1995SpatialDiagrams}
\begin{barticle}
\bauthor{\bsnm{Senechal}, \binits{M.}},
\bauthor{\bsnm{Okabe}, \binits{A.}},
\bauthor{\bsnm{Boots}, \binits{B.}},
\bauthor{\bsnm{Sugihara}, \binits{K.}}:
\batitle{{Spatial Tessellations: Concepts and Applications of Voronoi
  Diagrams}}.
\bjtitle{The College Mathematics Journal}
\bvolume{26}(\bissue{1}),
\bfpage{79}
(\byear{1995})
\doiurl{10.2307/2687299}
\end{barticle}
\endbibitem

\bibitem[\protect\citeauthoryear{Talischi
  et~al.}{2012}]{Talischi2012PolyMesher:Matlab}
\begin{barticle}
\bauthor{\bsnm{Talischi}, \binits{C.}},
\bauthor{\bsnm{Paulino}, \binits{G.H.}},
\bauthor{\bsnm{Pereira}, \binits{A.}},
\bauthor{\bsnm{Menezes}, \binits{I.F.M.}}:
\batitle{{PolyMesher: a general-purpose mesh generator for polygonal elements
  written in Matlab}}.
\bjtitle{Structural and Multidisciplinary Optimization}
\bvolume{45}(\bissue{3}),
\bfpage{309}--\blpage{328}
(\byear{2012})
\doiurl{10.1007/s00158-011-0706-z}
\end{barticle}
\endbibitem

\end{thebibliography}

\end{document}